\documentclass[a4paper,11pt,french,leqno]{smfart}

\usepackage{hyperref}
\usepackage{amsmath}
\usepackage{amssymb}
\usepackage{amsthm}
\usepackage{amsopn}
\usepackage{amscd}
\usepackage{fullpage}
\usepackage[utf8]{inputenc}
\usepackage{graphics, xcolor}
\usepackage{textcomp} 
\usepackage{verbatim}
\usepackage[francais]{babel}
\usepackage{tikz}

\usepackage{tikz-cd}

\def\commentaire#1{\ifx\cachercommentaires\undefined \textcolor{red}{#1}\else \fi } 

\newcommand{\frI}{\mathfrak{I}}

\newcommand{\frb}{\mathfrak{b}}

\newcommand{\frg}{\mathfrak{g}}

\newcommand{\frk}{\mathfrak{k}}
\newcommand{\frl}{\mathfrak{l}}

\newcommand{\frn}{\mathfrak{n}}

\newcommand{\frqqq}{\mathfrak{q}}

\newcommand{\frt}{\mathfrak{t}}
\newcommand{\fru}{\mathfrak{u}}
\newcommand{\frv}{\mathfrak{v}}

\newcommand{\bbA}{\mathbb{A}}

\newcommand{\bbC}{\mathbb{C}}

\newcommand{\bbN}{\mathbb{N}}

\newcommand{\bbR}{\mathbb{R}}

\newcommand{\bbZ}{\mathbb{Z}}

\newcommand{\caB}{\mathcal{B}}

\newcommand{\caH}{\mathcal{H}}

\newcommand{\caL}{\mathcal{L}}

\newcommand{\caQ}{\mathcal{Q}}
\newcommand{\caR}{\mathcal{R}}

\newcommand{\caT}{\mathcal{T}}

\newcommand{\caV}{\mathcal{V}}

\newcommand{\caX}{\mathcal{X}}

\def\U{\mathbf{U}}

\newcommand{\GL}{\mathbf{GL}}

\newcommand{\Ad}{\mathrm{Ad}}

\newcommand{\Ind}{\mathrm{Ind}}

\newcommand{\tr}{\mathrm{Tr}\, }
\newcommand{\Id}{\mathrm{Id}}

\renewcommand{\Ind}{\mathrm{Ind}}

\newcommand{\SL}{\mathbf{SL}}

\newcommand{\SO}{\mathbf{SO}}
\newcommand{\Sp}{\mathbf{Sp}}
\newcommand{\Or}{\mathbf{O}}

\newcommand{\Tr}{\mathrm{Tr}}

\newcommand{\bil}[2]{\langle  #1,#2 \rangle }

\newcommand{\Trans}{\mathrm{Trans}}
\newcommand{\St}{\mathbf{St}}
\newcommand{\Speh}{\mathbf{Speh}}
\newcommand{\sgn}{\mathbf{sgn}}
\newcommand{\Triv}{\mathbf{Triv}}

\newcommand{\Std}{\mathbf{Std}}

\theoremstyle{plain}
\newtheorem{thm}{Théorème}[section]
\newtheorem{lemme}[thm]{Lemme}
\newtheorem{cor}[thm]{Corollaire}
\newtheorem{prop}[thm]{Proposition}

\theoremstyle{definition}
 \newtheorem{defi}[thm]{Définition}
 
\newtheorem{rmqs}[thm]{Remarques}
\newtheorem{rmq}[thm]{Remarque}
\newtheorem{exemple}[thm]{Exemple}

\def \dem {\noindent \underline{\sl Démonstration}. }

\begin{document}

\numberwithin{equation}{subsection}

\title{Sur les paquets d'Arthur des groupes classiques réels}
  \author{Colette Moeglin}
 \address{CNRS, Institut Mathématique de Jussieu } 
 \email{colette.moeglin@imj-prg.fr}

  \author{David Renard  }
 \address{Centre de Mathématiques
  Laurent Schwartz,  Ecole Polytechnique} 
\email{david.renard@polytechnique.edu}

\date{\today}

\thanks{Le deuxième auteur a bénéficié d'une aide de  l'agence nationale de la recherche 
ANR-13-BS01-0012 FERPLAY}

\begin{abstract} Cet article s'insère dans un projet d'étude des paquets d'Arthur pour les groupes classiques réels commencé dans \cite{AMR}.
Notre but est de donner une description explicite de ces paquets et d'établir la propriété de multiplicité un pour ceux-ci (qui est connue
pour les groupes $p$-adiques et complexes). Le résultat principal est ici une construction de paquets à partir de paquets
unipotents de facteurs de Levi ($c$-Levi) par induction cohomologique. Un outil important de la démonstration est un énoncé de commutativité
entre induction cohomologique et transfert endoscopique spectral.
\medskip
 
\noindent {\bf  Abstract.}\, ---\, 
This article is part of a project  started in \cite{AMR} which consists of investigating  Arthur packets for real classical groups.
 Our goal is to give an explicit
description of these packets and to establish the multiplicity one property (which is known to hold for $p$-adic and complex groups).
The main result  in this paper is a construction
of packets from unipotent packets on  $c$-Levi factors using cohomological induction.   An important tool used in the argument is a statement
of commutativity between cohomological induction and spectral endoscopic transfer.
\end{abstract}

\bigskip
\maketitle

\section{Introduction}

Soit $G$ un groupe classique réel, c'est-à-dire un groupe symplectique ou bien un groupe spécial orthogonal.
Notre but dans cet article, ainsi que dans une série d'articles compagnons (\cite{MR4}, \cite{noteMR}) est de décrire
de manière aussi explicite que possible les représentations irréductibles de $G$ qui sont composantes locales d'une représentation
automorphe de carré intégrable. Plus précisément, soit $F$ un corps de nombres et soit $\mathbf G$ un groupe algébrique réductif 
connexe défini sur $F$.
 Les représentations automorphes de carré intégrables de $\mathbf G$ sont les 
représentations irréductibles du groupe $\mathbf G(\bbA_F)$ ($\bbA_F$ est le groupe des adèles de $F$)
apparaissant comme sous-représentations dans $L^2(\mathbf G(F)\backslash \mathbf G(\bbA_F))$.
Pour les décrire, J. Arthur a introduit des paramètres locaux 
\[\psi_v:\, W'_{F_v}\times \SL_2(\bbC)\longrightarrow {}^LG=\widehat G\rtimes W'_{F_v},\]
 où $v$ décrit l'ensemble des  places de $F$, et $W'_{F_v}$ est le groupe de Deligne-Weil de $F_v$
(si $v$ est une place archimédienne, $W'_{F_v}=W_F$ est le groupe de Weil de $F_v$). Ces paramètres, en plus de leur 
propriétés locales, doivent satisfaire à une propriété de compatibilité globale qui ne nous intéresse pas directement ici; disons simplement
que pour les groupes classiques quasi-déployés et leurs variantes, elle n'est décrite qu'artificiellement dans \cite{Art13}, et que pour l'exprimer
plus généralement, il faudra attendre la description des groupes tannakiens associés aux corps de nombres.

Soit donc $\psi_v$ un paramètre d'Arthur local et notons $A(\psi_v)$ le groupe des composantes connexes du centralisateur de $\psi_v$
dans $\widehat G$ (en fait, il faut en général passer à un revêtement {\sl cf.}  \cite{Art13}, chapitre 9, mais ceci est inutile dans le cas des groupes 
quasi-déployés, et  des groupes classiques même non  quasi-déployés).
Arthur suggère qu'au paramètre $\psi_v$ est attaché une  combinaison linéaire de représentations  irréductibles   
de $\mathbf G(F_v)$ à coefficients dans l'espace des fonctions à valeurs complexes sur le groupe  $A(\psi_v)$, 
fonctions invariantes par conjugaison. On note 
 ${\pi}^A(\psi_v)$ cette combinaison linéaire. Ces objets  ${\pi}^A(\psi_v)$ doivent être  compatible à l'endoscopie. 
 Cela ne suffit pas à les définir dans le cas quasi-déployé. Pour compléter la définition dans ce cas et pour les groupes classiques, 
 Arthur ajoute la compatibilité à l'endoscopie tordue et cela suffit alors.  Limitons nous dans ce qui suit aux cas des groupes classiques.
  Le groupe $A(\psi_v)$ est alors abélien (c'est même un $2$-groupe), les fonctions invariantes par conjugaison sur ce groupe
   sont donc directement des
   combinaisons linéaires à coefficients complexes de caractères de ce groupe  et on cherche donc une combinaison linéaire à 
   coefficients complexes de représentations irréductibles de $\mathbf G(F_v)\times A(\psi_v)$.

Quand le groupe classique est  quasi-déployé, Arthur montre dans \cite{Art13} que 
 ${\pi}^A(\psi_v)$ est une représentation unitaire de $\mathbf G(F_v)\times A(\psi_v)$, c'est-à-dire que la combinaison
  linéaire est à coefficients  entiers positifs (au lieu d'être des nombres complexes), on a donc une représentation semi-simple (par construction) 
  et la propriété supplémentaire est que cette représentation est unitaire.
 Remarquons
que les  ${\pi}^A(\psi_v)$ dépendent du choix des facteurs de transfert géométrique  et  en suivant Kottwitz et Shelstad,
 les facteurs de transfert sont normalisés par des choix de données de Whittaker.
 
 Pour les formes intérieures pures des groupes classiques quasi-déployés sur $F_v$, la compatibilité à l'endoscopie 
 suffit encore à caractériser  ${\pi}^A(\psi_v)$; ce qui n'est pas clair, c'est que 
  ${\pi}^A(\psi_v)$ soit une  représentation unitaire de $\mathbf G(F_v)\times A(\psi_v)$ (\cite{Art13}, Conjecture 9.4.2)
   comme dans le cas quasi-déployé. Ceci est établi pour les paramètres {\sl unipotents} dans 
\cite{noteMR}, en utilisant les mêmes méthodes globales qu'Arthur (c'est-à-dire la stabilisation de la formule des traces (tordue)).
  Les résultats de cet article et de \cite{MR4} montrent alors que cette propriété est en fait vraie pour tous les  paramètres $\psi_v$. 
 
Une question importante qui se pose alors est  celle de la {\sl décomposition} de la représentation ${\pi}^A(\psi_v)$ 
en représentations irréductibles et du  calcul des coefficients. 
Pour les groupes classiques et une place archimédienne complexe, ce problème est résolu dans \cite{MR2}, le cas  
des places réelles avec l'hypothèse supplémentaire que le caractère infinitésimal est entier et régulier est traité en \cite{AMR} 
et dans ces deux cas les coefficients sont un et les représentations irréductibles apparaissant sont simples à décrire. 
Pour le cas des places réelles, cela se fait avec l'induction parabolique ordinaire et l'induction cohomologique. 
Disons tout de suite, avant de détailler, que cet article généralise les méthodes et les résultats de \cite{AMR} et que
 l'on obtient le même type de description très explicite mais sous l'hypothèse que les paramètres des séries discrètes
  intervenant dans $\psi_v$ sont grands les uns par rapport aux autres et grand par rapport à la partie unipotente. 
  Donc par rapport à \cite{AMR}, on montre que de \og rajouter\fg \,  une partie unipotente ne change pas les résultats. 
  Pour avoir le cas général,  on utilise les techniques usuelles de  translation du caractère infinitésimal, c'est l'objet de l'article 
  \cite{MR4},    mais là on perd la description précise puisque la translation traverse des murs. 
  Soyons plus précis sur le contenu de cet article.

 Nous ne considérons donc que le cas des places archimédiennes réelles, 
 et dans la suite de cette introduction, $\mathbf G$ est maintenant un groupe classique défini sur $\bbR$ dont on note 
 $G=\mathbf G(\bbR)$ le groupe de ses points réels. On considère la représentation standard du  $L$-groupe de $G$, 
 $\Std_G:\, {}^LG\rightarrow \GL_N(\bbC)$, et par composition avec cette représentation standard, on peut voir
 un paramètre d'Arthur pour $G$ comme une représentation de $W_\bbR\times \SL_2(\bbC)$. Cette représentation est complètement réductible
 et dans sa décomposition en irréductibles, on va distinguer une partie {\sl de bonne parité}
 et une partie de {\sl mauvaise parité}. De plus, dans la partie de bonne parité, nous distinguons une partie {\sl discrète}, et une partie
  {\sl unipotente}. Une première réduction, qui fera l'objet d'un article ultérieur, est de montrer que la description des représentations
   ${\pi}^A(\psi)$  se ramène au cas où les paramètres sont de bonne parité, par une induction parabolique qui préserve l'irréductibilité
    des composantes. On se limite donc ici au cas de bonne parité.

 
 La réduction établie dans cet article dans le théorème  \ref{thm:Princ} donne une description de la représentation, 
 $\pi^A(\psi)$ cherchée comme une induite cohomologique d'un caractère d'un groupe unitaire convenable et de l'analogue 
 de $\pi^A(\psi)$ associée à un paramètre spécial unipotent (cf. \cite{pourhowe}) d'un groupe de même type que 
 $G$. Le problème est que l'on ne sait pas décrire la décomposition de cette induite cohomologique sans une hypothèse 
 de régularité. Pour les experts, disons que l'induite cohomologique se fait dans le weakly fair range et que l'on sait bien
  décrire le résultat dans le weakly good range. Dans le cas où la condition de régularité que l'on décrit ci-dessous, est 
  satisfaite, on a alors une description exacte en irréductibles de $\pi^A(\psi)$. En particulier cette représentation est sans multiplicité et même  la composante isotypique  d'une représentation irréductible $\pi$ de $\mathbf G(\bbR)$  dans $\pi^A(\psi)$ est réduite à $\pi$.
Expliquons cette condition de régularité.
  
 On suppose que l'on a un paramètre 
 d'Arthur $\psi_G$ pour notre groupe classique $\mathbf G$, qui après composition avec la représentation
 standard, donne une  représentation $\psi$ de $W_\bbR\times \SL_2(\bbC)$ de dimension $N$, et que celle-ci admet une décomposition 
 $\psi=\psi_d\oplus \psi'$. On suppose que $\psi_d$ est une représentation irréductible de la forme $V(0,t)\otimes R[c]$.
 Ici  $R[c]$ est une représentation algébrique irréductible  de $\SL_2(\bbC)$ de dimension $c$ et 
 $V(0,t)$ est une représentation irréductible de dimension  2 de $W_\bbR$, paramètre de Langlands d'une série discrète auto-duale
 de $\GL_2(\bbR)$ notée $\delta(\frac{t}{2}, -\frac{t}{2})$, de caractère infinitésimal $(\frac{t}{2}, -\frac{t}{2})$, $t$ étant un entier supérieur ou égal 
 à $c$. Ainsi $\psi_d$ est le paramètre d'Arthur d'une représentation de Speh, notée  
 $\Pi_{\psi_d}^\GL=\Speh(\delta\left(\frac{t}{2}, -\frac{t}{2}\right), c)$. La condition de bonne parité devant être vérifiée ici est 
 $t+c-1$ pair si le groupe dual de $\mathbf G$ est orthogonal, impair s'il est symplectique.
 A torsion près par le caractère signe de $W_\bbR$, le paramètre $\psi'$ est obtenu à partir d'un paramètre $\psi_{G'}$ d'un 
 groupe classique $\mathbf G'$ de même
 type que $\mathbf G$, de rang $n-c$ si $n$ est le rang de $\mathbf G$, par composition avec la représentation standard de $\mathbf G'$.
 La propriété de régularité exigée est ici que $t$ soit suffisamment grand
par rapport au caractère infinitésimal du paramètre $\psi'$ (la condition exacte est  (\ref{hyplemmegen2})).
Remarquons que l'on ne suppose pas ici une forme particulière pour $\psi'$, même si comme nous l'avons expliqué
ci-dessus, l'application principale sera de partir de paramètres unipotents de bonne parité et d'ajouter un à un 
des blocs discrets $\psi_d$ pour obtenir tous les paramètres de bonne parité réguliers.

 L'idée est de proposer une formule pour la représentation ${\pi}^A(\psi_G)$
 à partir de la représentation  ${\pi}^A(\psi_{G'})$ et des entiers $t$ et $c$, et de montrer qu'elle  satisfait bien les propriétés
 voulues relativement à l'endoscopie et à l'endoscopie tordue qui caractérisent la représentation d'Arthur (théorème \ref{thm:main}).
La formule donnant  ${\pi}^A(\psi_G)$ fait intervenir les foncteurs d'induction cohomologique de Vogan-Zuckerman
à partir de sous-groupe de Levi de $G$ ($c$-Levi) de la forme 
\[L_i\simeq \U(i,c-i)\times G'_i\]
où $\U(i,c-i)$ est un groupe unitaire réel attaché à une signature $(i,c-i)$ et $G'_i$ est une forme intérieure pure
du groupe $\mathbf G'$.
La vérification de l'identité endoscopique tordue vers le groupe tordu $(\GL_N(\bbR), \theta_N)$ se fait essentiellement comme dans 
\cite{AMR}, en reprenant de manière   cruciale certains résultats qui y ont été établis. Dans \cite{AMR}, les identités endoscopiques
ordinaires n'apparaissent pas car elles avaient déjà été établies par Adams et Johnson mais ce n'est 
le cas ici, et leur vérification est le point technique principal de l'article. Elle se fait grâce à un principe de commutation entre transfert endoscopique
et induction cohomologique, qui nous semble présenter un intérêt intrinsèque et constitue à ce titre le deuxième résultat important de l'article. Un tel principe doit être valide en général, mais ce que montre le cas des groupes classiques pour lequel nous 
l'établissons ici est que l'énoncé en est nécessairement subtil. Ceci se devine du fait que si l'induction parabolique ordinaire commute au transfert, 
elle ne commute pas aux foncteurs d'induction cohomologique de Vogan-Zuckerman. Pour avoir des formules de commutation, il faut introduire
certains caractères de torsion sur les sous-groupes de Cartan (voir l'appendice) que l'on retrouve dans la proposition \ref{propquicompte}.
Pour avoir des formules agréables, nous préférons renormaliser
les foncteurs de Vogan-Zuckerman pour y inclure les torsions.  Ces foncteurs renormalisés, dans les cas des groupes orthogonaux, commutent
aux tranferts endoscopiques  (Eq. (\ref{gto})) et pour les groupes symplectiques, il reste une torsion.

 Donnons maintenant rapidement une idée de l'organisation de l'article. La section \ref{notetgen} fixe les notations et donne quelques rappels
 sur les paramètres et la classification de Langlands, puis sur les paramètres d'Arthur et les propriétés  des représentations qui leur sont 
 conjecturalement attachées, en particulier les identités de transfert endoscopiques ordinaires (Eq. (\ref{IdEnd})).
 La section \ref{grcla} introduit les groupes classiques quasi-déployés 
  considérés dans cet article : groupes symplectiques et groupes spéciaux orthogonaux
 quasi-déployés. La section \ref{PaqArtGC} explique une partie des résultats de \cite{Art13}, en particulier
 l'existence de la représentation ${\pi}^A(\psi_G)$ attachée à un paramètre d'Arthur, et sa caractérisation par les identités endoscopiques 
  (\ref{IdEnd}) (pour les données endoscopiques elliptiques du groupe $\mathbf G$) et l'identité endoscopique tordue (\ref{IdEndT}). 
 La section \ref{FIP} est une simple remarque sur l'extension dans \cite{noteMR}
  des résultats d'Arthur aux groupes classiques non nécessairement quasi-déployés.
 La section \ref{motivation} a pour but  
 d'expliquer avec plus de détails qu'au début de   cette introduction les motivations de notre travail et les résultats obtenus.
  La façon dont les paramètres se décomposent   la définition des parties de bonne parité, de mauvaise parité, discrète et unipotente 
  sont donnés en \ref{decpar}.
  On donne ensuite les principaux énoncés de réduction obtenus dans cet article et les articles attenants ainsi que ce qui est connu dans le cas unipotent.
 Dans la section \ref{resgen}, nous revenons à un groupe algébrique réductif connexe défini sur $\bbR$ quelconque, pour rappeler quelques
  résultats de la littérature qui nous seront utiles, en particulier le lien entre donnée de Whittaker et paire fondamentale de type Whittaker, la notion de 
 $c$-sous-groupe de Levi (ce sont les sous-groupes de Levi à partir desquels se fait l'induction cohomologique de Vogan-Zuckerman),
  la paramétrisation des séries discrètes selon Adams et Shelstad, etc.
 Dans la section \ref{StabIC} on énonce un résultat de stabilité de certaines représentations virtuelles obtenue par induction cohomologique
 (proposition \ref{cor:stable2}) essentiel à la formulation même de notre résultat de commutation entre transfert endoscopique et induction 
 cohomologique. Toutefois, pour éviter de trop alourdir le texte, nous ne le démontrerons que pour 
  les groupes classiques un peu plus loin dans le texte, section \ref{demprop}.
Dans la section \ref{ICgrcla} nous revenons au cadre des groupes classiques, et nous spécialisons à ceux-ci les discussions sur les 
séries discrètes et $c$-Levi, ce qui fait apparaître naturellement les (classes d'équivalence de)
formes intérieures pures $\SO(p,q)$ des groupes spéciaux orthogonaux quasi-déployés.
Dans les sections \ref{DEE} et \ref{cLeviH} nous rappelons de manière explicite les données endoscopiques elliptiques des groupes classiques, et 
les $c$-Levi maximaux des groupes endoscopiques correspondants.  Ceci détermine le cadre pour notre énoncé de commutation
 entre  l'induction cohomologique et transfert endoscopique qui est établi dans la section \ref{ICetTrans} (proposition \ref{propquicompte})
 et après renormalisation des foncteurs d'induction cohomologique dans la section \ref{SecNormIC}, dans le théorème \ref{ThmICetTrans}.
 L'extension aux  groupes classiques  non quasi-déployés est discutée dans la section \ref{TICFIP}.
 La section \ref{principal} décrit le résultat  de l'ajout d'un bloc discret de paramètre
 suffisamment grand à un paramètre d'Arthur, que nous avons déjà décrit ci-dessus. 
 Par une recurrence immédiate, mise en place dans la section \ref{PABPR},  on en tire le théorème de réduction \ref{thm:Princ}.
 Enfin l'appendice, tiré en grande partie de \cite{AV},
 introduit les caractères de torsion apparaissant dans notre énoncé de commutation entre induction cohomologique et transfert endoscopique
 et les calcule explicitement pour les groupes classiques.

\section{Notations et généralités}\label{notetgen}

On note $\Gamma=\{1,\sigma\}$ le groupe de Galois de $\bbC/\bbR$.
Si $\mathbf G$ est  un  groupe  algébrique réductif connexe défini sur $\bbR$, on identifie $\mathbf  G$ au groupe de ses points complexes, 
et l'on note  $\sigma_{G}$ l'action de $\sigma$ sur $\mathbf G$. On note $G=\mathbf G(\bbR)$ le groupe des points réels de $\mathbf G$.
 On fixe alors  une involution de Cartan $\tau_G$, qui commute
  avec $\sigma_G$.   On pose $\mathbf K_G=\mathbf G^{\tau_G}$ et $K_G=G^{\tau_G}$, c'est un sous-groupe compact maximal de $G$. 
  Lorsqu'il n'y a pas d'ambigüité, on écrit  simplement 
  $\tau$ et $K$ plutôt que $\tau_G$ et $K_G$.
   Par représentation de $G$, nous entendons un  $(\frg,K_G)$-module de longueur finie, et par représentation virtuelle, un élément du 
   groupe de Grothendieck de la  catégorie des représentations de $G$. On note $[\pi]$ l'image d'une représentation $\pi$ dans ce
    groupe de Grothendieck. On note $\Pi(G)$ l'ensemble des classes d'équivalence de représentations irréductibles de $G$.

Introduisons un invariant important d'un tel groupe $\mathbf G$. On pose 
\begin{equation} \label{qG}
q(G)=\frac{1}{2}(\dim \mathbf G-\dim \mathbf K)-c(G)
\end{equation}
où $c(G)$ est la moitié de la dimension de la partie déployée d'un sous-groupe de Cartan fondamental ({\sl i.e.} maximalement compact)
de $G$. C'est un entier, et si les rangs de $G$ et $K$ sont égaux, cet entier est égal à la moitié de la dimension de l'espace symétrique
$G/K$.

Le signe de Kottwitz attaché au groupe $\mathbf G$ est 
\begin{equation} \label{eG}
e(G)=(-1)^{q(G)-q(G^*)}.
\end{equation} où  $\mathbf G^*$ est une forme intérieure quasi-déployée de $G$.

 \begin{exemple}\label{qUiai} Pour un groupe unitaire réel  $\U(a,b)$ attaché à un espace hermitien de signature $(a,b)$, on 
 a $q(\U(a,b))=ab$ et pour un groupe spécial orthogonal    réel  $\SO(a,b)$ attaché à un espace quadratique  de signature $(a,b)$, on 
 a $q(\SO(a,b))=ab$. 
 \end{exemple}

Lorsque l'on a un  groupe algébrique défini sur $\bbR$ ou $\bbC$, on note par la lettre gothique correspondante
l'algèbre de Lie du groupe de ses points complexes.

\subsection{Paramètres de Langlands}
\label{LangL}

Dans cette section, nous rappelons brièvement la classification de Langlands utilisant le $L$-groupe, 
principalement dans le but d'introduire les notations utilisées dans le reste de l'article.
Nous renvoyons le lecteur à \cite{Bor} pour plus de détails. 

Le groupe de  Weil de  $\bbC$ est $W_\bbC= \bbC^\times$. Le groupe de  Weil  $W_\bbR$ de $\bbR$ 
est une  extension  non  scindée   de $\bbZ/2\bbZ$ par  $\bbC^\times$ :
\begin{equation}\label{WR1}
1 \longrightarrow \bbC^\times \longrightarrow W_\bbR \longrightarrow \bbZ/2\bbZ \longrightarrow 1.
\end{equation}
En identifiant  $\bbC^\times$ avec son image dans  $W_\bbR$ et  $\bbZ/2\bbZ$ avec  $\{\pm 1\}$, on voit que  
  $W_\bbR$ est engendré par  $\bbC^\times$ et un élément  $j$ qui se projette sur  $-1 \in \bbZ/2\bbZ$, 
  avec les relations :
\begin{equation} \label{WR2} j^2=-1\in \bbC^\times, \qquad jzj^{-1}=\bar z, \quad (z\in \bbC^\times).\end{equation}

\bigskip

Soit $\mathbf{G}$ un groupe algébrique réductif connexe défini sur $\bbR$, et soit $\widehat G$ 
le dual de Langlands de $\mathbf{G}$. On fixe un épinglage $\mathbf{spl}_{\hat G}=(\caB, \caT,  \{\caX_\alpha \})$
de $\widehat G$. On construit le $L$-groupe de $\mathbf G$ comme un produit semi-direct
\begin{equation} 
{}^L G=\widehat G\rtimes W_\bbR,
\end{equation}
où l'action de $W_\bbR$ sur $\widehat G$ se factorise par la projection $p_{W_\bbR}: \, W_\bbR\rightarrow \Gamma$, et l'action de 
$\sigma\in \Gamma$ sur $\widehat G$, que l'on note $\sigma_{\widehat G}$ laisse stable $\mathbf{spl}_{\hat G}$.

\begin{defi} \label{DefParLang}
Un paramètre  de Langlands est un morphisme continu : 
\[ \phi: \; W_\bbR \rightarrow  {}^LG  \]
tel que 
\begin{itemize}
\item[a)] $p_{W_\bbR}\circ \phi =\Id_{W_\bbR}$,

\item[b)] la restriction $\phi_{\vert \bbC^ \times }: \;\bbC^ \times  \rightarrow  \widehat G\times \bbC^\times$ 
à pour image des   éléments dont la composante dans  $\widehat G$ est semi-simple.
\end{itemize}

Le groupe $\widehat G$ agit par conjugaison sur l'ensemble des paramètres de Langlands, et l'on note 
$\Phi(G)$ l'ensemble de ces classes de conjugaison. On note $\Phi_{\mathrm{temp}}(G)$  l'ensemble
 des classes de conjugaison de paramètres de Langlands d'image bornée.
\end{defi}

Là encore, nous commettrons fréquemment l'abus de langage consistant à ne pas distinguer
entre un paramètre de Langlands et l'élément de $\Phi(G)$ qu'il définit.  

\medskip 

Le théorème de classification de Langlands est l'existence d'une partition de $\Pi(G)$  : 
\begin{equation} \label{LanglPartition}\Pi(G)=\coprod_{\phi \in \Phi(G)}  \Pi(\phi,G) \end{equation}
où les  $\Pi(\phi,G)$,  appelés $L$-paquets, ou paquets de Langlands,  sont des ensembles finis
(de classes d'équivalence infinitésimales) de représentations irréductibles. Le point essentiel est bien entendu
l'ensemble des propriétés de cette partition. Donnons-en quelques unes.
Tous les éléments d'un paquet ont même caractère
infinitésimal. Pour les représentations tempérées, on a 
\begin{equation} \label{LanglTempPartition}\Pi_{\mathrm{temp}}(G)=\coprod_{\phi \in \Phi_{\mathrm{temp}}(G)}  \Pi(\phi,G). \end{equation}

 Pour tout  $\phi \in \Phi_{\mathrm{temp}}(G)$,  la représentation virtuelle
\begin{equation} \label{TempStableDistr}  \sum_{\pi \in \Pi(\phi,G)} [\pi] \end{equation}
est stable. Rappelons ce que cela signifie.
L'application  qui à une représentation de longueur finie 
associe son caractère passe au groupe de Grothendieck et induit une application linéaire injective du groupe de Grothendieck dans 
l'espace des distributions sur $G$ invariantes par conjugaison. Dire qu'une représentation virtuelle est stable signifie
 que son caractère est une distribution stablement invariante  sur $G$. Pour ce qui concerne la  notion de stabilité, nous renvoyons par exemple à \cite{Shelstad79} et \cite{Boua}.
 Ceci n'est plus vrai pour un paquet non tempéré. Soit $\phi \in \Phi(G)$, non nécessairement tempéré.  
 Les  représentations $\pi\in \Pi(\phi,G)$
sont obtenues dans la classification de Langlands comme uniques sous-représentations irréductibles
de représentations standard $I(\pi)$, obtenues par induction  parabolique à partir  d'un sous-groupe parabolique $P=MN$ de $G$,
et d'un paquet \og essentiellement tempéré\fg \, de $M$.

\begin{defi}\label{pseudopaquets}
Appelons \og pseudo-paquet\fg, 
et notons $\Pi^\sharp(\phi,G)$, l'ensemble
des $I(\pi)$ pour $\pi \in \Pi(\phi,G)$. Définissons la représentation virtuelle
\begin{equation}\label{StablePsP} [\pi(\phi,G)]
= \sum_{I(\pi) \in \Pi^\sharp(\phi,G)} [I(\pi)]. 
\end{equation}
\end{defi}

Le paquet $\Pi(\phi,G)$
est  aussi obtenu de la manière suivante dans la classification de Langlands : il existe un sous-groupe parabolique $\mathbf P=\mathbf M\mathbf N$ 
de  $\mathbf G$ défini sur $\bbR$, et un paramètre de Langlands $\phi_M$ discret qui factorise $\phi$ par l'inclusion
${}^LM\hookrightarrow {}^LG$ tels que les éléments du paquet $\Pi(\phi,G)$ sont les sous-représentations irréductibles des $\Ind_P^G(\pi_M)$,
où les $\pi_M$ parcourent le paquet de séries discrètes ({\sl i.e.} essentiellement de carré intégrable modulo le centre)  $\Pi(\phi_M,M)$.
On a alors 
\begin{equation}\label{StablePsP2} [\pi(\phi,G)]
= \sum_{\pi_M \in \Pi(\phi_M,M)} [\Ind_P^G(\pi_M)]. 
\end{equation}

Il est bien connu (\cite{Shelstad79}, \cite{AdJo},  Lemma 4.3) 
 que  $[\pi(\phi,G) ]$ est une  représentation virtuelle stable de $G$. 

\begin{rmqs}  \label{rmqtemp}
1. Remarquons que lorsque $\phi$ est tempéré,    $[\pi(\phi,G)]$  est aussi égale à (\ref{TempStableDistr}).

--- 2. Selon le contexte, on appelle \og représentations standard\fg \, les représentations $I(\pi)$ de (\ref{StablePsP})
ou les représentations $\Ind_P^G(\pi_M)$ de (\ref{StablePsP2}).

--- 3. Les $[\pi(\phi,G)]$ forment une base du sous-groupe du groupe de Grothendieck constitué des représentations virtuelles stables
\end{rmqs}

\subsection{Paramètres d'Arthur}\label{ArtP}

Les notations sont les m\^emes que dans la section précédente.

\begin{defi}\label{ArtPar}
Un paramètre d'Arthur pour  $G$ est un morphisme de groupes continu
\begin{equation*}
\psi: \, W_\bbR \times \SL_2(\bbC) \longrightarrow {}^ L G
\end{equation*}
tel que 
\begin{itemize}
\item[(i)] la restriction de $\psi$ à  $W_\bbR$ est un paramètre de Langlands tempéré,
\item[(ii)]  la restriction de  $\psi$ à $\SL_2(\bbC)$ est algébrique.
\end{itemize}

\medskip

Le groupe $\widehat G$ agit par conjugaison sur l'ensemble des paramètres d'Arthur, et l'on note 
$\Psi(G)$ l'ensemble de ces classes de conjugaison. On identifie  $\Phi_{\mathrm{temp}}(G)$ à l'ensemble
 des paramètres d'Arthur de restriction triviale à $\SL_2(\bbC)$.
\end{defi}

A tout paramètre d'Arthur  $\psi$, on associe un paramètre de Langlands 
\begin{equation}\label{ArtLang}
\phi_\psi: \, W_\bbR \longrightarrow {}^ L G, \qquad w \mapsto \psi(w,
  \left(  \begin{matrix}  \vert w \vert^{\frac{1}{2}} & 0\\ 
0 &   \vert w \vert^{-\frac{1}{2}}   \end{matrix} \right) , 
\end{equation}
où $w\mapsto \vert w \vert$ est le morphisme de groupe de $W_\bbR$ dans $\bbR^\times_+$ 
défini par $\vert j\vert =1$ et $\vert z\vert =z\bar z$ si $z\in \bbC^\times$.

\subsection{Paquets d'Arthur}\label{pacAr}
Dans \cite{Art84}, \cite{Art89}, J. Arthur conjecture l'existence de paquets $\Pi(\psi, G)$  attachés aux 
paramètres $\psi \in \Psi(G)$,   devant posséder certaines propriétés. Parmi les principales,
 citons le fait que les $\Pi(\psi,G)$ sont finis, constitués de (classes d'équivalence) de  représentations
{\sl  unitaires}, ayant toutes le m\^eme caractère infinitésimal.
 Le paquet d'Arthur $\Pi(\psi,G)$ contient
 le paquet de Langlands $\Pi(\phi_\psi, G)$. 
 Ils doivent  satisfaire les identités de caractères  attendues
dans la théorie de l'endoscopie (standard et tordue), c'est ce qui est appelé le {\sl transfert spectral}.
En revanche, ces  paquets ne sont pas disjoints, et ne sont pas des réunions de $L$-paquets.

\bigskip 

Donnons quelques précisions au sujet des conjectures d'Arthur, en faisant quelques hypothèses simplificatrices
qui seront vérifiées pour les groupes que nous allons étudier. Certains termes (\og forme intérieure pure\fg,\, \og  donnée de Whittaker \fg) 
seront précisés plus loin dans le texte. On suppose que $\mathbf G$ quasi-déployé, ou bien est une forme intérieure pure d'un groupe quasi-déployé.
  
Soit $\psi$ un paramètre d'Arthur pour le groupe $G$. Soit $S_\psi$ le centralisateur de l'image de $\psi$
dans $\widehat G$ et $S_\psi^0$ sa composante connexe neutre.
On pose 
\begin{equation} \label{Apsi} A(\psi)=S_\psi/S_\psi^0.\end{equation}
On suppose que les groupes $A(\psi)$ sont abéliens.
 Notons $\widehat{A(\psi)}$ le groupe des caractères  de $A(\psi)$.

Comme nous l'avons expliqué dans l'introduction, Arthur  conjecture l'existence d'une combinaison linéaire 
à coefficients complexes de représentations irréductibles de $G \times A(\psi)$, notée 
 ${\pi}^A(\psi,G)$ 
devant vérifier un  certain nombre de propriétés, dont nous allons détailler certaines. 
Le paquet  $\Pi(\psi,G)$ est alors l'ensemble des représentations irréductibles de $G$ qui apparaissent dans ${\pi}^A(\psi,G)$.

Pour tout $x \in A(\psi)$, on considère la représentation virtuelle obtenue en évaluant ${\pi}^A(\psi,G)$ en $s_\psi x$, à savoir
\begin{align} \label{distrstable}
[\pi(x,\psi,G)]= {\pi}^A(\psi)(s_\psi x)
\end{align}
où $s_\psi$ est l'image dans $A(\psi)$ de l'élément $\psi((1,-\Id))$, 
$(1,-\Id)\in W_\bbC\times \SL_2(\bbC)$. Ce sont ces 
représentations virtuelles qui apparaissent dans les
 identités de transfert endoscopiques.

Lorsque $x=1$, (\ref{distrstable})  doit être une  représentation virtuelle stable, que 
l'on  note simplement
 \begin{equation}\label{Xpsistable} [\pi(\psi,G)]=  [\pi(1,\psi,G) ]=  {\pi}^A(\psi)(s_\psi ) .
\end{equation}
  Lorsque $\psi=\phi$ est un paramètre tempéré,  l'élément $s_\phi$ est trivial, et  de plus 
il résulte des travaux de Shelstad  (voir \cite{She08II} et \cite{She08III}) que ${\pi}^A(\phi,G)$ est de la forme 
\begin{equation}\label{piAtemp}{\pi}^A(\phi,G)=\bigoplus_{\eta\in \widehat{A(\phi)}} \pi(\phi,\eta,G)\boxtimes \eta
\end{equation}
où les  $\pi(\phi,\eta,G)$ sont des représentations irréductibles tempérées de $G$ ou $0$ et celles qui sont  non nulles
décrivent exactement le paquet $\Pi(\phi,G)$ et sont donc en particulier non isomorphes deux à deux. 
 La notation (\ref{Xpsistable}) pour $\psi=\phi$ paramètre de Langlands tempéré est bien compatible avec    
(\ref{StablePsP}) grâce à la remarque \ref{rmqtemp}.

Soit $x$ un élément dans le centralisateur de $\psi$ dans $\widehat G$ tel que $x^2=1$, et 
notons encore $x$  l'image de $x$ dans $A(\psi)$.
Soit  $(\mathbf H, \caH, x,\xi)$  une donnée endoscopique elliptique de $\mathbf G$, qui factorise $\psi$. 
 Nous renvoyons à \cite{LS87} et \cite{KS99} pour la définition d'une donnée endoscopique. 
Nous considérons ici le transfert spectral, pour lequel, outre ces deux références, on peut aussi consulter \cite{She08III}.
 Pour simplifier, nous supposons que 
$\caH$ est isomorphe au $L$-groupe de $\mathbf H$, et que cet isomorphisme composé avec $\xi:\, \caH \rightarrow {}^LG$
donne un plongement de $L$-groupes $\iota_{H,G}:\, {}^LH \rightarrow {}^LG$. 

Le paramètre $\psi$ se factorise en 
$\psi=\iota_{H,G}\circ \psi_x$, ce qui définit le paramètre d'Arthur $\psi_x$ pour $\mathbf H$ dans la formule ci-dessous.
Le transfert est normalisé par le choix d'une donnée de Whittaker pour $G$.
L'identité de transfert  endoscopique s'écrit: 
\begin{equation}\label{IdEnd}
\Trans_H^G( [\pi(\psi_x,H)])= e(G)\;  \pi(x,\psi,G). \end{equation}

Remarquons que dans le cas d'un paramètre de Langlands tempéré $\phi$, l'identité endoscopique devient d'après 
Shelstad (\cite{She08II} et \cite{She08III}), avec les notations
introduites en (\ref{piAtemp}) 
\begin{equation}\label{IdEndTemp}
\Trans_H^G( [\pi(\phi_x,H)])=  e(G)\sum_{\eta \in \widehat{A(\phi)} }  \eta(x) \;   [\pi(\phi,\eta,G)] .
\end{equation}

\begin{rmq}\label{transGG*}
Remarquons que lorsque on prend $x=1$, le groupe endoscopique associé est la forme intérieure quasi-déployée de $\mathbf G$, notons-la
$\mathbf G^*$, et l'identité endoscopique devient, pour tout paramètre d'Arthur $\psi$
\[ \Trans_{G^*}^G( [\pi(\psi,G^*)])= e(G) [\pi(\psi, G)]. \] 
\end{rmq}

Lorsque $G$ est un groupe classique quasi-déployé, J. Arthur donne dans \cite{Art13} une définition des 
 ${\pi}^A(\psi,G)$. Il montre que ce sont des 
représentations unitaires de $G\times A(\psi)$, qui sont caractérisées par  les identités endoscopiques
(\ref{IdEnd}), pour toutes les données endoscopiques elliptiques de $G$, ainsi qu'une identité d'endoscopie tordue, où
 $\mathbf G$  apparaît comme un groupe endoscopique pour un  groupe tordu $(\GL_N,\theta_N)$, $\theta_N$ étant l'automorphisme 
extérieur de $\GL_N$.
Nous en dirons plus sur  tout ceci dans la section  \ref{PaqArtGC}.

\section{Les groupes classiques et leurs paquets d'Arthur}

\subsection{Les  groupes classiques}\label{grcla}
  
 Les \og groupes classiques \fg \;  quasi-déployés que nous considérons sont ceux qui apparaissent
 dans les données endoscopiques elliptiques simples des groupes tordus $(\GL_N,\theta_N)$ selon la terminologie d'Arthur {\sl cf.}
 \cite{Art13}, \S I.2 (voir  le paragraphe suivant).   Le plus commode est encore d'en faire la liste, et de fixer quelques notations pour
  pouvoir s'y référer  facilement. Pour $n\in \bbN^\times$, on considère les groupes de rang $n$ suivants: 
 
 \medskip
 
 \begin{itemize}
 \item[\bf A.] Le groupe symplectique $\Sp_{2n}$.  
 
 C'est un groupe déployé. Son dual de Langlands est $\SO(2n+1,\bbC)$ et son $L$-groupe est le 
 produit direct $\SO(2n+1,\bbC)\times W_\bbR$. 
 
  \item[\bf B.] Le groupe spécial orthogonal impair $\SO_{2n+1}$. 
  
   C'est un groupe déployé. Son dual de Langlands 
  est $\Sp(2n,\bbC)$ et son $L$-groupe est le   produit direct $\Sp(2n,\bbC)\times W_\bbR$.
   
    \item[\bf C.]  Le groupe spécial orthogonal pair déployé  $\SO_{2n}^d$, $n\geq 2$.
  
   C'est un groupe déployé. Son dual de Langlands  est $\SO(2n,\bbC)$ et son $L$-groupe est le 
 produit direct $\SO(2n,\bbC)\times W_\bbR$.

    \item[\bf D.]  Le groupe spécial orthogonal pair quasi déployé,  non déployé $\SO_{2n}^{qd}$, $n\geq 1$. 
  
   C'est un groupe quasi-déployé. Son dual de Langlands  est $\SO(2n,\bbC)$ et son $L$-groupe est le 
 produit semi-direct $\SO(2n,\bbC)\rtimes W_\bbR$.
   
   \end{itemize}

\bigskip 

Pour chacun de ces groupes, on dispose d'une représentation naturelle du $L$-groupe dans un $\widehat \GL_N$ : 
\begin{equation}\label{DefStdG} \Std_G: \; {}^LG \longrightarrow  \GL_N(\bbC) \end{equation}
Dans le cas $\mathbf A$, elle est donnée par l'inclusion de $\SO(2n+1,\bbC)$ dans $\GL_{2n+1}(\bbC)$, dans le cas 
$\mathbf B$, par l'inclusion de $\Sp(2n,\bbC)$ dans $\GL_{2n}(\bbC)$, dans le cas $\mathbf C$,
par  l'inclusion de $\SO(2n,\bbC)$ dans $\GL_{2n}(\bbC)$. Le cas $\mathbf D$, est plus délicat
car le groupe étant non déployé, le $L$-groupe est un produit semi-direct non trivial
$\SO(2n,\bbC)\rtimes W_\bbR$. Mais l'action de $W_\bbR$ sur $\SO(2n,\bbC)$ est donnée par l'action d'un élément de 
$\Or(2n,\bbC)$, de sorte que l'on a un morphisme
\[{}^L\SO_{2n}^{qd}=  \SO(2n,\bbC)\rtimes W_\bbR  \rightarrow  \Or(2n,\bbC), \]
et la composition avec  l'inclusion de $\Or(2n,\bbC)$ dans $\GL_{2n}(\bbC)$ nous donne la représentation voulue.
On a donc $N=2n$ dans les cas $\mathbf B$, $\mathbf C$, $\mathbf D$, et  $N=2n+1$ dans le cas $\mathbf A$. 

\medskip

Dans la section  \ref{cLeviH} ci-dessous, nous allons considérer des $c$-sous-groupes de  Levi des groupes classiques et de 
leurs données endoscopiques elliptiques, qui vont avoir des groupes unitaires comme facteurs.
A la liste ci-dessus, on rajoute donc  : 

\medskip 

 \begin{itemize}
 \item[\bf U.] Le groupe unitaire $\U_{n}$.   C'est un groupe quasi-déployé. Son dual de Langlands  est $\GL(n,\bbC)$ et son $L$-groupe est le 
 produit semi-direct $\GL_n(\bbC)\rtimes W_\bbR$.

\end{itemize}

\subsection{Paquets d'Arthur des groupes classiques}\label{PaqArtGC}

Comme nous l'avons expliqué à la fin de la section \ref{pacAr}, lorsque $G$ est un groupe classique quasi-déployé
 et $\psi_G$ est un paramètre, J. Arthur établit dans \cite{Art13}  l'existence de
  ${\pi}^A(\psi_G,G)$ vérifiant les propriétés voulues et montre que c'est une représentation unitaire 
  de $G\times A(\psi)$.  On la 
 décompose en somme de produits tensoriels extérieurs de représentations de ces deux groupes:
\begin{equation}\label{GA}
 {\pi}^A(\psi,G)=\bigoplus_{\eta\in \widehat{A(\psi)}}  \pi(\psi,\eta,G)\boxtimes \eta=\bigoplus_{\pi \in \Pi(\psi, G)}   \pi \boxtimes   \rho_{\pi}.
 \end{equation}
Les résultats  d'Arthur nous disent en particulier  que la représentation virtuelle $[\pi(\psi_G,G)]$ définie en (\ref{Xpsistable})
 est stable. Il  montre que les identités endoscopiques
  (\ref{IdEnd}), pour toutes les données endoscopiques elliptiques de $G$, ainsi qu'une égalité endoscopique
  tordue que nous allons décrire ci-dessous   caractérisent ${\pi}^A(\psi_G,G)$.
  Remarquons simplement que l'identité endoscopique   (\ref{IdEnd}) se réécrit en tenant compte de (\ref{GA}) sous la forme
\begin{equation}\label{IdEnd2}
\Trans_H^G( [\pi(\psi_H,H)])=  e(G)\sum_{\eta \in \widehat{A(\psi)} }  \eta(s_\psi x) \;   [\pi(\psi,\eta,G)] 
=e(G)  \sum_{\pi \in \Pi(\psi,G) }
  \Tr(\rho_{\pi}(s_\psi x)) \;   [\pi]. \end{equation}

  Décrivons maintenant l'identité endoscopique tordue.
  On  note $\tau : \, g\mapsto {}^tg^{-1}$   l'involution de Cartan de $\GL_N(\bbR)$. 
 Soit $J_N\in \GL_N(\bbR)$ la matrice antidiagonale  {\tiny $ \begin{pmatrix}
 .&.&   &&& 1  \\   
.&. &&&-1&  \\
.&&& 1 && \\
&&.&&&\\
&.&&&&&\\
(-1)^{N-1}&.&.&&.&.\end{pmatrix}$ }
  
  On note $\theta_N :\, \mathbf{GL}_N\rightarrow  \mathbf{GL}_N$ l'automophisme défini par   $g\mapsto J_N\, ({}^tg^{-1})J_N^{-1}$.
On définit le produit semi-direct $\mathbf{G}_N^+=\mathbf{GL}_N \rtimes \langle \theta_N\rangle$.
 C'est un groupe algébrique réductif non connexe. L'automorphisme $\theta_N$ étant d'ordre 2,
 ce groupe compte deux composantes connexes, et l'on note $\widetilde{\mathbf{G}}_N=\mathbf{GL}_N \rtimes  \theta_N$ celle qui
  ne contient pas l'élément neutre.  On note  $\widetilde G_N$ l'ensemble des points réels de $\widetilde{\mathbf{G}}_N$.
On obtient ainsi un espace tordu au sens de Labesse \cite{LaWa}. 

La donnée de $(\mathbf G,\Std_G)$ est celle d'une {\sl  donnée endoscopique tordue}
elliptique pour $\widetilde G_N$. Nous renvoyons le lecteur à \cite{KS99} et \cite{Art13}
pour tout ce qui concerne la théorie de l'endoscopie tordue. Rappelons seulement
que dans une telle situation, Kottwitz et Shelstad définissent un  facteur de transfert, permettant de définir une 
application $\mathrm{Trans_{geo}}$ (\og {\sl transfert géométrique} \fg) entre l'espace des intégrales orbitales  
sur $\widetilde G_N$ et  l'espace des intégrales  intégrales orbitales stables sur $G$.
 Ceci est démontré par D. Shelstad dans \cite{Shel}.
Ce facteur de transfert n'est défini {\sl a priori } qu'à une constante multiplicative près, mais le choix de la donnée
 de Whittaker sur $\GL_N(\bbR)$ ({\sl cf.} \cite{AMR}  section 5.2) permet de fixer cette constante 
 (\cite{KS99}, section 5.3).

\medskip

Par dualité,  ce transfert d'intégrales orbitales définit un transfert spectral entre représentations virtuelles
stables de $G$ et représentations virtuelles de $\widetilde G_N$, noté
\begin{equation}\label{Trans}
\mathrm{Trans}_G^{\widetilde G_N}\end{equation}

Soit $\psi_G \in \Psi(G)$ un paramètre d'Arthur pour le groupe $G$. Posons 
$\psi=\Std_G \circ \psi_G$. C'est un paramètre d'Arthur (auto-dual) pour le groupe 
$G_N$. Soit $\Pi_\psi^\GL$ la représentation irréductible autoduale de $\GL_N(\bbR)$ associée à ce paramètre
({\sl cf.} \cite{AMR}, \S 3.1).

L'identité endoscopique tordue est alors 
\begin{equation}\label{IdEndT}
\mathrm{Trans_G^{\widetilde G_N}} \left( [\pi(\psi_G, G)]\right)=\Tr_{\theta_N}(\Pi_\psi^\GL),
\end{equation}  
  le membre de droite étant la trace tordue, normalisée par la donnée de Whittaker ({\sl cf.}
 \cite{AMR},   Définition 5.4).

Lorsque le paramètre $\psi_G$ est trivial sur le facteur $\SL_2(\bbC)$, 
le paquet $\Pi(\psi_G,G)$ est un paquet de Langlands tempéré, la représentation  virtuelle $[\pi(\psi_G, G)]$
est la somme des représentations du paquet,  et l'identité (\ref{IdEndT}) est alors
démontrée par P. Mezo \cite{mezo}, à un facteur multiplicatif près. Nous avons démontré dans \cite{AMR} 
qu'en fait  l'identité (\ref{IdEndT}) est valide, autrement dit que le facteur multiplicatif restant à déterminer
dans Mezo est $1$.
 Le résultat de Mezo 
 et le fait que  le transfert endoscopique commute à l'induction entraîne le résultat suivant
 pour les pseudo-paquets ({\sl cf.} Définition \ref{pseudopaquets}). 
\begin{prop}\label{TransPsPaq}Soient $G$ un groupe classique quasi-déployé comme ci-dessus et  $\phi_G$ un paramètre de Langlands pour $G$.
Posons $\phi=\Std_G \circ \phi_G$. On a alors 
 \[ \mathrm{Trans}( [\pi(\phi_G, G) ])=\Tr_{\theta_N}(\Pi_\phi^\GL). \]
  \end{prop}
  
 \subsection{Formes intérieures pures}\label{FIP}
Dans  \cite{noteMR} les résultats d'Arthur sont étendus aux formes intérieures pures 
 des groupes spéciaux  orthogonaux quasi-déployés. Soit $\mathbf G$ l'une de ces formes intérieures pures, que l'on peut voir comme un 
 groupe $\SO(p,q)$ (voir sections \ref{seriesdiscretescla} et \ref{subsec51}). 
Pour les groupes symplectiques, il n'y a pas de forme intérieure pure de $\Sp_{2n}$ qui ne soit pas équivalente à 
  $\Sp_{2n}$.

Les ${\pi}^A(\psi,G)$ sont alors simplement  caractérisées par les identités endoscopiques (\ref{IdEnd}).
Remarquons que lorsque $\mathbf G$ est la forme quasi-déployée, l'identité endoscopique avec $x=1$
est tautologique. Ce n'est pas le cas lorsque l'on a une forme intérieure non quasi-déployée, 
l'identité endoscopique avec $x=1$ (le groupe endoscopique $\mathbf H$
est  alors la forme intérieure quasi-déployée de $\mathbf G$) donne une relation supplémentaire qui remplace celle donnée
par le transfert endoscopique tordu et suffit à compléter la caractérisation de   ${\pi}^A(\psi,G)$.

Dans ce cadre, les ${\pi}^A(\psi,G)$ ne sont pas supposées {\sl a priori} être des représentations unitaires
de $G\times A(\psi)$, mais on montre dans {\sl loc. cit.} que tel est bien le cas pour les paramètres unipotents. Les 
résultats de réduction énoncés dans la section suivante permettent au final d'étendre ce résultat à tous les paramètres. Notons ce résultat.
  
  \begin{prop}\label{pipsiFIP}
  Soit $\mathbf G$ une forme intérieure pure d'un groupe orthogonal quasi-déployé, et soit $\psi$
  un paramètre d'Arthur pour $\mathbf G$. Alors ${\pi}^A(\psi,G)$ est une représentation unitaire de $G\times A(\psi)$.
  \end{prop}

  \section{Enoncé des résultats de réduction}\label{motivation}
 Comme nous l'avons expliqué dans les deux  sections
 précédentes, les  ${\pi}^A(\psi_G,G)$ attachées  aux paramètres $\psi_G$
 des groupes classiques
 sont caractérisées par des identités endoscopiques, mais l'on aimerait en savoir plus sur  celles-ci.
 En particulier, on aimerait déterminer la décomposition des ${\pi}^A(\psi_G,G)$  et montrer que les coefficients sont 1  
 (propriété de multiplicité un).
 Pour cela, nous obtenons des résultats de réduction que nous allons maintenant décrire.

   \subsection{Décomposition des paramètres}\label{decpar}

Soit $\mathbf G$ l'un des groupes classiques de la section \ref{grcla}.  

\begin{defi} On définit la {\sl bonne parité  pour le groupe} $G$ comme étant $1\mod 2$ si le groupe dual
est symplectique (cas  {\bf B}), et $0 \mod 2$  si le groupe dual est orthogonal (cas {\bf A}, {\bf C}, \bf{D}). 
\end{defi}

Soit $\psi_G$ un paramètre d'Arthur de $\mathbf G$.  Posons 
\[\psi=\Std_G\circ \psi_G:\, W_\bbR\times \SL_2(\bbC)\longrightarrow \GL_N(\bbC). \]

Cette représentation de $W_\bbR\times \SL_2(\bbC)$  est semi-simple.
Rappelons que les représentations irréductibles de $W_\bbR$ sont de deux types ({\sl cf. }\cite{AMR}, \S 3.1 ) : les représentations 
$W(s,\epsilon)$, $s\in \bbC$, $\epsilon\in \{0,1\}$, de dimension $1$ (c'est le paramètre de Langlands de la représentation 
 $\gamma(s,\epsilon) : \, x\mapsto \sgn(x)\, \vert x\vert ^s$ de $\GL_1(\bbR)$), et les représentations $V(s,t)$, $s\in \bbC$, $t\in \bbN^\times$
(c'est le paramètre de Langlands de la représentation essentiellement de carré intégrable modulo le centre 
 $\delta(s,t)$ de caractère infinitésimal $\left(    \frac{s+t}{2},    \frac{s-t}{2} \right)$ de    $\GL_2(\bbR)$).
On note $R[a]$ un représentant de l'unique classe d'équivalence de représentations algébriques de dimension $a$ de $\SL_2(\bbC)$.

Avec ces notations, on peut donc écrire la 
 décomposition en irréductibles de $\psi$ sous  la forme générale (avec des ensembles d'indices disjoints) :
\[ \psi= \bigoplus_i \left(  W(s_i,\epsilon_i) \boxtimes R[a_i] \oplus   W(-s_i,\epsilon_i) \boxtimes R[a_i] 
\right) 
\oplus  \bigoplus_j \left(  V(s_j,t_j) \boxtimes R[a_j] \oplus   V(-s_j,t_j) \boxtimes R[a_j] \right) \]
\[\oplus \bigoplus_k \left(  W(0,\epsilon_k) \boxtimes R[a_k] \oplus   W(0,\epsilon_k) \boxtimes R[a_k] \right) 
 \oplus  \bigoplus_\ell \left(  V(0,t_\ell) \boxtimes R[a_\ell] \oplus   V(0,t_\ell) \boxtimes R[a_\ell] \right) \]
\[\oplus \bigoplus_m   W(0,\epsilon_m) \boxtimes R[a_m] \
 \oplus  \bigoplus_r   V(0,t_r) \boxtimes R[a_r]\]

Dans la première somme,  les $s_i$ sont dans $i\bbR^\times$ et les $\epsilon_i$ dans $\{0,1\}$.
Dans la deuxième  somme, les $s_j$ sont dans $i\bbR^\times$ et les $t_j$ dans $\bbN^\times$.
Dans la troisième  somme, les $\epsilon_k$ sont  dans $\{0,1\}$ et la parité des $a_k-1$ est mauvaise.
Dans la quatrième  somme,  les $t_\ell$ sont dans $\bbN^\times$ et la parité de $t_\ell+a_\ell-1$ est mauvaise.
Dans la cinquième  somme, les $\epsilon_m$ sont  dans $\{0,1\}$ et la parité des $a_m-1$ est bonne.
Dans la sixième  somme,  les $t_r$ sont dans $\bbN^\times$ et la parité de $t_r+a_r-1$ est bonne.
\medskip 

En effet, le paramètre est $\theta_N$-stable. La contribution des  facteurs irréductibles non $\theta_N$-stable 
apparaît alors sous la forme de la première et deuxième somme.
Ce qu'il faut remarquer ensuite, c'est que la multiplicité d'un facteur $W(0,\epsilon) \boxtimes R[a] $
(resp. $V(0,t) \boxtimes R[a]$)  dans $\psi$ est   paire lorsque la parité de $a-1$ (resp.  $t+a-1$)
est mauvaise. La contribution de ces facteurs apparaît alors sous la forme de la troisième et quatrième somme.

On note $\psi_{mp}$ le paramètre formé avec les quatre premières sommes, et $\psi_{bp}$ celui formé avec les deux dernières.
On a alors $\psi=\psi_{mp}\oplus \psi_{bp}$. On décompose aussi $\psi_{bp}$ en $\psi_{bp,u}$ (cinquième somme)
et $\psi_{bp,disc}$ (sixième somme).

On remarque que $\psi_{mp}$ s'écrit sous la forme $\psi_{mp}=\rho\oplus \rho^*$ ($\rho^ *$ 
est la contrégrédiente de $\rho$). Une telle décomposition n'est pas unique.

\subsection{Réduction aux paramètres de bonne parité} \label{reducbp} 
 Les résultats de cette section seront démontrés   ailleurs. 
Soit $\mathbf G$ l'un des groupes classiques de la section \ref{grcla}.  Soit $\psi_G$ un paramètre d'Arthur pour $\mathbf G$, et 
comme précédemment posons $\psi=\Std_G\circ \psi_G$.
Considérons   une décomposition de $\psi$ de la forme :
\begin{equation}\label{psipsi}
       \psi=\rho\oplus \rho^*\oplus \psi' 
 \end{equation}
 où, dans $\rho$, il n'apparaît que des facteurs de mauvaise parité. 
 Le paramètre $\psi'$ se factorise par le $L$-groupe d'un groupe classique quasi-déployé 
$\mathbf G'$ de même type que $\mathbf G$. Soit $\psi_{G'}$ le paramètre d'Arthur pour le groupe $\mathbf G'$ tel que $\psi'=\Std_{G'}
\circ\psi_{G'}$. Notons $N_\rho$ la dimension de la représentation $\rho$
de $W_\bbR \times \SL_2(\bbC)$, et soit $ \Pi^{\GL}_\rho$
la représentation de $\GL_{N_\rho}(\bbR)$ de paramètre d'Arthur $\rho$ ({\sl cf. } \cite{AMR},   \S 3.1).
 Le groupe $G$ admet un sous-groupe de Levi maximal standard $M$ isomorphe
à $\GL_{N_\rho}(\bbR)\times G'$, et ceci fournit une injection
\begin{equation}\label{injGG}
\iota: \;  \GL_{N_\rho}(\bbR) \times {}^LG'\hookrightarrow {}^LG
\end{equation}
de sorte que $\psi_G=\iota\circ (\psi_{G'},\rho)$.

\begin{rmq} Les groupes $A(\psi_G)$ et $A(\psi_{G'})$ sont naturellement isomorphes.        
\end{rmq}

\begin{prop} \label{reducPaqArt}
Soit $\eta \in \widehat{A(\psi_G)}$ et soient $\pi(\psi_G,\eta,G)$ et $\pi(\psi_{G'},\eta, G')$ les représentations semi-simples de $G$ et $G'$
respectivement attachées par Arthur ({\sl cf.} (\ref{GA}), où pour $\pi(\psi_{G'},\eta,G')$ on tient compte de la remarque 
ci-dessus). On a alors
\begin{equation}\label{IndpsiG}
\pi (\psi_G,\eta,G) =   \Ind_P^G  \left( \Pi^{\GL}_\rho\otimes \pi(\psi_{G'},\eta, G')\right),
  \end{equation}
\end{prop}
où $P$ est un sous-groupe parabolique standard  maximal de $G$ de facteur de Levi $M$.

\begin{thm}\label{thm:reducPaqArt} Soit $\pi' \in \Pi(\psi',G')$. Alors $  \Ind_P^G  \left( \Pi^{\GL}_\rho\otimes \pi '\right) $ est irréductible. 
\end{thm}

\begin{cor}\label{mult1etcie} Si les $ \pi(\psi_{G'},\eta,G')$ ont la propriété de multiplicité $1$ et sont disjointes, il en est de même des 
 $\pi(\psi_{G},\eta,G)$. 
Ainsi la propriété de multiplicité 1 pour ${\pi}^A(\psi_G,G)$ découle de celle pour  ${\pi}^A(\psi_{G'},G')$. 
D'autre part, si la décomposition en irréductibles de    ${\pi}^A(\psi_{G'},G')$ est connue, alors elle  l'est en principe
pour ${\pi}^A(\psi_G,G)$.
\end{cor}

Evidemment, la formulation de la seconde partie est un peu vague, connaître la décomposition en irréductible de 
 ${\pi}^A(\psi_{G'},G')$, cela veut dire savoir donner les paramètres des composantes irréductibles des $\pi(\psi_{G'},\eta,G')$
 dans une classification connue. Il faut savoir ensuite ce que donne une induction parabolique irréductible
 dans cette classification, ce qui est le cas pour les classifications usuelles (\cite{KV}, Chapter 11).

Si $\mathbf G$ est maintenant un groupe spécial orthogonal  non quasi-déployé de rang $n$, donc une forme intérieure pure
de l'un des groupes de la section \ref{grcla}, cas {\bf B,C,D}, 
 il faut légèrement adapter la formulation de ces résultats 
 (voir la section \ref{FIP}).
  Si $N_\rho> \inf(p,q)$, on a 
${\pi}^A(\psi_G,G)=0$. Si $N_\rho\leq \inf(p,q)$,
$G$ admet un sous-groupe de Levi maximal standard $M$ isomorphe
à $\GL_{N_\rho}(\bbR)\times G'$, où $G'$ est un groupe spécial orthogonal  de rang $n-N_\rho$.  
Si l'on fait   l'hypothèse que  ${\pi}^A(\psi_{G'},G')$ est une représentation unitaire de $ G' \times  A(\psi_{G'})$, alors les 
$\pi(\psi_{G'},\eta, G')$ sont bien définis pour tout $\eta\in  \widehat{ A(\psi_{G'} ) } $ et la formule (\ref{IndpsiG})
définit    une représentation unitaire $\pi (\psi_G,\eta,G)$ de $G$ pour tout $\eta \in \widehat{A(\psi_{G})} \simeq \widehat{A(\psi_{G'})}$.
On en déduit une représentation unitaire   de $G \times A(\psi_{G})$ : 
\[{\pi}^A(\psi_G,G)=\bigoplus_{\eta\in A(\psi_{G})}  \pi (\psi_G,\eta,G)\boxtimes \eta,  \]
et l'on vérifie qu'elle satisfait bien aux identités endoscopique voulues.
La proposition, le théorème et le corollaire ci-dessus sont  donc  toujours valides, mais ils sont  précédés de 
\begin{prop} Si ${\pi}^A(\psi_{G'},G')$ est une représentation unitaire  de $G' \times  A(\psi_{G'})$, alors 
${\pi}^A(\psi_G,G)$ est une représentation unitaire de $G\times  A(\psi_{G})$.
\end{prop}

\medskip
Ainsi, la description des paquets d'Arthur est réduite au cas des paramètres de bonne parité.

\subsection{Cas des paramètres unipotents de bonne parité}\label{parunip}

Dans cette section, $\mathbf G$ désigne un groupe classique non nécessairement  quasi-déployé, c'est-à-dire 
l'un des groupes quasi-déployés  de la section \ref{grcla}
ou bien une forme intérieure pure d'un groupe spécial orthogonal. On considère le cas des paramètres unipotents de bonne parité.
On suppose donc que $\psi_G$  est tel que $\psi=\psi_{bp,u}$.

\begin{thm} Soit $\psi_G$ un paramètre d'Arthur pour $G$ tel que   $\Std_G\circ \psi_G=\psi=\psi_{bp,u}$.
Alors  ${\pi}^A(\psi_G,G)$ est une représentation unitaire de $G\times A(\psi_G)$ qui   vérifie la propriété de multiplicité un.
\end{thm}

Pour les groupes classiques quasi-déployés (les groupes classiques sur $\bbC$ sont aussi traités), on sait d'après Arthur 
que  ${\pi}^A(\psi_G,G)$ est une représentation unitaire de $G\times A(\psi_G)$, ceci a déjà été 
dit.  La propriété de multiplicité un est établie par  
le premier auteur dans  \cite{pourhowe}.   Le cas des  groupes classiques non quasi-déployés est réglé en \cite{noteMR}.

Le problème de classification est lui aussi essentiellement résolu dans  \cite{pourhowe}, en termes de correspondance de Howe.
D'autre part, on a le résultat suivant.
\begin{thm}  Soit $\pi \in \Pi(\psi_G,G)$. Alors $\pi$ est une representation faiblement unipotente au sens de 
\cite{KV}, Chapter XII.
\end{thm}
Ceci est  démontré dans \cite{MR4}.

\subsection{Réduction du cas des paramètres de bonne parité réguliers au paramètres unipotents}
Soit $G$ un groupe classique,  et soit $\psi_G$ un paramètre d'Arthur pour $G$.
Supposons que $\psi_G$ soit de bonne parité, on décompose donc $\psi=\Std_G\circ \psi_G$ en 
\[  \psi=  \left(\bigoplus_{r=1,\ldots R} V(0,t_r)\boxtimes R[a_r]\right) \oplus   \left(  \bigoplus_{m=1,\ldots M}
 W(0,\epsilon_m) \boxtimes R[a_m] \right) \]
où $\psi_{bp,u}: =  \bigoplus_{m=1,\ldots M} W(0,\epsilon_m)  $ est unipotent de bonne parité.
On suppose que $\psi_G$ vérifie la {\sl condition de régularité} suivante :
\begin{equation}\label{conditionregul0}  \forall r=1,\ldots ,R-1, \;
t_r-a_r+1>t_{r+1}+a_{r+1}-1 \; \text{ et }  t_{R} -a_{R}+1> \max_{m=1,\ldots M } \{a_m-1\}
\end{equation}

Le théorème \ref{thm:Princ} donne une formule pour la représentation unitaire  ${\pi}^A(\psi_G,G)$ utilisant l'induction 
cohomologique et les représentations ${\pi}^A(\psi_{G'},G_i')$ attachée à la partie unipotente du paramètre.
Ceci résout le problème de classification et  celui de multiplicité un pour de tels paramètres.

\subsection{Réduction du cas des paramètres de bonne parité au cas régulier}
On se place dans le même contexte que la section précédente, mais l'on ne suppose plus
la condition de régularité (\ref{conditionregul0}).
Dans \cite{MR4}, nous obtenons une formule pour  ${\pi}^A(\psi_G,G)$ en introduisant un paramètre 
 $\psi_G^{reg}$ de $G$ vérifiant
\[  \psi^{reg}=  \Std_G\circ \psi_G^{reg}=\left(\bigoplus_{r=1,\ldots R} V(0,t'_r)\boxtimes R[a_r]\right) \oplus   \left(  \bigoplus_{m=1,\ldots M}
 W(0,\epsilon_m) \boxtimes R[a_m] \right) \]
 et la condition de régularité  (\ref{conditionregul0}).
D'après le paragraphe précédent, les problèmes de classification et de multiplicité un sont résolus pour ${\pi}^A(\psi_G^{reg},G)$.
La  formule pour ${\pi}^A(\psi_G,G)$ est obtenue en appliquant les techniques usuelles de translation du caractère infinitésimal
(foncteur de translation de Zuckerman). Malheureusement, il apparaît des réductibilités difficiles à contrôler, ce qui nous  empêche
pour le moment de conclure sur les problèmes de décomposition en irréductibles  et de multiplicité un.

\section{Quelques résultats généraux}\label{resgen}
 
\subsection{Donnée de Whittaker et paires de Borel fondamentales}
   
  Si $\mathbf G$ est  quasi-déployé, 
on fixe un épinglage $\mathbf{spl}_{\mathbf G}=(\mathbf B_d, \mathbf T_d, \{X_\alpha \})$ stable par   $\sigma_{G}$.
On note $\mathbf N_d$ le radical unipotent de $\mathbf B_d$, 
$N_d=\mathbf N_d(\bbR)=\mathbf N_d^{\sigma_ G}$ le groupe de ses points réels, et l'on fixe un caractère
$\chi$ de $N_d$ de sorte que $(N_d, \chi)$ est une donnée de Whittaker de $G$. Le rôle d'une telle donnée dans 
la théorie d'Arthur-Langlands-Shelstad est  de distinguer
certaines représentations, celle qui admettent  un modèle de Whittaker, et  de normaliser les facteurs de transfert
endoscopiques  de \cite{LS87}  et \cite{KS99}.

Nous allons plutôt utiliser une donnée équivalente, celle d'une {\sl paire de Borel  fondamentale de type Whittaker}
({\sl cf. }\cite{She15}, \S 2.2, 2.3, 2.4).
Rappelons qu'une paire de Borel  $(\mathbf B,\mathbf T)$  de $\mathbf G$ est un couple constitué 
d'un tore maximal $\mathbf T$ de $\mathbf G$ et d'un sous-groupe de Borel $\mathbf B$.

\begin{defi}
Une paire de Borel 
 $(\mathbf B_*,\mathbf T_*)$ de $\mathbf G$ est dite fondamentale si les conditions suivantes sont réalisées : 

$(i)$ $\mathbf T_*$ est stable sous $\sigma_{G}$ et 
$T_*=\mathbf T_*(\bbR)=\mathbf T_*^{\sigma_{G}}$ est un sous-groupe de Cartan fondamental
({i.e.} maximalement compact) de $G$, 

$(ii)$ $\mathbf B_*$ est un sous-groupe de Borel contenant $\mathbf T_*$ tel que l'ensemble des 
racines de $\mathbf T_*$ dans $\mathbf B_*$ soit stable par $-\sigma$ (cette condition est automatique lorsque 
$\mathbf T_*$ est anisotrope car toute les racines sont imaginaires. Comme $\mathbf T_*$ est fondamental, il n'y a pas de racines
réelles, et c'est donc une condition sur les racines complexes).

 Une paire de Borel fondamentale 
 $(\mathbf B_*,\mathbf T_*)$ de $\mathbf G$ est dite de type Whittaker  
si de plus  la condition suivante est réalisée  :

$(iii)$ les racines simples
imaginaires de  $\mathbf T_*$ dans $\mathbf B_*$ sont imaginaires non compactes. 
\end{defi}

Dans \cite{AV}, les auteurs  utilisent  la terminologie \og large\fg\, pour la  propriété $(iii)$.
Le fait que $\mathbf G$ possède une paire  de Borel fondamentale de type Whittaker  est équivalent
au fait que $\mathbf G$ soit quasi-déployé ({\sl cf.} \cite{AV}, Prop. 6.24). Le choix d'une telle paire fondamentale est équivalent 
au choix d'une donnée de Whittaker $(N_d,\chi)$, la correspondance entre les deux
étant réalisée de la manière suivante : une paire de Borel  fondamentale de type Whittaker
 $(\mathbf B_*,\mathbf T_*)$ de $\mathbf G$ détermine, en se fixant un caractère infinitésimal entier, 
 une représentation générique de la série fondamentale. Cette représentation admet donc un modèle de Whittaker, et 
 détermine donc à conjugaison près une donnée de Whittaker $(N_d,\chi)$. 
 Cette construction induit une bijection entre classes de conjugaison sous $G$ de données de Whittaker 
 et classes de conjugaison de  paire de Borel  fondamentales de type Whittaker.

On fixe donc dans la suite une paire de Borel fondamentale  $(\mathbf B_*,\mathbf T_*)$ de $\mathbf G$, et si $\mathbf G$ est
quasi-déployé, on suppose que cette paire est de type Whittaker et  est 
compatible avec le modèle de Whittaker choisi. On peut supposer, ce que l'on fera, que $\mathbf T_*$ est $\tau$-stable.
La condition $(ii)$ est alors équivalente au fait que $\mathbf B_*$ est aussi $\tau$-stable.

\subsection{$c$-Levi}\label{cLevi}

Soit ${}^d Q= {}^d L{}^d V$ un sous-groupe parabolique standard de $\widehat G$.
Ceci signifie que $\caB\subset {}^dQ$ et $\caT\subset {}^dL$.
Les racines simples de $\caT$ dans $\caB$ sont soit  dans ${}^dL$, soit dans  ${}^dV$, et l'on fait l'hypothèse 
que l'action de $\sigma_{\widehat G}$ sur les racines simples préserve cette partition.

Soit $(\mathbf B, \mathbf T)$ une paire de Borel fondamentale dans $\mathbf G$. 
L'identification entre racines simples  de $\mathbf T$ dans $\mathbf B$ et coracines simples 
de $\caT$ dans $\caB$ permet de définir un sous-groupe parabolique
$\mathbf Q=\mathbf L\mathbf V$ de $\mathbf G$ contenant $\mathbf B$. 

\begin{lemme} Le groupe $\mathbf L$ est stable sous $\sigma_{G}$, autrement dit, 
c'est un sous-groupe de $\mathbf G$ défini sur $\bbR$. Si de plus, $\mathbf T$ est $\tau$-stable, il en est de même de 
$\mathbf Q$ et $\mathbf L$.
\end{lemme}
\dem Le sous-groupe $\mathbf T$ est stable par $\sigma_{G}$ par définition, il suffit donc de démontrer
que l'ensemble des racines de $\mathbf T$ dans $\mathbf L$ est stable sous l'action de $\sigma_{G} $.
Soit $n_0$ un représentant de l'élément le plus long du groupe de Weyl.
 Alors  $ \Ad(n_0) \circ \sigma_{G}$ préserve $(\mathbf B, \mathbf T)$, et quitte à modifier $n_0$
 en le multipliant à gauche par un élément de $T$, on peut trouver un épinglage 
 $\mathbf{spl}=(\mathbf B, \mathbf T, \{X_\alpha\})$, préservé par $ \Ad(n_0) \circ \sigma_{G}$. 
 On a alors par hypothèse $\mathbf L$ préservé par  $ \Ad(n_0) \circ \sigma_{G}$, et donc, comme $\Ad(n_0)$ envoie une racine $\alpha$ 
 sur $-\alpha$ et que l'ensemble des racines de $\mathbf T$ dans $\mathbf L$ est stable par $\alpha\mapsto -\alpha$, on voit que 
  $\mathbf L$ est préservé par  $ \sigma_{G}$. Si $\mathbf T$ est $\tau$-stable, l'action de $\tau$ sur les racines
  coïncide avec celle de $-\sigma_G$, et la seconde assertion s'en déduit.
\qed

\bigskip 
Suivant Shelstad \cite{She15}, 
nous appelons un sous-groupe $\mathbf L$ de $\mathbf G$ défini sur $\bbR$ obtenu de la sorte un $c$-Levi de $\mathbf G$.

\bigskip

Considérons l'ensemble $\Sigma_{{}^dQ}^G$ des classes de conjugaison sous $G=\mathbf G(\bbR)$
de paires $(\mathbf Q,\mathbf L)$ ainsi obtenues. Pour deux paires   $(\mathbf Q,\mathbf L)$ et $(\mathbf Q',\mathbf L')$, on 
voit, puisque leur $L$-groupes sont les mêmes, que $\mathbf L$ et $\mathbf L'$ sont des formes intérieures.

Pour chaque élément dans    $\Sigma_{{}^dQ}^G$, on peut trouver un représentant ($\mathbf Q,  \mathbf L)$ où  
$\mathbf Q$ et  $ \mathbf L$ sont  $\tau_{G}$-stables, ce que l'on supposera toujours.
Si $\mathbf G$ est quasi-déployé, on trouve facilement un représentant  quasi-déployé. En effet, 
partons de la paire fondamentale de type Whittaker $(\mathbf  B_*, \mathbf T_*)$
de $\mathbf G$. On construit  le sous-groupe parabolique $\mathbf L_*=  \mathbf L_*\mathbf V_*$ contenant $\mathbf B_*$
comme ci-dessus. Alors nous avons vu que  $\mathbf L_*$
est $\sigma_{G}$-stable et   $(\mathbf  B_{*,L}= \mathbf B_*\cap \mathbf L_*, \mathbf T_*)$
est une  paire fondamentale de type Whittaker de $\mathbf L_*$.
 Le groupe $\mathbf L_*$ possédant une 
paire fondamentale de type Whittaker, il est quasi-déployé, et nous avons vu que ceci détermine aussi
une classe de conjugaison de donnée de Whittaker pour $\mathbf L_*$.

Dans tous les cas, on a choisi  une paire fondamentale  $(\mathbf  B_*, \mathbf T_*)$, et l'on  obtient  grâce à ce choix 
une bijection entre  $\Sigma_{{}^dQ}^G$ et 
\[   W(T_*,G)\backslash W(\mathbf T_*,\mathbf G)^\tau/W( \mathbf T_*,\mathbf L_*)\]
(\cite{AdJo}, Section 10). Lorsque $\mathbf T_*$ est compact, c'est-à-dire lorsque le rang de $\mathbf G$ est celui 
de $\mathbf K_{G}$, on a $ W(\mathbf T_*,\mathbf G)^\tau= W(\mathbf T_*,\mathbf G)$.

\medskip

Construisons un $L$-groupe en constatant que $\mathbf{spl}_{{}^dL}=(\caB_L:= \caB\cap {}^dL, \caT, \{\caX_\alpha\})$ est un épinglage de 
${}^dL$ préservé par $\sigma_{\widehat G}$.
On construit alors ${}^LL={}^dL\rtimes W_\bbR$ en étendant l'action de $\Gamma$ sur ${}^dL$ via la projection de $W_\bbR$ sur $\Gamma$.

Définissons le  plongement de $L$-groupes suivant 
\begin{equation}\label{iotaLG}
\iota_{L,G}: \; {}^LL\longrightarrow {}^LG  
\end{equation}
qui prolonge l'inclusion de ${}^dL=\widehat L$ dans $\widehat G$.
Il suffit donc de donner la valeur de $\iota_{L,G}$ sur ${1}\rtimes W_\bbR\in {}^LL$.
Pour $z\in \bbC^\times$, on pose 
\[\iota_{L,G}   \prod_{\alpha\in R(\caT,{}^dV)}    \check \alpha  \left(\frac{z}{\vert z \vert}\right)\rtimes z  
=\frac{\check \rho_{{}^dV}(z)}{\check \rho_{{}^dV}(\bar z)}\rtimes z.\]
Ici $\check \rho_{{}^dV}$ est la demi-somme des racines de $\caT$ dans ${}^dV$. 
On note $w_G$ l'élément le plus long de $W_G$ et $w_L$ celui de $W_L$. On note 
\[n: W_G=W(\widehat G,\caT)\rtimes W_\bbR\longrightarrow N({}^LG,\caT)=N(\widehat G,\caT)\rtimes W_\bbR\]
la section ensembliste définie dans \cite{LS87}, \S 2.1. 
Il est démontré dans \cite{Tai} que 
$$n(w_Gw_L\rtimes j)^2=2\check \rho_{{}^dV}(-1).$$ 
On pose 
$\iota_{L,G}(j)= n(w_Lw_G\rtimes j)$,  et on obtient un morphisme bien défini. On a aussi $n(w\rtimes j)=n(w)\rtimes j$
pour tout $w\in W_G$.

\begin{lemme}\label{Transitivite} Supposons que nous soyons dans la situation ci-dessus, et que l'on fixe pour $\mathbf L$ et 
$\caL$ les paires de Borel $(\mathbf B_L=\mathbf B\cap \mathbf L,\mathbf T)$ et $(\caB_\caL=\caB\cap \caL,\caT)$.
Supposons que  $\caQ_1=\caL_1\caV_1$ et $\mathbf Q_1=\mathbf L_1\mathbf V_1$ soient des sous-groupes
paraboliques de $\caL$ et $\mathbf L$ obtenus comme ci-dessus, mais relativement à $\caL$ et $\mathbf L$.
On a alors des plongements de $L$-groupes
\[  \iota_{L,G}: \; {}^LL\longrightarrow {}^LG   \quad \text{et}\quad \iota_{L_1,L}: \; {}^LL_1\longrightarrow {}^LL . \]
On a aussi des sous-groupes paraboliques $\caQ'=\caL_1\caV_1\caV$ et $\mathbf Q'=\mathbf L_1\mathbf V_1\mathbf V$
de $\widehat G$ et $\mathbf G$ respectivement, et un plongement de $L$-groupes
\[   \iota_{L_1,G}: \; {}^LL_1\longrightarrow {}^LG \]
 On a alors $\iota_{L_1,G}=\iota_{L,G}\circ\iota_{L_1,L}$.
\end{lemme}

\dem Tout se joue sur l'élement $j\in W_\bbR$. Or 
\[n(w_{L_1}w_L)n(w_Lw_G)=t(w_{L_1}w_L, w_Lw_G) n(w_{L_1}w_G)\]
où $w\mapsto t(w)$ est un certain cocycle, donné par la formule de \cite{LS87}, Lemma 2.2. On vérifie facilement que 
$t(w_{L_1}w_L, w_Lw_G)=1$.\qed

\subsection{Séries discrètes} \label{seriesdiscretes}

Supposons maintenant que le sous-groupe de Cartan $\mathbf T_*$ de la paire fondamentale 
fixée $(\mathbf B_*,\mathbf T_*)$ soit anisotrope. Alors le groupe $G$ admet des séries discrètes
qui sont des représentations de carré intégrable.

 On applique les construction de la section \ref{cLevi} au cas 
 ou   ${}^d Q= {}^d L{}^d V=\caB$, ${}^dL=\caT$. Dans ce cas le sous-groupe parabolique 
$\mathbf Q=\mathbf L\mathbf V$ de $\mathbf G$ associé à 
une   paire de Borel fondamentale $(\mathbf B, \mathbf T)$  est le sous-groupe de Borel 
$\mathbf B$ lui-même (et donc $\mathbf L=\mathbf T$). A conjugaison près dans $G$, on peut supposer que 
$\mathbf T=\mathbf T_*$, et bien sûr la paire fondamentale  
$(\mathbf B_*, \mathbf T_*)$  fixée au départ nous donne un élément distingué 
d'un système de représentant de $\Sigma_\caB^G$. Le choix de ce représentant distingué
fournit une bijection 
\[  \Sigma_\caB^G\simeq W(G,T_*)\backslash W(\mathbf G,\mathbf T_*).\]

Soit $\phi_G:\, W_\bbR\rightarrow {}^LG$ un paramètre de Langlands discret.
Il se factorise par le $L$-plongement 
\[ \iota_{T_*,G}:\,  {}^L T_*\longrightarrow  {}^LG \]
définit en (\ref{iotaLG}), c'est-à-dire que $\phi_G=\iota_{ T_*,G}\circ \phi_{ T_*}$, pour  certains paramètres de Langlands
$\phi_{T_*}:\, W_\bbR\rightarrow {}^L T_*$. La classification de Langlands pour les tores (\cite{Bor}, \cite{LangTor})
associe à  un tel $\phi_{ T_*}$ un caractère $\xi=\xi_{\phi_{T_*}}$ de $ T_*$.

Si $(\mathbf B, \mathbf T_*)$  est une paire de Borel fondamentale comme ci-dessus,
il existe un paramètre de Langlands  $\phi_{T_*}(\mathbf B)$ comme ci-dessus telle que la différentielle du caractère $\xi_{\phi_{T_*}(\mathbf B)}$ 
soit   dominante relativement  à $\mathbf B$. Notons 
$\caR^d_{\frb, T_*,G}$ le foncteur d'induction cohomologique de Vogan-Zuckerman, en degré
\begin{equation} d=d(G,T)=\frac{1}{2}(\dim \mathbf G-\dim \mathbf T_*)-(q(G)-q( T_*))\end{equation}
   Alors 
 $\caR^d_{\frb, T_*,G}(\xi_{\phi_{ T_*} (\mathbf B)})$ est une série discrète du paquet $\Pi(\phi_G,G)$, et cette construction fournit  une bijection
entre  $\Sigma_\caB$ et $\Pi(\phi_G,G)$.
La condition de dominance de $\xi_{\phi_{ T_*} (\mathbf B)}$ relativement à $\mathbf B$ entraîne que l'on est dans le \og good range \fg \, 
pour l'induction cohomologique (\cite{KV}, Definition 0.49).
Faisons maintenant le lien avec une autre paramétrisation des séries discrètes. Nous suivons les idées de \cite{Adsd}.

\medskip

 Si  $\phi_G:\, W_\bbR\rightarrow {}^LG$ est un paramètre de Langlands, grâce aux travaux de Shelstad
 \cite{She82}, \cite{She08III}, les éléments du pseudo-paquet
 $\Pi^\sharp(\phi_G,G)$ vont être paramétrés par certains caractères de  $A(\phi_G)$, 
  le groupe des composantes connexes du groupe $\mathrm{Centr}(\widehat G,\phi_G)$, 
 le centralisateur dans $\widehat G$ de $\phi_G$. Expliquons ceci pour un paramètre discret.
 Dans ce cas,   $\Pi^\sharp(\phi_G,G)=\Pi(\phi_G,G)$
est un paquet de séries discrètes et $A(\phi_G)= \mathrm{Centr}(\widehat G,\phi_G)$ 
est un 2-groupe fini.

Dans la suite, nous supposons de plus $\mathbf G$ quasi-déployé et une paire fondamentale de type Whittaker
$(\mathbf B_*,\mathbf T_*)$ fixée.

Notons $ T_*[2]$  le sous-groupe des élements d'ordre 2 de $ T_*$.
Nous allons expliquer comment  paramètrer les éléments de $\Pi(\phi_G, G)$  par certains éléments de $ T_*[2]$.
Le lien avec les caractères du groupe  $A(\phi_G)$ provient d'un pairing parfait canonique ({\sl cf.} \cite{Adsd})
\begin{equation}\label{pairingAdams}
 T_*[2]\times A(\phi_G)\rightarrow \bbC^\times, 
\end{equation}
qui identifie $ T_*[2]$ à $\widehat{A(\phi_G)}$.

Pour  tout $t\in T_*[2]$, $\Ad(t)$ définit un automorphisme intérieur 
involutif de $\mathbf G$, qui  commute avec $\sigma_{G}$.  Posons 
$\sigma_t=\sigma_G\circ \Ad(t)$ : c'est un automorphisme antiholomorphe involutif de $\mathbf G$, c'est-à-dire une forme réelle
intérieure de $\mathbf G$. C'est même une forme intérieure \og pure \fg \, au sens de Kottwitz ({\sl cf.} \cite{Kal}) (les formes intérieures
\og non pures \fg \, étant obtenues de la même manière en partant d'éléments $t\in T$ d'ordre fini et dont le carré est dans le centre de $\mathbf G$).
Notons $\mathbf G_t$ le groupe défini sur $\bbR$ dont le groupe des points complexes est $\mathbf G$ et dont
la conjugaison complexe est $\sigma_t$. Remarquons que $T_*$ est un sous-groupe de Cartan commun à tous les $G_t$.

Notons $\caR^{d_t}_{\frb_*, T_*,G_t}$ le foncteur d'induction cohomologique de Vogan-Zuckerman, de la catégorie 
des représentations du tore $ T_*$ vers la catégorie des répresentations de $G_t$, en degré 
\begin{equation} d_t=d(G_t,T_*)=\frac{1}{2}(\dim \mathbf G-\dim \mathbf T_*)-(q(G_t)-q( T_*)). \end{equation}

On obtient donc pour chaque $t\in T_*[2]$ une représentation $\pi(\phi_G,t):=\caR^{d_t}_{\frb_*, T_*,G_t}(\xi_{\phi_{ T_*} (\mathbf B_*)})$
de la série discrète de $G_t$. L'ensemble $\{\pi (\phi_G,t)=\caR^{d_t}_{\frb_*, T_*,G_t}(\xi_{\phi_{ T_*} (\mathbf B_*)}), \, t\in T_*[2]\}$
constitue le  \og super-paquet \fg\, associé à $\phi_G$,   selon une terminologie consacrée.
 Remarquons la différence
avec ce que l'on a fait ci-dessus, où l'on restait dans le groupe $G$ mais on faisait varier le sous-groupe de Borel $\mathbf B$.
Ici, on fixe le sous-groupe de Borel $\mathbf B_*$, et l'on fait varier les formes réelles $G_t$.

On peut regrouper les formes réelles $\mathbf G_t$ par classe d'équivalence. Deux formes réelles $G_t$ et $G_{t'}$ sont équivalentes
s'il existe $g\in \mathbf G$ tel que $\sigma_t=\Ad(g)\circ \sigma_{t'}\circ \Ad(g)^{-1}$.
Ceci définit une relation d'équivalence   $t\sim_G t' $ sur $ T_*[2]$, et l'on note $[t]_G$ la classe d'équivalence 
d'un élément $t\in T_*[2]$. La classe d'équivalence de $1\in T_*[2]$ est donc constitué des éléments $t$
tels que $\mathbf G_t$ est équivalente à $\mathbf G$. On peut donc par conjugaison intérieure
identifier les représentations d'un tel  $G_t$ avec les représentations de $G$. On a alors 
\begin{equation}\label{Pisd}
\Pi(\phi_G,G)= \{\pi(\phi_G,t)=\caR^{d_t}_{\frb_*, T_*,G_t}(\xi_{\phi_{ T_*} (\mathbf B^*)}), \, t\in T_*[2], \; [t]_G=[1]_G \}.\end{equation}
On peut faire la même chose en fixant un système de représentants des classes d'équivalence de formes réelles, et en se ramenant
aux représentations de ces représentants. On obtient ainsi les paquets de séries discrètes des formes intérieures pures de $\mathbf G$.
Soit $t_0 \in T_*[2]$, les représentations d'un groupe $G_t$ avec $[t]_G=[t_0]_G$ s'identifient par conjugaison intérieure aux représentations
de  $G_{t_0}$, et 
\begin{equation}\label{PisdFI}
\Pi(\phi_G,G_{t_0})= \{\pi(\phi_G,t)=\caR^{d_t}_{\frb_*, T_*,G_t}(\xi_{\phi_{ T_*} (\mathbf B^*)}), \, t\in T_*[2], \; [t]_G=[t_0]_G \}.\end{equation}

\begin{rmq}[Limites de séries discrètes]\label{lsd} Ce que l'on a fait ci-dessus pour un paramètre 
$\phi_G$ discret s'adapte facilement au cas d'un paramètre de limites de séries discrètes. Les inductions
cohomologiques sont alors dans le weakly good range, et la différence avec le cas discret est que certaines des 
représentations $\pi(\phi_G,t):=\caR^{d_t}_{\frb_*, T_*,G_t}(\xi_{\phi_{ T_*} (\mathbf B_*)})$ peuvent être nulle.
La paramétrisation des paquets est alors comme en (\ref{Pisd}) et (\ref{PisdFI}), à ceci près qu'il faut rajouter la condition
  $\pi(\phi_G,t)\neq 0$.
  Le groupe $A(\phi_G)$ est un quotient du groupe $A(\phi_G^{reg})$ où $\phi_G^{reg}$ est un paramètre discret, et 
  $\widehat{A(\phi_G^{reg})}\simeq T_*[2]$ par le pairing (\ref{pairingAdams}), et donc $\widehat{A(\phi_G)}$ s'identifie bien
  à un sous-groupe de $T_*[2]$, dont les éléments sont les $t$ tels que  $\pi(\phi_G,t)\neq 0$.
\end{rmq}

\subsection{Stabilité et induction cohomologique} \label{StabIC}
Revenons au contexte de la section \ref{cLevi}. On suppose donc que l'on a un  sous-groupe parabolique standard 
${}^d Q= {}^d L{}^d V$  de $\widehat G$,  compatible avec   l'action de $\sigma_{\widehat G}$.

On suppose  que $\mathbf G$ est quasi-déployé et que l'on a fixé une paire fondamentale de type Whittaker
$(\mathbf B_*,\mathbf T_*)$. 
Elle  nous fourmit un  sous-groupe parabolique $\mathbf Q_*=\mathbf L_*\mathbf V_*$, et 
$\mathbf L_*$ est stable par  $\sigma_G$ et $\tau_G$.

Supposons aussi que $\mathbf T_*$ soit anisotrope.
Soit $\phi_G:\, W_\bbR \rightarrow {}^LG$ un paramètre de Langlands discret, que nous factorisons en 
$\phi_G=\iota_{T_*, G}\circ \phi_{T_*}$ avec $\phi_{T_*}$ choisi de sorte que le caractère $\xi_{\phi_{T_*}}$ qui lui correspond soit dominant pour 
$\mathbf B_*$.
Notons $\phi_L=\iota_{T_*,L_*}\circ \phi_{T_*}$, où $\iota_{T_*,L_*}$ est le $L$-plongement défini par $(\caB_L,\caT)$.
C'est un paramètre de Langlands discret pour $\mathbf L_*$.

Nous avons vu dans la section précédente une paramétrisation des séries discrètes de paramètre $\phi_G$  par les éléments de
$ T_*[2]$. Remarquons que cette paramétrisation s'applique à $\mathbf G$, mais aussi à $\mathbf L_*$.
En particulier, pour tout $t\in T_*[2]$, $\mathbf L_*$ est stable sous $\sigma_t$, et l'on a donc une forme intérieure 
$\mathbf L_{t}$ de $\mathbf L_*$.

Notons $\caR^{d(G_t,L_t)}_{\frqqq_1,L_{t},G_t}$ le foncteur d'induction cohomologique de Vogan-Zuckerman, de la catégorie 
des représentations de $L_{t}$ vers la catégorie des représentations de $G_t$, en degré
\begin{equation} d(G_t,L_t) =\frac{1}{2}(\dim \mathbf G-\dim \mathbf L_*)-(q(G_t)-q(L_{t}))  .  \end{equation}
 Remarquons que par le lemme  \ref{Transitivite}  de transitivité  des $L$-plongements, on a aussi $\phi_L=\iota_{T_*, L}\circ \phi_{T_*}$.
Par transitivité de l'induction cohomologique ({\sl cf.}  \cite{Vgreen}, Cor. 6.3.10), on a 
\begin{equation}\label{ICsd}
\caR^{d(G_t,L_t)}_{\frqqq_*,L_{t},G_t}( \pi(\phi_L,t))
=\caR^{d(G_t,L_t)}_{\frqqq_*,L_{t},G_t}(\caR^{d(L_t,T_*)}_{\frb_{L_*} , T_*,L_{t}}(\xi_{\phi_{ T_*} (\mathbf B_*)}) )=  
\caR^{d(G_t, T_*)}_{\frb_{*} ,T_{*},G_t }(\xi_{\phi_{ T_*} (\mathbf B_*)})
=\pi(\phi_G,t) \end{equation}
où 
\begin{equation} d(L_t, T_*)=\frac{1}{2}(\dim \mathbf L_*-\dim \mathbf T_*)-(q(L_t)-q( T_*))  .  \end{equation}
Remarquons que ces inductions cohomologiques sont dans le  good range.

Fixons un élément $t_0$ de $T_*[2]$ 
et un   système de représentants $(\mathbf Q_i,\mathbf L_i)$ de  $\Sigma_{{}^dQ}^{G_{t_0}}$
 que l'on suppose $\tau_{G_{t_0}}$-stables. On a 
pour tout $i$, un  foncteur d'induction cohomologique
\[ \caR^{d_i}_{\frqqq_i, L_i,G} \]
en degré
\begin{equation} d_i=d(G,L_i)=\frac{1}{2}(\dim \mathbf G-\dim \mathbf L_*)-(q(G_{t_0})-q(L_{i}))\end{equation}
 de la catégorie des représentations de $L_i$ vers la catégorie des représentations de $G_{t_0}$.

Pour chaque $L_i$, on a un paquet de séries discrètes $\Pi(\phi_L,L_i)$ et la représentation virtuelle 
 stable associée $ [\pi(\phi_L,L_i)] $ (la somme des éléments du paquet).
 \begin{prop}\label{cor:stable}
On a 
\[ \sum_{(\mathbf Q_i,\mathbf L_i)\in \Sigma_{{}^dQ}^{G_{t_0}}}    \; \caR^{d_i}_{\frqqq_i,L_i,G_{t_0}}\left(  [\pi(\phi_L,L_i) ]\right)
 =   [\pi(\phi_G,G_{t_0}) ]  \]
en particulier cette représentation virtuelle est stable.
\end{prop}
\dem  On a introduit la relation d'équivalence $t\sim_G t'$ sur $T_*[2]$ ci-dessus, et l'on peut définir la relation d'équivalence
analogue, mais relativement au groupe $\mathbf L_*$ et que l'on note $\sim_L $.
C'est une relation d'équivalence plus faible que $\sim_G$.
Le membre de droite de l'égalité à démontrer est 
\[   [\pi(\phi_G,G_{t_0}) ]= \sum_{t\in T_*[2], t\sim_G t_0}  [\pi(\phi_G,t)] \] 
Fixons un système de représentants $(t_i)_i$ pour la relation $\sim_L $ dans l'ensemble  $\{ t\in T_*[2]\, \vert\,   t \sim_G  t_0  \}$.
On a donc en utilisant (\ref{ICsd})
\begin{align*}   [\pi(\phi_G,G_{t_0}) ]&= \sum_i \sum_{t\in T_*[2], t\sim_L t_i}  [\pi(\phi_G,t)] = \sum_i \sum_{t\in T_*[2], t\sim_L t_i} 
\left[ \caR^{d(G_{t},L_{t})}_{\frqqq_*,L_{t},G_{t}}( \pi(\phi_L,t)) \right] \\
 &=\sum_i \caR^{d(G_{t_i},L_{t_i})}_{\frqqq_*,L_{t_i},G_{t_i}} (  [\pi(\phi_L, L_{t_i})] ).
\end{align*}
Par conjugaison intérieure des $\mathbf G_{t_i}$ avec $\mathbf G_{t_0}$ les $L_{t_i}$ s'identifient aux $L_i$ de
 $\Sigma_{{}^dQ}^{G_{t_0}}$  et  l'on retrouve bien le membre de gauche dans l'égalité de la proposition.
\qed

\medskip

Revenons maintenant au cas d'un groupe $\mathbf G$ quelconque, et à la donnée d'un 
sous-groupe parabolique standard 
${}^d Q= {}^d L{}^d V$  de $\widehat G$,  compatible avec   l'action de $\sigma_{\widehat G}$, et d'un 
système de représentants $(\mathbf Q_i,\mathbf L_i)$ de  $\Sigma_{{}^dQ}^{G}$
 que l'on suppose $\tau_{G}$-stables.
Soit  $\phi_L: \, W_\bbR\rightarrow {}^LL$  un paramètre de Langlands, et posons toujours $\phi_G=\iota_{L,G}\circ \phi_L$. 
On suppose  qu'on est dans le  good range.
On a les  représentations virtuelles 
 stables $ [\pi(\phi_L,L_i)] $. 

 \begin{prop}\label{cor:stable2}
Il existe un paramètre de Langlands $\phi_G'$ pour $G$ de même caractère infinitésimal que $\phi_G$ tel que    
\[ \sum_{(\mathbf Q_i,\mathbf L_i)\in \Sigma_{{}^dQ}^{G}}      \; \caR^{d_i}_{\frqqq_i,L_i,G} \left( [\pi(\phi_L,L_i) ]\right)  =   [\pi(\phi'_G,G) ]  \]
en particulier cette représentation virtuelle est stable.
\end{prop}

Attention,  $\phi'_G$ n'est pas égal à $\phi_G$, il y a une torsion,
due au fait que l'induction cohomologique ne commute pas à l'induction parabolique,  voir l'appendice.
Evidemment, ici nous n'avons pas précisé ce qu'est $\phi'_G$, mais nous n'en avons pas besoin pour le moment car nous ne visons
que l'énoncé de stabilité. 

La démonstration dans  le cas général demandant un formalisme un peu lourd, nous allons nous contenter 
d'une démonstration dans les cas qui nous servent dans cet article, c'est-à-dire  les groupes classiques.
Nous repoussons la démonstration à la section \ref{demprop}, lorsque nous aurons introduits quelques objets et notations
spécifiques à ces cas. Pour ceux-ci, le calcul de la torsion pour passer du paramètre $\phi_G$ au paramètre $\phi_G'$ est explicite.

\medskip 

On peut appliquer ceci au cas des  formes intérieures pures $\mathbf G_t$
d'un groupe quasi-déployé $\mathbf G$, situation qui a été considérée ci-dessus. On a alors un 
$c$-sous-groupe de Levi $\mathbf L_*$ de $\mathbf G$. 
Comme toute représentation virtuelle stable sur $L_*$ est une combinaison linéaire de 
 $[\pi(\phi_L,L_*) ]$, et qu'une représentation virtuelle stable de $L_*$ détermine par transfert endoscopique
 une  représentation virtuelle stable de toute forme intérieure $L_*$, on obtient 
 ainsi que $\sum_{(\mathbf Q_i,\mathbf L_i)\in \Sigma_{{}^dQ}^{G_{t}}}    e(L_i)  \; \caR^{d_i}_{\frqqq_i,L_i,G_t}  $ envoie les 
 \og représentations virtuelles super-stables de $\mathbf L_*$\fg \,  
 dans le good range   
 sur les représentations virtuelles stables de $G_t$. 
 
 
\begin{cor}\label{cor:stable3} Dans le contexte décrit ci-dessus, soit $ [X_{L_*}^{st}]$ une représentation virtuelle stable de $L^*$, et soit 
$ [X_{L_i}^{st}]$ la représentation stable obtenue par transfert endoscopique de $L_*$ vers sa forme intérieure $L_i$. Alors 
 \[ \sum_{(\mathbf Q_i,\mathbf L_i)\in \Sigma_{{}^dQ}^{G_{t}}}     e(L_i)  \; \caR^{d_i}_{\frqqq_i,L_i,G_t} \left( [X_{L_i}^{st}]\right)  \]
 est stable.
 De plus, lorsqu'on considère toutes les formes intérieures pures $\mathbf G_t$ simultanément, ceci définit une 
 représentation virtuelle superstable, dans le sens où 
 \[  \Trans_{G}^{G_t} \left(\sum_{(\mathbf Q_i,\mathbf L_i)\in \Sigma_{{}^dQ}^{G}}     e(L_i)  \; \caR^{d_i}_{\frqqq_i,L_i,G} \left( [X_{L_i}^{st}]\right) 
 \right)=
 \sum_{(\mathbf Q_i,\mathbf L_i)\in \Sigma_{{}^dQ}^{G_{t}}}     e(L_i)  \; \caR^{d_i}_{\frqqq_i,L_i,G_t} \left( [X_{L_i}^{st}]\right).   \]
\end{cor}
Ceci est aussi établi pour les groupes classiques dans la section \ref{demprop}.

\section{Formes intérieures pures, $c$-Levi et données endoscopiques elliptiques des   groupes classiques}\label{ICgrcla}
    
  \subsection{Paramétrisation des séries discrètes et formes intérieures pures} \label{seriesdiscretescla}
  Les groupes classiques quasi-déployés de la section \ref{grcla} admettent tous des séries discrètes, à l'exception 
 de $\SO_{2n}^d$ pour $n$ impair, et de   $\SO_{2n}^{qd}$ pour $n$ pair, que nous excluons donc  de la discussion jusqu'à  la section suivante.
 Nous allons spécialiser la discussion de la section \ref{seriesdiscretes} au cas où $\mathbf G$ est l'un de ces groupes classiques.
 Rappelons que nous avons fixé une paire fondamentale de type Whittaker $(\mathbf T_*,\mathbf B_*)$, et qu'ici, $\mathbf T_*$ est anisotrope.

 Fixons $\phi_G:\, W_\bbR\rightarrow {}^LG $,  un paramètre  de Langlands discret. Alors  $\Pi^\sharp(\phi_G,G)=\Pi(\phi_G,G)$
est un paquet de séries discrètes.

On identifie le sous-groupe de Cartan $\mathbf T_*$ à $\U(1)^n$ grâce aux racines simples de $\mathbf T_*$ dans $\mathbf B_*$
qui fournissent une base $(e_i)_{i=1, \ldots n}$ de $X^*(\mathbf T_*)$.
Les éléments d'ordre 2 de $ T_*$ seront alors noté $t=(\pm 1, \ldots ,\pm 1)$ et il y en a $2^n$.
Les éléments d'un super-paquet de séries discrètes  $\Pi(\phi_G)$
 sont donc paramétrés par les éléments $t=(\pm 1, \ldots, \pm 1)$ de $ T_*[2]$.
 Rappelons que nous les avons notées $\pi(\phi_G,t)$. 
La série discrète   $\pi(\phi_G,1)$ est générique et admet le modèle de Whittaker fixé au départ.

Pour tout  $t=(\pm 1, \ldots, \pm 1)\in T_*[2]$, on obtient  via le pairing (\ref{pairingAdams})  un caractère 
$\eta_t$ de $A(\phi_G)$. Le caractère trivial de  $A(\phi_G)$ est le paramètre de la série discrète $X(\phi_G,1)$
admettant le modèle de Whittaker.

Il reste à identifier dans chaque cas quelles sont les formes intérieures qui interviennent, et comment les $2^n$ séries discrètes
se répartissent entre ces formes intérieures. Pour cela, posons 
\begin{equation} \label{choixt*} t_*=((-1)^n   \ldots ,1, ,-1)\in T_*[2].\end{equation} 
Pour tout $t\in T_*[2]$, notons  $n_{1}(t)$  (resp.  $n_{-1}(t)$)  le nombre de coordonnées égales à $1$ (resp.  $-1$)  dans $tt_*$. 
\bigskip 

Cas {\bf A} :  $\mathbf G=\Sp_{2n}$. 
Toute les formes intérieures $\mathbf G_t$ sont équivalentes à $\mathbf G$. On obtient ainsi un 
paquet  $\Pi(\phi_G,G)$ avec $2^n$ éléments.
\bigskip 

Cas {\bf B} :  $\mathbf G=\SO_{2n+1}$.  Deux formes $G_t$ et $G_{t'}$ sont équivalentes si et seulement si $n_1(t)=n_1(t')$.
Nous allons noter $\SO(p,q)$  la classe d'équivalence  classe d'équivalence déterminée par la formule 
$p=2n_1(t)+1$ (et donc $q=2n_{-1}(t)$).
 \bigskip

Cas {\bf C} et {\bf D} :  $\mathbf G=\SO_{2n}^\alpha$, avec $\alpha=d$ si $n$ est pair et $\alpha=qd$ si $n$ est impair. 
 Deux formes $G_t$ et $G_{t'}$ sont équivalentes si et seulement si $n_1(t)=n_1(t')$.
Nous allons noter $\SO(p,q)$  la classe d'équivalence déterminée par la formule 
$p=2n_1(t)$ (et donc $q=2n_{-1}(t)$).

 \bigskip

Cas {\bf U} :  $\mathbf G=\U_{n}$. Deux formes intérieures pures $G_t$ et $G_{t'}$ sont équivalentes si et seulement si $n_1(t)=n_1(t')$.
Nous allons noter $\U(p,q)$  la classe d'équivalence donnée par  $n_1(t)=p$ (et donc $q=n_{-1}(t)$).

 \begin{rmqs} \label{rmqs}
 Pour les groupes orthogonaux, on remarque que la classe du groupe quasi-déployé de départ $\mathbf G$ est celle obtenue pour $t=1$, et le 
 choix du $t_*$ en (\ref{choixt*}) nous donne dans le cas {\bf B},  $\SO(n+1,n)$ si $n$ est pair, et $\SO(n,n+1)$ si $n$ est impair, et dans les cas 
   {\bf C} et {\bf D}, $\SO(n,n)$ si $n$ est pair et $\SO(n-1,n+1)$ si $n$ est impair.
 On peut interpréter cette construction de
  la manière suivante : on fixe une suite de $n$ plans réels 
  muni d'une forme quadratique de signature $(2,0)$ ou $(0,2)$, alternativement, {\sl en finissant}
  par un plan de signature $(0,2)$.  Dans le cas {\bf B}, on complète par une droite muni d'une forme quadratique définie positive, 
   pour former par somme directe un espace quadratique $V$
  dont $\mathbf G$ est le groupe de symétries. Le tore compact $T_*$ est le stabilisateur de cette décomposition de $V$.
  Pour tout élément $t\in T_*[2]$, le groupe $\mathbf G_t$ est le groupe de symétrie de l'espace obtenu en changeant 
  la forme du $i$-ième plan en son opposée à chaque fois que la $i$-ème coordonnée de $t$ est $-1$.

Dans les cas  {\bf C} et {\bf D}, remarquons que les classes de formes intérieures pures $\SO(p,q)$ et $\SO(q,p)$ sont distinguées si 
$p\neq q$, alors que les groupes sous-jacents sont isomorphes.
   On distingue aussi les formes intérieures $\U(p,q)$ et $\U(q,p)$ si $p\neq q$.
Par exemple, pour $n=1$, un super-paquet de séries discrètes contient 2 éléments, l'un est une représentation de la forme $\U(1,0)$, 
l'autre  est une représentation de la forme 
$\U(0,1)$. La classe du groupe quasi-déployé $G$ de départ est donc $\U(n/2,n/2)$ si $n$ est pair, et $\U\left(\frac{n-1}{2},\frac{n+1}{2}\right)$
si $n$ est impair. On a une interprétation similaire de $\mathbf G$ et de ses formes intérieures pures $\mathbf G_t$ 
en termes d'un espace hermitien formé par somme directe de droites complexes hermitiennes.
 \end{rmqs}

\begin{rmq}
Cette classification fait apparaître de manière naturelle les  classes d'équivalence de formes intérieures pures de notre groupe de départ.
Rappelons qu'en général celle-ci sont classifiée par le groupe de 
cohomologie  $H^1(\Gamma,\mathbf G)$ et que dans les cas des groupes orthogonaux ou unitaires, 
l'identification de ce groupe avec les classes d'équivalence de formes bilinéaires ou hermitiennes est bien connue.
\end{rmq}

\begin{rmq}\label{ChoiX}
 La formulation d'une partie des  résultats qui suivent va dépendre, de manière évidemment inessentielle, du choix de $t_*$
 en (\ref{choixt*}). Nous avons choisi une des quatre possibilités pour alterner des $\pm 1$ : commencer par la gauche
 ou par la droite, avec un $+1$ ou bien avec un $-1$, en commençant avec $-1$ par la droite. Par exemple, si l'on fait un autre choix, 
 pour $\mathbf G$ groupe orthogonal ou unitaire, la classe d'équivalence de formes intérieures de $\mathbf G$ peut être donnée
 par un $\SO(p,q)$ ou $\U(p,q)$ différent (mais toujours quasi-déployé, {\sl i.e.} $\vert  p-q \vert\leq 1$ ).  Voir aussi la remarque 
 \ref{repdist} pour un énoncé qui dépend de ce choix.
 \end{rmq}

\subsection{Séries fondamentales et formes intérieures pures des groupes orthogonaux sans séries discrètes}\label{subsec51}
Il est utile de faire apparaître les formes intérieures pures des groupes manquants, à
savoir $\SO_{2n}^d$, $n$ impair, et $\SO_{2n}^{qd}$, $n$ pair. Plaçons nous dans le premier cas.
Réalisons ce groupe comme groupe d'isotropie pour un espace quadratique réel $(V,q)$ de dimension 
$2n$ et de signature  $(n,n)$.
Fixons un plan hyperbolique $P$ dans $V$, et notons $W$ son orthogonal. Alors  $\SO_{2(n-1)}^d$ 
est naturellement réalisé comme  groupe d'isotropie de  $W$, et le stabilisateur de 
$P\oplus W$ dans $\SO_{2n}^d$  est isomorphe à $\bbR^\times \times   \SO_{2(n-1)}^d$, qui apparait 
comme sous-groupe de Levi de $ \SO_{2n}^d$. Fixons une paire fondamentale de type Whittaker $(\mathbf B_*, \mathbf T_*)$
de $\SO_{2n}^d$ de sorte que l'intersection avec $\SO_{2(n-1)}^d$ soit une 
 paire fondamentale  $(\mathbf B'_*, \mathbf T'_*)$ de $\SO_{2(n-1)}^d$. Alors $\mathbf T_*=\bbR^\times\times  \mathbf T'_*$ et  
$\mathbf T'_*$ est anisotrope et $\SO_{2(n-1)}^d$ admet des séries discrètes, on peut donc lui appliquer les considérations de la section précédente.
 Chaque $t\in \mathbf T'_*[2]$ donne une forme intérieure $\mathbf G'_t$  de $\SO_{2(n-1)}^d$ et 
une forme intérieure  $\mathbf G_t$  de $\SO_{2n}^d$,  
La  classe d'équivalence de la forme intérieure $\mathbf G'_t$ est $ \SO(p-1,q-1)$, où $p$ et $q$ sont déterminés dans le paragraphe précédent, 
avec ici $q$ impair,  et l'on note $ \SO(p,q)$ la classe d'équivalence de $\mathbf G_t$. 

 Les séries fondamentales de  $ \SO_{2n}^d$ et de ses formes intérieures pures $\mathbf G_t$ sont obtenues
par induction parabolique à partir des séries discrètes de $\SO_{2(n-1)}^d$ et de ses  formes intérieures $\mathbf G'_t$
et d'un caractère du facteur $\bbR^\times$, et ceci préserve 
les paquets de Langlands. Ainsi, les paquets de séries  fondamentales de $ \SO_{2n}^d$ et de ses formes intérieures pures $\mathbf G_t$
 sont paramétrés par les éléments  de $T'_*[2]$.

Le cas  $\SO_{2n}^{qd}$, $n$ pair se traite de la même manière à partir un espace quadratique réel $(V,q)$ de dimension 
$2n$ et de signature  $(n-1,n+1)$.

Pour avoir des notations uniformes qui s'appliquent aussi à ces groupes, on pose  $T_*[2]:=T'[2]$ pour les groupes de cette section. 
Donc $T_*[2]$ n'est pas ici l'ensemble des élément d'ordre 2 de $T_*$, mais l'ensemble des éléments d'ordre 2 de son facteur anisotrope $T'_*$.
Dans ce contexte, si $\phi_G$ est un paramètre de Langlands de séries fondamentales, on a aussi un pairing parfait 
$T_*[2]\times A(\phi_G)\rightarrow \bbC^\times $ qui identifie $T_*[2]$ et $\widehat{A(\phi_G)}$ et une paramétrisation 
des éléments du super-paquet correspondant par  $T_*[2]$ ou $\widehat{A(\phi_G)}$.

 \subsection{$c$-sous-groupes de Levi des groupes classiques}\label{cLevicla}
 Soit $\mathbf G$ un groupe classique quasi-déployé de la section \ref{grcla}

Soit ${}^d Q= {}^d L{}^d V$ un sous-groupe parabolique standard de $\widehat G$,
compatible avec l'action de $\sigma_{\widehat G}$.
On va supposer de plus que ${}^d Q$ est un sous-groupe parabolique maximal de $\widehat G$, c'est-à-dire qu'il est associé
au choix d'une racine simple. 
Considérons le  sous-groupe parabolique $\mathbf Q_*=\mathbf L_*\mathbf V_*$ obtenu à partir de la paire de Borel 
fondamentale de type Whittaker $(\mathbf B_*, \mathbf T_*)$. Nous avons vu que $\mathbf L_*$ est quasi-déployé.
Il est isomorphe à un produit 
\[\U_c\times \mathbf G_{n-c}\]
où $\U_c$ est un groupe unitaire quasi-déployé de rang $c$, et $\mathbf G_{n-c} $ est un 
groupe classique quasi-déployé de même type que $\mathbf G$ (sauf dans le cas $c$ impair ou les types {\bf C} et {\bf D} sont échangés).
L'entier $c$ se détermine facilement selon la place de la racine simple $\alpha$ 
mentionnée ci-dessus dans le diagramme de Dynkin.

Rappelons qu'à tout $t\in T_*[2]$, on a associé une conjugaison complexe $\sigma_t$ et des groupes 
$\mathbf G_t$ et $\mathbf L_t$. Nous avons vu que les formes intérieures 
$G_t$ et $G_{t'}$ sont équivalentes si et seulement si $n_1(t)=n_1(t')$, en dehors du cas 
{\bf A} où elles sont toutes équivalentes. En particulier, les formes
équivalentes à $\mathbf G$ sont celles telles que $n_1(t)=\lfloor \frac{n}{2}\rfloor$ (en dehors du cas 
{\bf A}).
L'équivalence des  formes intérieures $\mathbf L_t$ de $\mathbf L_*$, est-elle donnée
par un invariant de même type, mais où l'on sépare les coordonnées de $t$ en deux, les $c$
premières, qui vont correspondre au facteur $\U_c$, et les $n-c$ dernières, qui correspondent au facteur 
 $\mathbf G_{n-c}$. La recette est donc la suivante : on note $n_{1,L}(t)=(n_{1,U}(t),n'_{1}(t))$ le couple dont la première coordonnée est le nombre de 
$1$ dans les $c$ premières coordonnées de $tt_{*}$ et la seconde le nombre de 
$1$ dans les $n-c$ dernières coordonnées de $tt_{*}$. On définit de même de manière évidente
 $n_{-1,L}=(n_{-1,U},n'_{-1})$. Deux formes $\mathbf L_t$ et $\mathbf L_{t'}$
 sont alors équivalentes si et seulement si $n_{1,L}(t)=n_{1,L}(t')$, sauf dans le cas {\bf A}
où elles sont équivalentes  si et seulement si $n_{1,U}(t)=n_{1,U}(t')$.

Remarquons aussi que les $c$ premières coordonnées de $t_*$ sont $((-1)^n, \ldots ,(-1)^{n-c+1})$ et que ceci ne coïncide
avec le choix du $t_*$ pour les groupes unitaires que si $n-c$ est pair.
La forme du groupe unitaire qui apparaît est donc $\U(n_{1,U}(t), n_{-1,U}(t))$ si  $n-c$ est pair, mais c'est 
$\U(n_{-1,U}(t), n_{1,U}(t))$ si $n-c$ est impair.

On obtient des systèmes de représentants $(\mathbf Q_i=\mathbf L_i\mathbf V_i)_{i=0,\ldots,c}$ de $\Sigma_{{}^d Q}^G$ explicites, avec de plus  :
\begin{equation}\label{Li}  L_i\simeq \U(i,c-i)\times G'_i \end{equation}
où $G'_i$ est une forme intérieure d'un groupe  $\mathbf G'$ de même type que $G$ (les cas {\bf C} et {\bf D} sont échangés si $c$ est impair).
C'est-à-dire que pour chaque $i=1,\ldots,c$, on fixe un $t_i\in \mathbf T_*[2]$ tel que $\mathbf L_{t_i}$ est isomorphe à $L_i$.
Les isomorphismes (\ref{Li}) peuvent être rendu explicites, nous ne le faisons pas, mais nous nous en servons pour identifier
leurs  deux membres.

\begin{rmq}\label{repdist}
Le représentant distingué $\mathbf Q_*=\mathbf L_*\mathbf V_*$ est obtenu pour la valeur suivante de $i$ : 

- si $c$ est pair,  $i=c/2$ 

- si $c$ est impair, et $n-c$ est pair,  $i=\frac{c-1}{2}$, 

- si $c$ est impair, et $n-c$ est impair, $i=\frac{c+1}{2}$

 Le $t_i$ correspondant est alors égal à $1$.
\end{rmq}

Fixons maintenant, pour les groupes orthogonaux, une classe d'équivalence $\SO(p,q)$ de formes intérieures pures $\mathbf G_t$. On obtient 
aussi facilement 
un  système de représentants $(\mathbf Q_i=\mathbf L_i\mathbf V_i)$ de $\Sigma_{{}^d Q}^{\SO(p,q)}$ explicite, avec de plus  :
\begin{equation}\label{LiFi}  L_i\simeq \U(i,c-i)\times G'_i \end{equation}
où là encore $G'_i$ est une forme intérieure d'un groupe  $\mathbf G'$ de même type que $G$ 
(les cas {\bf C} et {\bf D} sont échangés si $c$ est impair).
Les indices $i$ qui apparaissent sont ceux vérifiant $2i\leq p$, $2(c-i)\leq q$, et la classe de $G'_i$ 
est alors $\SO(p',q')$ avec $p=p'+2i$, $q=q'+2(c-i)$.

\subsection{Démonstration de la proposition \ref{cor:stable2}
et du corollaire \ref{cor:stable3} pour les groupes classiques}\label{demprop}

On se place dans le le même contexte que la section précédente. On considère un paramètre de Langlands $\phi_L$
pour le groupe $\mathbf L_*$ et ses formes intérieures $\mathbf L_i$, où 
 $(\mathbf Q_i=\mathbf L_i\mathbf V_i)_i$ 
décrit un système de représentants de $\Sigma_{{}^d Q}^{G}$ de la forme (\ref{Li}). On pose 
$\phi_G=\iota_{L,G}\circ \phi_L$.   On suppose que l'on se trouve dans le  good range  
pour l'induction cohomologique.  Les paquets $\Pi(\phi_L,L_i)$ sont associés à une classe de conjugaison 
de sous-groupes de Cartan dans le groupe $L_i$, ou bien le paquet est vide.
  Dans le premier cas fixons un représentant $D_i$ pour chaque $i$, le sous-groupe de Cartan 
$D_i$ est alors isomorphe comme tore réel à un produit de facteurs de la forme $\U(1)$, $\bbR^\times$ ou $\bbC^\times$ : 
\[D_i\simeq\U(1)^{r_1}\times (\bbC^\times)^{m_1}\times  \U(1)^{r_2}\times (\bbC^\times)^{m_2}\times (\bbR^\times)^{s_2},\]
le facteur $\U(1)^{r_1}\times (\bbC^\times)^{m_1}$ correspondant à un sous-groupe de Cartan du facteur unitaire $\U(i,c-i)$
dans la décomposition (\ref{Li}), et le facteur $\U(1)^{r_2}\times (\bbC^\times)^{m_2}\times (\bbR^\times)^{s_2}$ à un sous-groupe de Cartan 
du facteur $G'_i$. On a donc 
$r_1+2m_1=c$,  $r_2+2m_2+s_2=n-c$, ces entiers ne dépendant que de $\phi_L$ et pas de l'indice $i$ 
(le cas où  $L_i$ ne possède pas de sous-groupe de Cartan
de cette forme est celui ou le paramètre $\phi_L$ n'est pas \og relevant \fg\, pour $L_i$, c'est-à-dire que le paquet 
$\Pi(\phi_L,L_i)$ est vide). Le centralisateur de la partie déployée de $D_i$ définit un sous-groupe de Levi cuspidal 
$\mathbf M_{L_i}$ de $\mathbf L_i$ et l'on a 
\begin{equation}\label{MLi} M_{L_i}\simeq  \U(i-m_1,c-i-m_1) \times (\bbC^\times)^{m_1}\times   G'_{r_2}\times \GL_2(\bbR)^{m_2}\times (\bbR^\times)^{s_2},  
\end{equation}
où $\mathbf G'_{r_2}$ est un   groupe classique de même type que $\mathbf G$ et de rang $r_2$.
Les paquets $\Pi(\phi_L,L_i)$ sont obtenus à partir des  paquets de séries discrètes  $\Pi(\phi_{M_L}, M_{L_i})$ attaché à un paramètre discret 
$\phi_{M_L} $ et les représentations virtuelles stables associées s'écrivent comme en (\ref{StablePsP2}) : 
\[  [\pi(\phi_L,L_i)]=\sum_{\pi_{M_{L_i}}  \in  \Pi(\phi_{M_L}, M_{L_i})} [\Ind_{M_{L_i}}^{L_i}(\pi_{M_{L_i}})]  \]

D'autre part, le groupe $G$ admet une classe de conjugaison de sous-groupes de Cartan isomorphes à 
\[ \U(1)^{r_1+r_2}\times (\bbC^\times)^{m_1+m_2}\times (\bbR^\times)^{s_2} .\]
Fixons un représentant $D$ de cette classe de conjugaison, et soit $\mathbf M$
le sous-groupe de Levi cuspidal de $\mathbf G$ obtenu en prenant le centralisateur de la partie déployée de $D$.
On a 
\[  M\simeq  G'_{r_1+r_2}\times \GL_2(\bbR)^{m_1+m_2}\times (\bbR^\times)^{s_2},  \]
où $\mathbf G'_{r_1+r_2}$ est un groupe classique de même type que $\mathbf G$ et de rang $r_1+r_2$.

On applique les considérations des sections \ref{cLevi} et \ref{cLevicla} aux groupes $\mathbf M$ et $\mathbf G'_{r_1+r_2}$ : 
il  existe un système de représentants de classes de conjugaison de $c$-paraboliques 
\[  (\mathbf Q_k^\sharp,\mathbf L_k^\sharp)_k \]
de $\mathbf M$, avec 
\begin{equation} \label{Lksharp}
L_k^\sharp\simeq \U(k,c-k)\times G'_k\times (\bbC^\times)^{m_1}\times \GL_2(\bbR)^{m_2}\times (\bbR^\times )^{s_2}, 
\end{equation}
où les $G'_k$ sont des formes intérieures pures d'un groupe classique de même type que $\mathbf G$ et de rang
$r_1+r_2-c$. De même, et  de manière compatible (dans un sens évident), il existe  
un système de représentants de classes de conjugaison de $c$-paraboliques 
\[  (\mathbf Q_k,\mathbf L_k)_k \]
de $\mathbf G_{r_1+r_2}$, avec 
\[  L_k\simeq \U(k,c-k)\times G'_k. \]
Pour tout indice $i$, le sous-groupe $M_{L_i}$ de  (\ref{MLi}) est conjugué dans $G$ à un sous-groupe 
$L_k^\sharp$ de (\ref{Lksharp}), pour un unique indice $k$, ce qui donne une bijection entre 
l'ensemble des indices $i$ et l'ensemble des indices $k$.

Dans ce contexte, le membre de gauche dans l'égalité de la proposition \ref{cor:stable2} est donc 
\[  \sum_{(\mathbf Q_i,\mathbf L_i)\in \Sigma_{{}^dQ}^{G}}      \; \caR^{d_i}_{\frqqq_i,L_i,G} \left( [\pi(\phi_L,L_i) ]\right)  = 
 \sum_{(\mathbf Q_i,\mathbf L_i)\in \Sigma_{{}^dQ}^{G}}    \sum_{\pi_{M_{L_i}}  \in  \Pi(\phi_{M_L}, M_{L_i})} 
  \; \caR^{d_i}_{\frqqq_i,L_i,G} \left(     [\Ind_{M_{L_i}}^{L_i}(\pi_{M_{L_i}})]  \right). \]
En utilisant la proposition \ref{TorsionI} de l'appendice, la représentation 
$\caR^{d_i}_{\frqqq_i,L_i,G} \left(     [\Ind_{M_{L_i}}^{L_i}(\pi_{M_{L_i}})]  \right)$ est de la forme 
$\Ind_P^G(\pi_M)$ pour une certaine série discrète $\pi_M$ de $M$, 
et lorsque $\pi_{M_{L_i}}$ décrit le paquet $\Pi(\phi_{M_L}, M_{L_i})$, 
$\pi_M$ décrit un ensemble de la forme  
 \[\left\{ \caR^{d_k}_{\frqqq_k,L_k^\sharp,M} (\pi_{L_k^\sharp}), \; 
\pi_{L_k^\sharp} \in \Pi( \phi_{L^\sharp}, L_k^\sharp ) \right\} \]
où $\phi_{L^\sharp} $ est un paramètre discret qui s'obtient à partir de $\phi_{M_L}$ par une torsion
très simple sur le facteur $(\bbR^\times)^{s_2}$ (voir section \ref{calculTorsionI} dans l'appendice).
On peut donc continuer le calcul du 
 membre de gauche dans l'égalité de la proposition \ref{cor:stable2}, que l'on écrit sous la forme 
\[ \sum_k   \sum_{\pi_{L_k^\sharp} \in \Pi( \phi_{L^\sharp}, L_k^\sharp ) } \left[   \Ind_M^G  \left(  \caR^{d_k}_{\frqqq_k,L_k^\sharp,M} 
(\pi_{L_k^\sharp} )  \right)\right]=
\Ind_M^G \left(  \sum_k   \sum_{\pi_{L_k^\sharp} \in \Pi( \phi_{L^\sharp}, L_k^\sharp ) }      \caR^{d_k}_{\frqqq_k,L_k^\sharp,M}  
(\pi_{L_k^\sharp} )    \right). \]

On déduit facilement du cas des séries discrètes établi dans la proposition  \ref{cor:stable}, que l'on applique au 
groupe $\mathbf G_{r_1+r_2}$ et aux $(\mathbf Q_k,\mathbf L_k)_k$, un énoncé analogue pour le groupe $\mathbf M$ et 
les $(\mathbf Q_k^\sharp,\mathbf L_k^\sharp)_k $, à savoir 
 \[ \sum_k   \sum_{\pi_{L_k^\sharp} \in \Pi( \phi_{L^\sharp}, L_k^\sharp ) }      \caR^{d_k}_{\frqqq_k,L_k^\sharp,M}  
(\pi_{L_k^\sharp} )  = \sum_k       \caR^{d_k}_{\frqqq_k,L_k^\sharp,M}  
([\pi(\phi_{L^\sharp},L_k^\sharp )] )= [\pi(\phi'_M,M)] \]
pour le paramètre discret $\phi'_M$ obtenu à partir de $\phi_{L^\sharp}$ par composition avec le plongement de $L$-groupe 
défini comme en (\ref{iotaLG}) mais pour $L^\sharp$ et $M$.

Ainsi, on peut conclure, le membre de gauche dans l'égalité de la proposition \ref{cor:stable2} est bien
\[ \Ind_M^G([ \pi(\phi'_M,M) ])= [\pi(\phi'_G,G)]\]
avec $\phi'_G=\iota_{M,G}\circ \phi_M'$.

La démonstration est la même, à quelques adaptations évidentes près, si on remplace le groupe classique
quasi-déployé  $\mathbf G$ par une forme intérieure
pure $\SO(p,q)$ d'un groupe orthogonal et le corollaire \ref{cor:stable3} en découle immédiatement.

\subsection{Données endoscopiques elliptiques}\label{DEE}

Soit $\mathbf G$ l'un des groupes classiques de la section \ref{grcla}.
Nous renvoyons à \cite [\S 2.1]{KS99}  pour la définition d'une donnée endoscopique $\underline H=(\mathbf H,\caH,x,\xi)$ 
de $\mathbf G$. En particulier, $\caH$ est une extension scindable de $W_\bbR$ par $\widehat H$, et $\xi:\, \caH\rightarrow {}^LG$ est un 
$L$-plongement. Nous suivons maintenant \cite[\S 1.8]{WaldForm} pour inclure dans la donnée endoscopique
un $L$-isomorphisme ${}^L\xi:\, {}^LH\rightarrow \caH$, qui a le bon goût de toujours exister pour les données
endoscopiques des groupes de la section \ref{grcla}. Ceci nous dispense d'introduire les $z$-extensions
de $\mathbf H$ de \cite{KS99}.  En composant $\xi$ avec cet $L$-isomorphisme, on obtient un plongement 
$\iota_{H,G}: \, {}^LH\rightarrow {}^LG$.  

Voici, pour les groupes de la section \ref{grcla}, une liste de représentants des classes d'équivalence de données
endoscopiques elliptiques. On donne ici seulement le groupe endoscopique $\mathbf H$, et l'on renvoit à 
\cite{WaldForm} pour les autres éléments de la donnée.

\medskip

 \begin{itemize}
 \item[\bf Cas  A] :   $\mathbf G=\Sp_{2n}$, $\mathbf H=\Sp_{2a}\times \SO_{2b}^\alpha $ où $a+b=n$, 
 $\alpha\in \{d,qd\}$ et $(b,\alpha)\neq (1,d)$ et si $b=0$, alors $\alpha=d$.
 
 \medskip

  \item[\bf Cas B] :  $\mathbf G=\SO_{2n+1}$,  $\mathbf H=\SO_{2a+1}\times \SO_{2b+1} $ où $a+b=n$.
    On obtient une donnée endoscopique équivalente en échangeant $a$ et $b$.
    \medskip
   
    \item[\bf Cas C,  D] :    $\mathbf G=\SO^\alpha_{2n}$, $\mathbf H=\SO^\beta_{2a}\times \SO_{2b}^\gamma $ où $a+b=n$, 
 $\alpha,\beta,\gamma \in \{d,qd\}$, $\beta=d$ si $a=0$, $\gamma=d$ si $b=0$,
$(a,\beta)\neq (1,d)$,  $(b,\gamma)\neq (1,d)$ et 
  $\alpha=\beta\gamma$ 
 où l'on met sur $\{d,qd\}$ la structure de groupe à deux éléments avec $d$ comme élément neutre.
 Les couples $((a,\beta), (b,\gamma))$ et  $( (b,\gamma), (a,\beta))$
 définissent des   données endoscopiques équivalentes.
    \end{itemize}

\medskip 


\begin{rmq}\label{vp1}
Le groupe endoscopique $\mathbf H$ est écrit ci-dessus sous la forme
\[ \mathbf H=\mathbf H_1\times \mathbf H_2.\]
On peut supposer que le $x$ de la donnée endoscopique soit un élément d'ordre $2$ de $\widehat G$.
Les groupes $\mathbf H_1$ et $\mathbf H_2$ ci-dessus correspondent alors aux 
valeurs propres $\pm 1$ de l'action adjointe de $x$, mais on ne suppose pas que $\mathbf H_1$ corresponde à la valeur
propre $1$.
\end{rmq}

\begin{rmq}\label{plongH}
La relation entre les plongements $\iota_{H,G}$ définis par Waldspurger
et les représentations standard des $L$-groupes est la suivante. Dans tous les cas, sauf le cas $\bf A$
ci-dessus avec $\mathbf H=\Sp_{2a}\times \SO_{2b}^{qd} $, la représentation $\Std_G\circ \iota_{H,G}$ de ${}^LH$
est conjuguée à $(\Std_{H_1}\oplus\Std_{H_2})\circ \iota_{H,H_1\times H_2}$, où $\iota_{H,H_1\times H_2}$ est l'inclusion de 
${}^LH$ dans ${}^LH_1\times {}^LH_2$. Dans le cas $\mathbf G=\Sp_{2n}$, $\mathbf H=\Sp_{2a}\times \SO_{2b}^{qd} $,
$\Std_G\circ \iota_{H,G}$ est conjuguée à 
$\left((\Std_{\Sp_{2a}}\otimes \sgn_{W_\bbR})\oplus \Std_{\SO_{2b}^{qd}}\right)\circ\iota_{H,H_1\times H_2}$ où 
$\sgn_{W_\bbR}$ est le caractère quadratique non trivial de $W_\bbR$.
\end{rmq}

\subsection{$c$-sous-groupes de Levi maximaux des groupes endoscopiques des groupes classiques}\label{cLeviH}

Soit $\mathbf G$ un groupe classique quasi-déployé de la section \ref{grcla} et $\underline H=(\mathbf H,\caH,x,\xi_{L,H})$ 
  une donnée endoscopique elliptique de $\mathbf G$ comme dans la section \ref{DEE}.
  
 Soit $\mathbf{spl}_{\widehat H}=(\caB_H=\caB\cap \widehat H, \caT,  \{X_{\alpha}\} )$ 
 l'épinglage de $\widehat H$ obtenu à partir de l'épinglage  $\mathbf{spl}_{\widehat G}$ déjà fixé
 (On suppose le $x$ de la donnée endoscopique dans $\caT$).
  
  On fait pour $\mathbf H$ la construction de la section \ref{cLevi} :
  On fixe  ${}^d Q_H= {}^d L_H{}^d V_H$ un sous-groupe parabolique standard de $\widehat H$, compatible avec 
l'action de $\sigma_{\widehat H}$ sur les racines simples.

Toute   paire de Borel fondamentale $(\mathbf B_H, \mathbf T_H)$  dans $\mathbf H$ 
permet de définir un sous-groupe parabolique
$\mathbf Q_H=\mathbf L_H\mathbf V_H$ de $\mathbf H$, avec $\mathbf L_H$ défini sur $\bbR$.

  On fixe une paire de Borel fondamentale et de type Whittaker $(\mathbf B_{*,H}, \mathbf T_{*,H})$
  de $\mathbf H$
  ce qui nous fournit $\mathbf Q_{*,H}=\mathbf L_{*,H}\mathbf V_{*,H}$, avec $\mathbf L_{*,H}$ quasi-déployé. 
    
    \medskip 
    
  On suppose que ${}^d Q_H$ est maximal. Voyons les différents cas possibles. On note $\alpha\mapsto \bar \alpha $
  la permutation non triviale de $\{d,qd\}$ et si $\alpha\in \{d,qd\}$, on note 
 $\alpha_c=\alpha$ si $c$ est pair et $\alpha_c=\bar \alpha$ si $c$ est impair.
Notons que c'est la remarque \ref{vp1} qui impose de distinguer les cas {\bf A}1
 et    {\bf A}2.
    \medskip 
  
  \begin{itemize}
  \item[{\bf A}1.]  $\mathbf G=\Sp_{2n}$, $\mathbf H=\Sp_{2a}\times \SO_{2b}^\alpha =\mathbf H_1\times \mathbf H_2$, 

 \noindent  $\mathbf  L_{*,H}=\U_c  \times\Sp_{2(a-c)}\times \SO_{2b}^\alpha =\U_c\times \mathbf  H'_1\times \mathbf H_2=
 \mathbf  L_{*,H_1}\times \mathbf H_2$, 
  
 \noindent   $\mathbf  L_*=\U_c  \times\Sp_{2(n-c)}=\U_c  \times \mathbf G'$.
  
  \item[{\bf A}2.]  $\mathbf G=\Sp_{2n}$, $\mathbf H= \SO_{2b}^\alpha \times \Sp_{2a}  =\mathbf H_1\times \mathbf H_2$, 

 \noindent $\mathbf  L_{*,H}=  \U_c  \times \SO_{2(b-c)}^{\alpha_c}\times  \Sp_{2a}= \U_c\times \mathbf H'_1\times \mathbf H_2 =
 \mathbf  L_{*,H_1}\times \mathbf H_2$,
    
     \noindent $\mathbf  L_*= \U_c\times \Sp_{2(n-c)}= \U_c  \times \mathbf G'$.
 
 \item[{\bf B}.]   $\mathbf G=\SO_{2n+1}$,  $\mathbf H=\SO_{2a+1}\times \SO_{2b+1} =\mathbf H_1\times \mathbf H_2$ 

   \noindent  $\mathbf L_{*,H}=\U_c  \times \SO_{2(a-c)+1}\times \SO_{2b+1}=\U_c\times \mathbf H'_1\times \mathbf H_2 =
 \mathbf  L_{*,H_1}\times \mathbf H_2$, 
 
   \noindent  $\mathbf L_*=\U_c  \times \SO_{2(n-c)+1}=\U_c  \times \mathbf G'$.
  
   \item[{\bf C},{\bf D}.]  $\mathbf G=\SO^{\alpha}_{2n}$, $\mathbf H=\SO^{\beta}_{2a}\times \SO_{2b}^{\gamma} =\mathbf H_1\times \mathbf H_2$ , 
 
    \noindent  $\mathbf L_{1,H}=\U_c  \times \SO_{2(a-c)}^{\beta_c}\times \SO_{2b}^{\gamma}=\U_c\times \mathbf  H'_1\times \mathbf  H_2 =
 \mathbf  L_{*,H_1}\times\mathbf  H_2$, 

   \noindent  $\mathbf L_{1}=\U_c  \times \SO_{2(n-c)}^{\alpha_c}= \U_c  \times \mathbf G' $.
  \end{itemize}
  
  \medskip
  
  Dans chaque cas, on considère un $c$-sous-groupe de Levi $\mathbf L_*$ de $\mathbf G$ correspondant, produit du même
  groupe unitaire que $\mathbf L_{*,H}$  de rang $c$ et du groupe classique quasi-déployé $\mathbf G'$ de même type que $\mathbf G$
 et  de rang $n-c$. 
  Le groupe $\mathbf  L_{*,H}$ est un $c$-sous-groupe de Levi de $\mathbf  H$, car 
  $\mathbf  L_{*,H_1}$ est un $c$-sous-groupe de Levi de $\mathbf  H_1$.
   
     Ce groupe $\mathbf L_*$ se retrouve  plus canoniquement de la manière suivante. 
  On a  ${}^dQ_H={}^dL_H{}^dV_H\subset \widehat H\subset \widehat G$, et comme ${}^dQ$ est maximal, 
  ce sous-groupe parabolique  standard correspond au choix d'une racine simple   de $\caT$ dans $\caB_{\widehat H}$.
  Le groupe $\widehat H$ est lui le centralisateur de $s$ dans $\widehat G$, il est obtenu en cassant en deux le 
  diagramme de Dynkin étendu de $\widehat G$. Les racines simples de $\caT$ dans $\caB_{\widehat H}$
  sont donc obtenues à partir des racines simples de $\caB_{\widehat G}$
  en en ajoutant une (le sommet supplémentaire du diagramme de Dynkin étendu), et en en  enlevant une, ce qui 
   casse le diagramme de Dynkin étendu 
  en deux. 
  
  La construction de la section \ref{cLevi} donne donc pour toute paire fondamentale $(\mathbf B,\mathbf T)$ de $\mathbf G$
  un sous-groupe parabolique $\mathbf Q=\mathbf L\mathbf V$, etc...
  En particulier, on a la paire fondamentale de type Whittaker $(\mathbf B_*,  \mathbf T_*)$ qui donne 
$\mathbf Q_*=\mathbf L_*\mathbf V_*$ avec $\mathbf L_*$ quasi-déployé.

    Le groupe $\mathbf L_{*,H}$ est un groupe endoscopique pour $\mathbf L_*$ (et ses formes intérieures pures $\mathbf L_{t}$, $t\in T_*[2]$),
    et ceci car $\mathbf H'_1\times \mathbf H_2$ est un groupe endoscopique pour $\mathbf G'$ et ses formes intérieures.
  Ces groupes endoscopiques se    complètent de manière évidente en une  donnée endoscopique.

  \section{Le transfert  endoscopique ne  commute pas à l'induction cohomologique, mais presque}\label{ICetTrans}
  
  \subsection{Considérations préliminaires}\label{Prel}
  
  On se place dans le contexte de la section \ref{cLeviH} : $\mathbf G$ est 
  un groupe classique quasi-déployé de la section \ref{grcla} et $\underline H=(\mathbf H,\caH,x,\xi_{L,H})$ 
  une donnée endoscopique elliptique de $\mathbf G$ comme dans la section \ref{DEE}.
  On fixe  ${}^d Q_H= {}^d L_H{}^d V_H$ un sous-groupe parabolique standard maximal  de $\widehat H$, compatible avec 
l'action de $\sigma_{\widehat H}$ sur les racines simples.
On fixe une paire de Borel fondamentale et de type Whittaker $(\mathbf B_{*,H}, \mathbf T_{*,H})$
  de $\mathbf H$ ce qui nous fournit $\mathbf Q_{*,H}=\mathbf L_{*,H}\mathbf V_{*,H}$, avec $\mathbf L_{*,H}$ quasi-déployé, ainsi
  qu'une paire  de Borel fondamentale et de type Whittaker $(\mathbf B_{L_{*,H}}, \mathbf T_{*,H})$ de $\mathbf L_{*,H}$.

On obtient aussi à partir de la paire de  Borel fondamentale et de type Whittaker $(\mathbf B_*, \mathbf T_*)$
  de $\mathbf G$ un $c$-sous-groupe de Levi quasi-déployé $\mathbf L_*$ de $\mathbf G$,  ainsi
  qu'une paire  de Borel fondamentale et de type Whittaker $(\mathbf B_{*,L}, \mathbf T_*)$ de $\mathbf L_*$.

  Notons $\{\mathbf Q_i=\mathbf L_i\mathbf V_i\}_{i=0,\ldots ,c}$ un système de représentants de $\Sigma_{{}^dQ}^G$ 
  comme en (\ref{Li})  et  
  $\{\mathbf Q_{k,H}=\mathbf L_{k,H}\mathbf V_{k,H}\}_{k=0,\ldots ,c}$ un système de représentants de $\Sigma_{{}^dQ_H}^H$ avec la même
  propriété, c'est-à-dire :
  \begin{equation}\label{LkH*} L_{k,H}\simeq \U(k,c-k) \times H'_k \times H_2. \end{equation}

  On a aussi les applications de transfert spectraux $\Trans_H^G$, et pour tout $i=0,\ldots,c$,    $\Trans_{L_{*,H}}^{L_i} $.
   Les transferts sont normalisés par 
  le choix de  données de Whittaker, c'est-à-dire comme nous l'avons vu, par les paires fondamentales de type Whittaker
  fixées. Le transfert $ \Trans_{L_{*,H}}^{L_i}$, via les identifications (\ref{Li}),  et (\ref{LkH*}) est égal à un signe près à 
   $\Trans_{\U_c}^{\U(i,c-i)}\boxtimes \Trans_{H'_1\times H_2}^{G'_i}$. Expliquons comment déterminer ce signe.
   Définissons tout d'abord un élément $x_d$ de $\bbZ/2\bbZ=\{\pm 1\}$ en posant 
    $x_d=1$ si  le facteur $\mathbf H_1$  du groupe endoscopique et qui contient $\U_c$ correspond à la valeur propre $1$
     du $x$ de la donnée endoscopique, et    $x_d=-1$  s'il correspond à la valeur $-1$  (voir la remarque \ref{vp1}).
Définissons ensuite un élément $S_i$ de $\widehat{\bbZ/2\bbZ}$.  Note $\sgn$ le caractère non-trivial de $\bbZ/2\bbZ$,
     et posons 
      \begin{equation}\label{Si} S_i=\begin{cases} \sgn^{i+\frac{c(c- 1)}{2}} \text{ si } n-c \text{ est pair} \\
  \sgn^{i+\frac{c(c+1)}{2}} \text{ si } n-c \text{ est impair} 
  \end{cases}  = \sgn^{i+\frac{c(c- 1)}{2}+(n-c)c}.\end{equation}
     Remarquons que $S_i$ va alterner selon la parité de $i$ et que la définition est faite de sorte que 
    pour   le représentant distingué  $\mathbf L_*$, on ait le caractère trivial, voir la remarque \ref{repdist}. 
          Avec cette définition, on obtient 
\begin{equation} \label{signtransfU}
  \Trans_{L_{*,H}}^{L_i} = S_i(x_d)   \;  \Trans_{\U_c}^{\U(i,c-i)}\boxtimes \Trans_{H'_1\times H_2}^{G'_i}   .\end{equation}
   Soit $[X_{L_{*,H}}^{st}]$ une représentation virtuelle stable de $L_{*,H}$, admettant un caractère infinitésimal dans le  good range  
  pour $\mathbf Q_{*,H}$. Par linéarité, on peut la supposer de la forme 
  \[  [X_{L_{*,H}}^{st}]= [X_{\U_c}^{st}] \boxtimes [X_{H'_1}^{st}]\boxtimes [X_{H_2}^{st}]\]
  
  Elle détermine par transfert endoscopique une représentation virtuelle stable  $[X_{L_{k,H}}^{st}]$ de tous les $L_{k,H}$, que l'on écrit aussi 
   \[  [X_{L_{k,H}}^{st}]= [X_{\U(k,c-k)}^{st}] \boxtimes [X_{H'_k}^{st}]\boxtimes [X_{H_2}^{st}].\]

On peut donc former: 
\[\sum_k   e(L_{k,H})   \; \caR^{d_k}_{\frqqq_{k,H}, L_{k,H},H} ([ X_{L_{k,H}}^{st}])=
 \left( \sum_k   e(L_{k,H_1})   \; \caR^{d_k}_{\frqqq_{k,H_1}, L_{k,H_1},H_1} \left( [X_{\U(k,c-k)}^{st}] \boxtimes [X_{H'_k}^{st}]\right)\
  \right) \boxtimes [X_{H_2}^{st}] , \]
qui d'après le corollaire \ref{cor:stable3} est une représentation stable de $H$. 
On a utilisé $e(L_{k,H})=e(L_{k,H_1})\times e(H_2)=e(L_{k,H_1})$ car $H_2$ est quasi-déployé.
On peut donc la transférer à $G$, pour obtenir, 
\begin{align} \label{LTransIC}
&\Trans_H^G\left(\sum_k  e(L_{k,H})   \;  \caR^{d_k}_{\frqqq_{k,H}, L_{k,H},H}( [ X_{L_{k,H}}^{st}]) \right)
\\ \nonumber
=&\Trans_H^G\left( \left( \sum_k   e(L_{k,H_1})   \; \caR^{d_k}_{\frqqq_{k,H_1}, L_{k,H_1},H_1} \left( [X_{\U(k,c-k)}^{st}] \boxtimes [X_{H'_k}^{st}]\right)\
  \right) \boxtimes [X_{H_2}^{st}] \right)  \end{align}

D'autre part, on peut aussi transférer de $L_{*,H}$ à  tous les $L_i$  la représentation virtuelle stable  $[X_{L_{*,H}}^{st}]$, pour obtenir 
en vertu de (\ref{signtransfU}) :
\begin{align*} \Trans_{L_{*,H}}^{L_i} \left( [X_{L_{*,H}}^{st}] \right) & =   S_i(x_d)  \Trans_{\U_c}^{\U(i,c-i)}( [X_{\U_c}^{st}] )  \boxtimes 
\Trans_{H'_1\times H_2}^{G'_i}  \left( [X_{H'_1}^{st}]\boxtimes [X_{H_2}^{st}]  \right)    \\
&   =    S_i(x_d)  [X_{\U(c,c-i)}^{st}] )  \boxtimes 
\Trans_{H'_1\times H_2}^{G'_i}  \left( [X_{H'_1}^{st}]\boxtimes [X_{H_2}^{st}]  \right)    \end{align*}
 On applique alors à ceci les foncteurs d'induction cohomologique 
de $L_i$ à $G$, et l'on suppose que l'on est là aussi dans le  good range   : 
\begin{align} \label{RTransIC}
&\sum_i  e(L_i \; )  \caR^{d_i}_{\frqqq_{i},L_i,G}  \left(\Trans_{L_{*,H}}^{L_i} \left( [X_{L_{*,H}}^{st}] \right)   \right)
\\ \nonumber 
=& \sum_i  e(L_i \; )     S_i(x_d) \;  \caR^{d_i}_{\frqqq_{i},L_i,G}  \left(   [X_{\U(i,c-i)}^{st}]   \boxtimes 
\Trans_{H'_1\times H_2}^{G'_i}  \left( [X_{H'_1}^{st}]\boxtimes [X_{H_2}^{st}]  \right) \right).\end{align}

Or (\ref{LTransIC}) et (\ref{RTransIC}) ne sont pas égaux, comme on s'en aperçoit aisément par le fait que le caractère infinitésimal
n'est pas préservé. En effet les foncteurs $ \caR^{d_k}_{\frqqq_{k,H}}$ de $L_H$ à $H$
translatent le caractère infinitésimal de $\rho_{V_H}=\rho_H-\rho_{L_H}$ et 
et les foncteurs  $\caR^{d_i}_{\frqqq_{i}}$ de $L$ à $G$ le translatent de $\rho_V=\rho_G-\rho_L$. Mais ce n'est pas la seule 
obstruction à l'égalité de ces deux termes, comme nous allons le voir.
On va donc introduire une  torsion $\varepsilon$ sur les représentations virtuelles stables de 
$L_{*,H}$ et de  ses formes intérieures $L_{k,H}$, et on va remplacer 
(\ref{LTransIC}) par l'expression obtenue de la même façon à partir de 
$\varepsilon ([X_{L_{1,H}}^{st}])$, à savoir 
 \begin{equation} \label{LTransIC2}
\Trans_H^G\left(\sum_k  e(L_{k,H})   \;  \caR^{d_k}_{\frqqq_{k,H, L_{k,H},H}}\left (  \varepsilon \left( [ X_{L_{k,H}}^{st}] \right)\right) \right).
\end{equation}
On aura alors 
\begin{prop} \label{propquicompte} Avec les notations qui précèdent 
\[\sum_i  e(L_i \; )  \caR^{d_i}_{\frqqq_{i},L_i,G}  \left(\Trans_{L_{*,H}}^{L_i} \left( [X_{L_{*,H}}^{st}] \right)   \right)
= \Trans_H^G\left(\sum_k  e(L_{k,H})   \;  \caR^{d_k}_{\frqqq_{k,H, L_{k,H},H}}\left (  \varepsilon \left( [ X_{L_{k,H}}^{st}] \right)\right) \right).
\]
\end{prop}

Cette torsion $\varepsilon$ est définie de la manière suivante. Tout d'abord, il suffit par linéarité ({\sl cf.} Rmq. \ref{rmqtemp}, -3)
de la définir sur les représentations
virtuelles stables  (\ref{StablePsP2}). Considérons donc un paramètre de Langlands 
$\phi_{L_H}$ du groupe $L_{*,H}$. 
Une représentation standard $X$ pour ce paquet  
s'écrit  
\[X=X_u\boxtimes X_1 \boxtimes X_2\]
où $X_u$, $X_1$ et $X_2$ sont des représentations standards  bien déterminés de
$U_c$, $H'_1$ et $H_2$ respectivement. 

Définissons la représentation standard
\[\varepsilon(X)=  (X_u\otimes \varepsilon_U) \boxtimes X_1\boxtimes \varepsilon_2(X_2)\]
de la manière suivante. Sur le facteur unitaire, $\varepsilon_U$ est un caractère, qui va corriger le caractère infinitésimal de la manière voulue.
 Comme les groupes unitaires ont leurs sous-groupes de Cartan connexes,  le caractère $\varepsilon_U$ est 
 déterminé par sa différentielle   $d\varepsilon_U$, et celle-ci vaut  $(\rho_G-\rho_L)-(\rho_H-\rho_{L_H})$.
 Selon les cas, on obtient : 

 \medskip 

 Cas {\bf A1}, {\bf B}, {\bf C},  {\bf D} :  $d\varepsilon_U=(b,b,\ldots,b)$.  Cas {\bf A2} :   $d\varepsilon_U=(a+1,\ldots,a+1)$.

\medskip

La description de  $\varepsilon_2(X_2)$ est plus subtile. Lorsque $c$ est pair,
$\varepsilon_2(X_2)=X_2$. Supposons maintenant que $c$ soit impair.
 La représentation standard $X_2$ est une induite parabolique, obtenue
à partir d'une série discrète  $X_{M}$
d'un sous-groupe de Levi cuspidal de $H_2$, disons $M_{H_2}$, isomorphe à un produit de la forme
\[M_{H_2}\simeq   (\bbR^\times)^s  \times \GL_{2}(\bbR)^m\times H''_2 \] 
où $H''_2$ est un groupe classique quasi-déployé  de même type que $H_2$. La représentation $\varepsilon_2(X_2)$ est alors la représentation
standard obtenue par induction parabolique  à partir du même  sous-groupe de Levi cuspidal  $M_{H_2}$, 
et de la représentation $X_M$, mais où sur les facteurs $\bbR^\times$, on a tensorisé par le caractère $\sgn_{\bbR^\times}$.

Cette opération sur les représentations standard du facteur $H_2$ s'interprète facilement lorsque 
$H_2$ est un groupe spécial orthogonal.  Si $G$ est un groupe spécial orthogonal  non compact, il admet deux composantes connexes. Notons 
$\sgn_G$ le caractère quadratique de $G$ valant $-1$ sur la composante connexe non neutre.
Si $G$ est un groupe orthogonal compact,  $\sgn_G$  est le  caractère trivial par convention.
On a alors $\varepsilon(X_2)= X_2\otimes \sgn_{H_2}$. Il est alors clair que lorsque $X_2$ décrit un pseudo-paquet pour $H_2$, 
il en est de même pour $\varepsilon_2(X_2)$. De même, lorsque $X$ décrit un pseudo-paquet pour $L_{*,H}$, 
il en est de même pour $\varepsilon(X)$. Cette propriété est encore vraie si $H_2$ est un groupe symplectique, bien 
que l'opération effectuée sur les représentations standard n'ait pas d'interprétation aussi simple que la tensorisation
par un caractère quadratique. On obtient un autre paquet avec même caractère infinitésimal.

On adapte facilement cette construction aux autres formes réelles $\mathbf  L_{k,H}$. Le facteur $H_2$ est commun à tous les 
 $\mathbf  L_{k,H}$ et la torsion sur ce facteur reste la même. Le facteur groupe unitaire varie parmi les formes intérieures
 pures $\U(k,c-k)$ de $\U_c$ et le caractère de torsion $\varepsilon_U$ admet la même description dans tous les cas. 
On a alors 
\begin{align} \label{LTransIC2bis}
&\Trans_H^G\left(\sum_k  e(L_{k,H})   \;  \caR^{d_k}_{\frqqq_{k,H, L_{k,H},H}} \left(  \varepsilon\left([ X_{L_{k,H}}^{st}] \right) \right) \right)
\\ \nonumber 
=&\Trans_H^G\left( \left( \sum_k   e(L_{k,H_1})   \; \caR^{d_k}_{\frqqq_{k,H_1}, L_{k,H_1},H_1} \left(( [X_{\U(k,c-k)}^{st}] \otimes \varepsilon_U)
  \boxtimes [X_{H'_k}^{st}]\right)\
  \right) \boxtimes \varepsilon_2 [X_{H_2}^{st}] \right) .  \end{align}

Ceci termine la description de $\varepsilon$, et nous pouvons maintenant commencer la démonstration de 
la proposition \ref{propquicompte}.

\bigskip 

\noindent  \underline{\sl Démonstration de la proposition \ref{propquicompte}}.
  Par linéarité ({\sl cf.} Rmq. \ref{rmqtemp}, -3), il suffit de considérer le cas où $[X_{L_{*,H}}^{st}] =[\pi(\phi_L, L_{*,H} ) ]$ pour un paramètre de 
  Langlands $\phi_{L_H}$. 
Plaçons nous tout d'abord dans le cas où $\phi_{L_H}$ est un paramètre discret.

  On a alors 
  \[  [X_{L_{*,H}}^{st}] = [\pi(\phi_{L_H}, L_{*,H} )]= \sum_{\pi \in \Pi(\phi_{L_H}, L_{*,H}) } \pi,  \]
  et 
   \[ [X_{L_{k,H}}^{st}] = e(L_{k,H})  [\pi(\phi_{L_H}, L_{k,H})]=  e(L_{k,H})  \;  \sum_{\pi\in \Pi(\phi_{L_H}, L_{k,H}) }\pi .\]

Occupons nous d'abord du membre de droite dans l'égalité de la proposition \ref{propquicompte} :
\begin{equation*}
\Trans_H^G\left(\sum_k   e(L_{k,H})   \;  \caR^{d_k}_{\frqqq_{k,H}} ( \varepsilon ([ X_{L_{k,H}}^{st}]  ) ) \right)
=\Trans_H^G\left(\sum_k    \; 
 \caR^{d_k}_{\frqqq_{k,H}, L_{k,H},H}
\left( \sum_{\pi \in \Pi(\phi_L, L_{k,H}) } \varepsilon \left (\pi  \right)  \right) \right).  \end{equation*}

 La paramétrisation des séries discrètes de la section
\ref{seriesdiscretes}, et les considérations de la section 
\ref{StabIC}, appliquées ici à $L_{*,H}$ et ses formes intérieures pures  $L_{t,H}$ 
(voir la fin de la section \ref{cLevicla}, adaptée au groupe endoscopique $H$, et plus précisément, à sa composante $H_1$),
 montrent que l'on a
\[ \sum_k    \; 
 \caR^{d_k}_{\frqqq_{k,H}, L_{k,H},H}
\left( \sum_{\pi \in \Pi(\phi_L, L_{k,H}) } \varepsilon (\pi ) \right) = \sum_{t\in T_*[2]}
 \caR^{d_t}_{\frqqq_{*,H}, L_{t,H},H} \left(  \varepsilon \left( \pi (\phi_{L_H},t) \right)  \right) .
\]

Le groupe $L_{t,H}$ est un produit d'un facteur unitaire et de deux  facteurs  groupes classiques, $H'_1$ et $H_2$
comme précédemment. Une série discrète
$\pi=\pi(\phi_{L_H},t)$ de $L_{t,H}$ se décompose en un produit tensoriel extérieur
$\pi=\pi_u \boxtimes \pi_1 \boxtimes \pi_2$ où $\pi_u$ est une série discrète du groupe unitaire et $\pi_1$ (resp. $\pi_2)$
 une série discrète de $H'_1$ (resp. $H_2$).
 Or, pour une série discrète $\pi_2$ du facteur classique 
$H_2$, on a $\varepsilon_2(\pi_2)= \pi_2$.
Ainsi $\varepsilon(\pi)= (\pi_u\otimes \varepsilon_U)\boxtimes \pi_1\boxtimes \pi_2$.
Le produit tensoriel par un caractère unitaire $\varepsilon_U$  préserve les paquets de séries discrètes pour le groupe unitaire. 
Lorsque $t$ décrit $ T_*[2]$, les  $ \varepsilon(\pi(\phi_{L_H},t)) $ décrivent donc un autre super-paquet
dont le paramètre est obtenu par torsion avec un $1$-cocycle $\mathbf a$ de $W_\bbR$ dans le centre 
du groupe dual de $L_{*,H}$, notons ce nouveau paramètre $\mathbf a\phi_{L_H}$ (et remarquons que la torsion est uniquement sur 
le facteur dual du groupe unitaire).
D'autre part, ceci n'affecte pas la paramétrisation interne des paquets, dans le sens où 
\[\pi(\mathbf a\phi_{L_H},t)=\varepsilon(\pi(\phi_{L_H},t)). \]
On obtient donc, en utilisant (\ref{ICsd}), 
\begin{align*} &\sum_k    \; 
 \caR^{d_k}_{\frqqq_{k,H}, L_{k,H},H}
\left( \sum_{\pi\in \Pi(\phi_L, L_{k,H}) } \varepsilon(\pi ) \right) = \sum_{t\in T_*[2]}
 \caR^{d_k}_{\frqqq_{*,H}, L_{t,H},H} \left( \pi(\mathbf a \phi_{L_H},t)  \right)
 \\&=
  \sum_{t\in T_*[2]}   \pi(\iota_{L_H,H}(\mathbf a \phi_{L_H}),t)  =  [\pi (\iota_{L_H,H} (\mathbf a \phi_{L_H})) ,H  )] .
\end{align*}
Le transfert endoscopique d'un paquet de séries discrètes est bien connu  d'après Shelstad \cite{She82}, \cite{She08III} et \cite{Adsd}.
Supposons tout d'abord que l'on se trouve dans le cas $G$-régulier, c'est-à-dire que
$ \iota_{H,G}(\iota_{L_H,H}(\mathbf a \phi_{L_H}))$ est un paramètre discret pour $G$. 
Ainsi, le    membre de droite dans l'égalité de la proposition \ref{propquicompte} devient,
\begin{equation}\label{LTransIC3} 
\Trans_H^G (  [\pi( \iota_{L_H,H}(\mathbf a \phi_{L_H})),H )] )= 
\sum_{t\in T_*[2],\,  [t]\sim_G[1]} \bil{x}{t}_{G}\;  \pi(\iota_{H,G}(\iota_{L_H,H}(\mathbf a \phi_{L_H})),t).   \end{equation}
Le pairing est celui de (\ref{pairingAdams}), entre  $ T_*[2]$ et $A(\iota_{L_H,H}(\mathbf a \phi_{L_H}))$ 
où $x$ désigne aussi l'image de l'élément $x$ de la donnée endoscopique dans 
$A(\iota_{L_H,H}(\mathbf a \phi_{L_H}))$.

Examinons maintenant le  membre de gauche dans l'égalité de la proposition \ref{propquicompte}. On a 
\begin{equation}\label{TLH} \Trans_{L_{*,H}}^{L_i} \left( [X_{L_{*,H}}^{st}]  \right)= e(L_i) \sum_{t\in T_*[2],\,  [t]\sim_L[t_i]}
 \bil{x}{t}_{L}\;  \pi(\iota_{L_H,L}( \phi_{L_H}),t).  \end{equation}
Le pairing est celui de (\ref{pairingAdams}), entre  $ T_*[2]$ et $A(\iota_{L_H,L}( \phi_{L_H}))$.
Si on applique à ceci le foncteur $\caR^{d_i}_{\frqqq_i,L_i,G}$, on obtient d'après (\ref{ICsd}),
\[  e(L_i)  \sum_{t\in T_*[2],\,  [t]\sim_L[t_i]}
 \bil{x}{t}_{L}\; \pi(\iota_{L,G}(\iota_{L_H,L}( \phi_{L_H})),t).   \]
 En sommant sur les formes intérieures $L_i$, on obtient finalement l'expression suivante pour
 le membre de gauche dans l'égalité de la proposition \ref{propquicompte}  : 
 \begin{equation}\label{sumisum}  \sum_i  \sum_{t\in T_*[2],\,  [t]\sim_L[t_i]}
 \bil{x}{t}_{L}\; \pi(\iota_{L,G}(\iota_{L_H,L}(\mathbf a \phi_{L_H})),t)= \sum_{t\in T_*[2]\sim_G[1]}
 \bil{x}{t}_{L}\; \pi(\iota_{L,G}(\iota_{L_H,L}( \phi_{L_H})),t).   \end{equation}
 La comparaison avec (\ref{LTransIC3}) est alors aisée car $\iota_{H,G}(\iota_{L_H,H}(\mathbf a \phi_{L_H}))
 =\iota_{L,G}(\iota_{L_H,L}( \phi_{L_H}))$
 (ce sont des paramètres discrets de même caractère infinitésimal), et 
 $A(\iota_{H,G}(\iota_{L_H,H}(\mathbf a \phi_{L_H})))$  est égal à $A(\iota_{L_H,L}( \phi_{L_H}))$ 
(pour les paramètres discrets $\phi$  des groupes  que nous considérons ci-dessus, le groupe 
$A(\phi)$ est de manière explicite isomorphe à $(\bbZ/2\bbZ)^n$, où $n$ est le rang du groupe, comme on le voit en (\ref{pairingAdams})).
 On obtient donc l'égalité entre (\ref{RTransIC}) et (\ref{LTransIC2})
dans le cas d'un paramètre discret et $G$ régulier. 
 
\bigskip

\begin{rmq} On a supposé que l'on était dans le good range pour les inductions cohomologiques, en particulier 
si $\phi_{L_H}$ est un paramètre discret, il en est de même de $\iota_{L_H,H}(\phi_{L_H})$. \end{rmq}

\bigskip 

Passons maintenant au cas d'un paramètre discret $\phi_{L_H}$, mais sans l'hypothèse de $G$-régularité:
 $ \iota_{H,G}(\iota_{L_H,H}(\mathbf a \phi_{L_H}))$ est un paramètre 
de limites de séries discrètes pour $G$ et d'après la remarque \ref{lsd}, le membre de droite de 
(\ref{LTransIC3}) a toujours un sens, mais certaines représentations 
 $\pi(\iota_{H,G}(\iota_{L_H,H}(\mathbf a \phi_{L_H})),t)$ peuvent maintenant être nulles, et le pairing 
 $\bil{x}{t}_{G}$ à un sens exactement lorsque  $\pi(\iota_{H,G}(\iota_{L_H,H}(\mathbf a \phi_{L_H})),t)\neq 0$.
  La formule de transfert (\ref{LTransIC3}) est toujours valable  
et se démontre par des arguments de continuation cohérente ({\sl cf.} \cite{She82}, \cite{She08III}).
Il en est de même pour le terme de droite dans  (\ref{sumisum}), et la comparaison entre (\ref{sumisum}) 
est (\ref{LTransIC3}) toujours valide.

\bigskip 

Passons maintenant au cas d'un paramètre  $\phi_{L_H}$ quelconque.
On se place dans le cas {\bf A1} de la liste, mais les autres cas se traitent de la même manière, il s'agit juste d'éviter 
d'introduire des notations trop lourdes. 
  On a donc  $\mathbf G=\Sp_{2n}$, $\mathbf H=\Sp_{2a}\times \SO_{2b}^\alpha =\mathbf H_1\times \mathbf H_2$, 
   \[ \mathbf L_{*,H}\simeq \U_c  \times\Sp_{2(a-c)}\times \SO_{2b}^\alpha = \U_c\times \mathbf H'_1\times \mathbf H_2.\]
Comme $H_2$ est un groupe orthogonal, la  torsion  $\varepsilon$ est dans ce cas
donnée par tensorisation par le caractère  $\varepsilon_{U}\boxtimes  \Triv_{G'} \boxtimes \sgn_{H_2}^c$, où 
$\varepsilon_{U}$ est le caractère de $\U_c$  de caractère infinitésimal $(b,b,\cdots, b)$.

 On part donc dans ce cas d'un paquet de $L_{*,H}$ associé à un sous-groupe de Cartan
$D=TA\simeq \U(1)^{r_1} \times (\bbC^\times)^{m_1}\times  \U(1)^{r_2} \times (\bbC^\times)^{m_2}\times (\bbR^\times)^{s_2}\times  
\U(1)^{r_3} \times (\bbC^\times)^{m_3}\times (\bbR^\times)^{s_3}$.
Ici le facteur  $\U(1)^{r_1} \times (\bbC^\times)^{m_1}$ est un sous-groupe de Cartan du facteur groupe unitaire (donc $r_1+2m_1=c$), 
le facteur   $\U(1)^{r_2} \times (\bbC^\times)^{m_2}\times (\bbR^\times)^{s_2}$
est un sous-groupe de Cartan du facteur $H'_1$ (donc $r_2+2m_2+s_2=a-c$)
et le facteur   $\U(1)^{r_3} \times (\bbC^\times)^{m_3}\times (\bbR^\times)^{s_3}$
est un sous-groupe de Cartan du facteur $H_2$ (donc $r_3+2m_3+s_3=b$). 
Le pseudo-paquet    $\Pi^\sharp (\phi_{L_H}, L_{*,H})$  est constitué d'induites paraboliques à partir du sous-groupe :
\[  M_{L_{*,H}}\simeq   \left( \U_{r_1} \times   (\bbC^\times)^{m_1} \right)\times 
\left(H'_{r_2}\times \GL_2(\bbR)^{m_2}\times   (\bbR^\times)^{s_2} \right)\times 
\left(H''_{r_3}\times \GL_2(\bbR)^{m_3}\times   (\bbR^\times)^{s_3} \right),\]
où on induit des séries discrètes relatives. Ici $\mathbf H'_{r_2}=\Sp_{2r_2}$ et 
$\mathbf H''_{r_3}=\SO_{2r_3}^\alpha$.
Le paquet correspondant  $\Pi^\sharp(\phi_H,H)$        de  $H$ est lui constitué 
 d'induites paraboliques à partir du sous-groupe :
\[  M_H \simeq   
\left(H'_{r_1+ r_2}\times \GL_2(\bbR)^{m_1+m_2}\times   (\bbR^\times)^{s_2} \right)\times 
\left(H''_{r_3}\times \GL_2(\bbR)^{m_3}\times   (\bbR^\times)^{s_3} \right),\]
avec $\mathbf H'_{r_1+r_2}=\Sp_{2(r_1+r_2)}$.
Le paquet correspondant $\Pi^\sharp(\phi_L,L_1)$     de $L_1$ est constitué 
 d'induites paraboliques à partir du sous-groupe :
\[  M_L\simeq     \left( \U_{r_1} \times   (\bbC^\times)^{m_1} \right)\times 
\left(G'_{r_2+ r_3}\times \GL_2(\bbR)^{m_2+m_3}\times   (\bbR^\times)^{s_2+s_3} \right)
\]
avec $\mathbf G'_{r_2+r_3}=\Sp_{2(r_2+r_3)}$.
Le paquet correspondant  $\Pi^\sharp(\phi_G,G)$  de  $G$ est lui constitué 
 d'induites paraboliques à partir du sous-groupe :
\[  M\simeq    G_{r_1+ r_2+r_3}\times \GL_2(\bbR)^{m_1+m_2+m_3}\times   (\bbR^\times)^{s_2+s_3} \]
avec $\mathbf G_{r_1+r_2+r_3}=\Sp_{2(r_1+r_2+r_3)}$.
Résumons la situation dans le diagramme suivant, 

\bigskip 

\begin{tikzcd}[row sep=scriptsize, column sep=scriptsize]
& M_{L_H}  \arrow[dl, "IP"]\arrow[rrrr,"IC"] \arrow[dddd, "Tr" ] & & && M_H \arrow[dl,"IP"] \arrow[dddd,"Tr"] \\
L_H \arrow[rrrr,"IC"] \arrow[dddd,"Tr"] & & & & H \arrow[dddd, "Tr"]\\
&&&&\\
&&&&\\
& M_L \arrow[dl,"IP"] \arrow[rrrr,"IC"] & &&&  M \arrow[dl,"IP"] \\
L \arrow[rrrr,"IC"] & & && G \\
\end{tikzcd}

où $IP$ indique une induction parabolique, $IC$ une induction cohomologique, 
$Tr$ un transfert endoscopique.
La face de devant est celle qui nous  intéresse. La face du fond est celle traitée dans le cas des séries discrètes, à des facteurs
$\bbR^\times$ et $\GL_2(\bbR)$ près, qui ne jouent aucun rôle.
Les faces du haut et du bas ne commutent pas, il faut introduire dans chaque cas une torsion pour que ça commute
(défaut de commutation entre induction parabolique et cohomologique, étudiée dans l'appendice).
Par contre, les faces de gauche et de droite commutent (le transfert et l'induction parabolique commutent). 
On va vérifier que le défaut de commutation des faces du haut et du bas induit bien le défaut de commutation de la face avant annoncée.

D'après l'appendice, le défaut de commutation de la face du bas est $\gamma_{\frn}\gamma_{\frn_L}^*$.
Cette torsion est triviale lorsque $c$ est pair et consiste à  tensoriser par   le caractère $\sgn_{\bbR^\times}$ 
sur chaque facteur $\bbR^\times$ de $D$  lorsque $c$ est impair.
Pour la face du haut, le défaut de commutation est analogue,  en appliquant le calcul de la torsion aux facteurs $H_1$ et $H_2$ du 
groupe endoscopique $H$.
Pour rendre la face avant commutative, on voit donc qu'il faut composer sur la flèche entre 
$L_H$ et $H$ par  la  torsion   $\varepsilon_2$ du facteur $H_2$. 
 Sur le facteur $H_1$, les torsions  des faces du bas et du haut sont les mêmes et se compensent.
\qed

\subsection{Normalisation de l'induction cohomologique}\label{SecNormIC}

Nous avons calculé dans la section précédente le défaut de commutativité entre transfert et induction cohomologique, 
où les foncteurs d'induction cohomologique sont ceux définis par Vogan et Zuckerman. D'autre normalisations des foncteurs
d'induction cohomologique sont possibles, par exemple dans \cite{KV}, Chapter 11, une normalisation des foncteurs
d'induction cohomologique est proposée, qui a l'avantage de respecter le caractère infinitésimal, mais a l'inconvénient
de faire intervenir des revêtements à deux feuillets des groupes réels. Pour nous, il n'est pas avantageux de faire intervenir de tels
groupes, qui ne sont plus des groupes classiques.
Ce que nous aimerions dans l'idéal,  c'est avoir une version de l'induction cohomologique qui commute au transfert endoscopique et préserve le caractère infinitésimal.
 Les résultats de la section précédente et les considérations ci-dessus montrent que ce n'est pas possible en général, mais nous allons 
 néanmoins renormaliser nos foncteurs d'induction cohomologique pour nous rapprocher le plus possible de cet idéal.

Les foncteurs   d'induction cohomologique qui apparaissent dans la section précédente sont, avec les notations de cette section,  les 
 $\caR^{d_i}_{\frqqq_{i},L_i,G} $ et les  $\caR^{d_k}_{\frqqq_{1,H}, L_{k,H},H}$.
 Commençons par les $\caR^{d_i}_{\frqqq_{i},L_i,G} $. Ils translatent le caractère infinitésimal de $\rho_G-\rho_L$. 
 Calculons ceci dans les coordonnées usuelles : 
 
 cas {\bf A}:  $\rho_G=(n,n-1,\ldots ,1)$ et $\rho_L=(\frac{c-1}{2}, \frac{c-3}{2}, \ldots , -\frac{c-1}{2}, n-c, n-c-1,\ldots, 1)$, d'où 
$\rho_G-\rho_L=(n-\frac{c-1}{2}, \ldots , n-\frac{c-1}{2}, 0,\ldots ,0) $.

 cas {\bf B}:  $\rho_G=(\frac{2n-1}{2},\ldots ,\frac{1}{2})$ et $\rho_L=(\frac{c-1}{2}, \frac{c-3}{2}, \ldots , -\frac{c-1}{2}, \frac{2(n-c)-1}{2},    \ldots,  \frac{1}{2})$, d'où 
$\rho_G-\rho_L=(n-\frac{c}{2}, \ldots , n-\frac{c}{2}, 0,\ldots ,0) $.

cas {\bf C} et {\bf D}:  $\rho_G=(n-1,\ldots, 0)$ et $\rho_L=(\frac{c-1}{2}, \frac{c-3}{2}, \ldots , -\frac{c-1}{2}, n-c-1,n-c-1,\ldots, 0)$, d'où 
$\rho_G-\rho_L=(n-1-\frac{c-1}{2}, \ldots , n-1-\frac{c-1}{2}, 0,\ldots ,0) $.

Rappelons comment sont normalisés les foncteurs de Vogan-Zuckerman. On a une version non normalisée 
des foncteurs d'induction cohomologique, que nous allons noter  ${}^n\caR^{d_i}_{\frqqq_{i},L_i,G} $, et qui sont obtenus comme foncteurs dérivés
d'un foncteur   ${}^n\caR^{0}_{\frqqq_{i},L_i,G} $, lui-même obtenu comme adjoint d'un foncteur d'oubli
 (des représentations de $G$ vers les représentations de $L_i$).   Si $Z$ est une représentation de $L_i$, 
\[ \caR^{d_i}_{\frqqq_{i},L_i,G}(Z) ={}^n\caR^{d_i}_{\frqqq_{i},L_i,G} (Z\otimes \textstyle \bigwedge ^{top} \frv_i) \]
où $\frv_i$ est le radical nilpotent de $\frqqq_i$, et  $\bigwedge ^{top} \frv_i$
son algèbre extérieure en degré maximal $\dim \frv_i$ est une représentation de dimension 1 de $L_i$
de caractère infinitésimal $2(\rho_G-\rho_L)$.

Pour avoir un foncteur qui préserve le caractère infinitésimal, il faudrait remplacer $\bigwedge ^{top} \frv_i$
par une représentation de dimension 1 de  $L_i$ de caractère infinitésimal $\rho_G-\rho_L$.
Or comme nous le voyons sur les formules ci-dessus, $\rho_G-\rho_L$ ne vit que sur le facteur 
unitaire de $L_i$ (et l'on remarque que c'est dans le centre), et pour un groupe unitaire de rang $c$, 
pour que $(e,e,\ldots, e)$ soit le caractère infinitésimal d'une représentation de dimension $1$, il faut que $e$ soit entier.

\begin{defi}\label{def:xiu}
Sur le facteur unitaire du groupe $L_i$, définissons  le caractère $\xi_u$, dont la différentielle est 

 cas {\bf A}: $ (n-\lfloor \frac{c-1}{2}\rfloor, \ldots , n-\lfloor \frac{c-1}{2}\rfloor)$. Ceci est égal à   $\rho_G-\rho_L$ si $c$ est impair, 
 et $\rho_G-\rho_L+\left(  \frac{1}{2}, \ldots  ,\frac{1}{2}, 0,\ldots ,0\right)$ si $c$ est pair.

 cas {\bf B}:  $(n-\lfloor \frac{c}{2} \rfloor , \ldots , n-\lfloor \frac{c}{2}\rfloor ) $. Ceci vaut $\rho_G-\rho_L$ si $c$ est pair, 
 et $\rho_G-\rho_L+\left(  \frac{1}{2}, \ldots  ,\frac{1}{2}, 0,\ldots ,0\right)$ si $c$ est impair.

cas {\bf C} et {\bf D}:  $(n-1-\lfloor  \frac{c-1}{2}\rfloor , \ldots , n-1-\lfloor \frac{c-1}{2}\rfloor ) $. Ceci vaut $\rho_G-\rho_L$ si $c$ est impair, 
 et $\rho_G-\rho_L+\left(  \frac{1}{2}, \ldots  ,\frac{1}{2}, 0,\ldots ,0\right)$ si $c$ est pair.
\end{defi}

L'autre facteur de $L_i$ est une forme intérieure pure  $G'_i$ d'un groupe 
$\mathbf G'$  de même type que $\mathbf G$. Si c'est un groupe orthogonal, tordons par le caractère 
quadratique $\xi'=\sgn_{G_i'}^c$ introduit dans la section précédente. 
Pour les groupes symplectiques, on ne peut rien faire de plus, et l'on définit alors $\xi'$ comme le caractère trivial.

Nous définissons les foncteurs d'induction cohomologiques normalisés de la manière suivante. 

\begin{defi}Si 
$Z=Z_u\boxtimes Z'$, où $Z_u$ est une représentation du facteur unitaire de $L_i$ et 
$Z'$ une représentation du facteur $G_i'$, posons
\begin{equation}\label{NormIC}  
\widetilde \caR^{d_i}_{\frqqq_{i},L_i,G}(Z) ={}^n\caR^{d_i}_{\frqqq_{i},L_i,G} ((Z_u\otimes \xi_u)\boxtimes (Z'\otimes \xi')). \end{equation}
 \end{defi}
 
 \begin{rmq}\label{bcshift}
 La translation sur le caractère infinitésimal est donc soit nulle, soit $(1/2,\ldots,1/2)$ sur le facteur unitaire.
  Le choix de la normalisation sur le facteur unitaire 
 est dicté par le fait que la translation sur le caractère infinitésimal qui reste
 correspond à celle des changements de base principal et antiprincipal ${}^L\U_c\rightarrow {}^L\GL_c(\bbC)$. 
\end{rmq}

On renormalise aussi de la même façon les foncteurs $\caR^{d_k}_{\frqqq_{1,H}, L_{k,H},H}$.
Pour les groupes orthogonaux  $L_{k,H}$ est un produit de trois facteurs, un groupe unitaire, et deux groupes orthogonaux  $H'_k$
et $H_2$. Pour une représentation $Z=Z_u\boxtimes Z_1\boxtimes Z_2$ de $L_{k,H}$, on définit 
\[  \widetilde \caR^{d_k}_{\frqqq_{1,H}, L_{k,H},H}(Z)=
\left( {}^n\caR^{d_k}_{\frqqq_{1,H}, L_{k,H_1},H_1}( Z_u\otimes \xi_u^H) \boxtimes (Z_1\otimes \xi_1)
\right) \boxtimes Z_2 .\]
Ici $\xi_u^H$ est défini de manière analogue à $\xi_u$, comme dans la définition \ref{def:xiu}, mais relativement à $H_1$, c'est-à-dire que 
$\xi_u$ est de différentielle $\rho_H-\rho_{L_H}$ ou bien $\rho_H-\rho_{L_H}+ \left(  \frac{1}{2}, \ldots  ,\frac{1}{2}, 0,\ldots ,0\right)$, selon les cas.
Dans le cas des groupes symplectiques, on prend la même formule si $H_1$ est orthogonal, et si  $H_1$ est symplectique, on ne met
 pas de torsion sur $Z_1$.
 On peut donc reformuler les résultats de la section précédentes, avec les foncteurs d'induction cohomologique normalisés.
 
 \begin{thm} \label{ThmICetTrans}
 Avec les notations de la section précédente, pour les groupes symplectiques,
 \begin{equation}\Trans_H^G\left(\sum_k  e(L_{k,H})   \; \widetilde \caR^{d_k}_{\frqqq_{k,H, L_{k,H},H}}( \tilde  \varepsilon
 \left([ X_{L_{k,H}}^{st}])\right) \right)
=
\sum_i  e(L_i \; ) \widetilde  \caR^{d_i}_{\frqqq_{i},L_i,G}  \left(\Trans_{L_{*,H}}^{L_i} \left( [X_{L_{*,H}}^{st}] \right)   \right).
\end{equation}
 et pour les groupes orthogonaux
  \begin{equation}\label{gto}
  \Trans_H^G\left( \sum_k  e(L_{k,H})   \; \widetilde \caR^{d_k}_{\frqqq_{k,H, L_{k,H},H}}   \left( [ X_{L_{k,H}}^{st}]  \right) \right)
=
\sum_i  e(L_i) \;  \widetilde  \caR^{d_i}_{\frqqq_{i},L_i,G}  \left(\Trans_{L_{*,H}}^{L_i} \left( [X_{L_{*,H}}^{st}] \right)   \right).
\end{equation}
 
 \end{thm}
 
 Ici, la correction $\tilde \varepsilon$ est analogue à la correction $\ \varepsilon$  de la section précédente. 
 On la détermine en donnant son effet sur une représentation standard de $L_{k,H}$ de la forme
 \[ X=X_u\boxtimes X_1\boxtimes X_2 .\]
Il n'y a plus d'effet sur $X_u$. Il n'y a pas d'effet sur $X_1$ si $H_1$ est symplectique, et si $H_1$ est orthogonal
et on tensorise par $\sgn_{H'_1}^c$. On applique à $X_2$ la même opération
$\varepsilon_2$ que celle que l'on avait définie
pour $\varepsilon$. On écrit :
 \[ \tilde \varepsilon(X)=X_u\boxtimes \tilde \varepsilon_1(X_1)\boxtimes \tilde \varepsilon_2 (X_2) .\]

\begin{rmq}
Dans les cas {\bf B},  {\bf C} et  {\bf D}  (groupes orthogonaux) nous avons
 un énoncé du type : le transfert commute à l'induction cohomologique (normalisée).
\end{rmq}

Reformulons encore le théorème. On écrit : 
\[[ X_{L_{k,H}}^{st}]=  [X_{\U(k,c-k)}^{st}]\boxtimes [ X_{H'_{1,k}}^{st}]  \boxtimes [ X_{H_2}^{st}], \]
et dans le cas des groupes symplectiques 
\[ \tilde  \varepsilon \left([X_{L_{k,H}}^{st}]\right)
 =  [X_{\U(k,c-k)}^{st}]\boxtimes\tilde  \varepsilon_1  ( [ X_{H'_{1,k}}^{st}] )
  \boxtimes \tilde  \varepsilon_2 \left( [ X_{H_2}^{st}]\right)  . \]

On a alors 
\[ \widetilde \caR^{d_k}_{\frqqq_{k,H, L_{k,H},H}} \left( \tilde  \varepsilon
 \left( [ X_{L_{k,H}}^{st}]\right)\right) = \widetilde \caR^{d_{1,k}}_{\frqqq_{k,H_1}, L_{k,H_1},H_1}
 \left( [X_{\U(k,c-k)}^{st}]\boxtimes  \tilde \varepsilon_1( [ X_{H'_{1,k}}^{st} ]) \right)    \boxtimes   \varepsilon_2 \left( [ X_{H_2}^{st}]\right),    \]

\[  \Trans_{L_{1,H}}^{L_i} \left( [X_{L_{1,H}}^{st}] \right)  = S_i(x_d) [X_{\U(i,c-i)}^{st}]\boxtimes  \Trans_{H'_1\times H_2}^{G'_i} 
\left(  [ X_{H'_{1,k}}^{st}]  \boxtimes [ X_{H_2}^{st}]\right)   \]

La formule du théorème est dans le cas des groupes symplectiques :
 \begin{align}\label{Formquisert} & \Trans_H^G \left(  \sum_k  e(L_{k,H_1}) 
   \; \widetilde \caR^{d_{1,k}}_{\frqqq_{k,H_1}, L_{k,H_1},H_1} \left(   [X_{\U(k,c-k)}^{st}]\boxtimes  \tilde \varepsilon_1 ([ X_{H'_{1,k}}^{st}] )  \right)
 \boxtimes  \tilde \varepsilon_2  \left( [ X_{H_2}^{st}]\right) \right)\\
\nonumber  &=\sum_i  e(L_i )   S_i(x_d)  \; \widetilde  \caR^{d_i}_{\frqqq_i,L_i,G}  \left( [X_{\U(i,c-i)}^{st}]\boxtimes  \Trans_{H'_1\times H_2}^{G'_i} 
\left(  [ X_{H'_{1,k}}^{st}]  \boxtimes [ X_{H_2}^{st}]\right)  \right).
\end{align}

Pour les groupes orthogonaux, c'est la même formule sans les   $\tilde \varepsilon $.

\subsection{Transfert et induction cohomologique pour les formes intérieures pures $\SO(p,q)$}\label{TICFIP}

Voyons ce qui change lorsqu'on remplace le groupe $\mathbf G$ par une forme intérieure pure $\SO(p,q)$ d'un groupe orthogonal.
On a maintenant  un système de représentants $\{\mathbf Q_i=\mathbf L_i\mathbf V_i\}$ de $\Sigma_{{}^dQ}^{\SO(p,q)}$ 
  comme en (\ref{LiFi})  et  un transfert spectral $\Trans_H^{\SO(p,q)}$. 
  On a aussi les transferts spectraux (\ref{signtransfU}) de $L_{*,H}$ vers les $L_i$.
Dans la proposition \ref{propquicompte} ce qui change est donc : dans le terme de droite, $\Trans_H^{G}$
qui est remplacé par  $\Trans_H^{\SO(p,q)}$, et dans le terme de gauche
les $L_i$ de (\ref{Li})  sont remplacés par ceux de (\ref{LiFi}).  On va toujours avoir l'égalité entre les  deux termes dans la 
proposition \ref{propquicompte}, à condition 
de multiplier un des membres par le signe de Kottwitz $e(G)$, et la démonstration est identique.

On en déduit l'analogue de la formule (\ref{gto}) pour les formes intérieures des groupes orthogonaux :
\begin{equation}\label{gtofi}
 e(\SO(p,q))\;  \Trans_H^{\SO(p,q)}\left( \sum_k  e(L_{k,H})   \; \widetilde \caR^{d_k}_{\frqqq_{k,H, L_{k,H},H}}  
   \left( [ X_{L_{k,H}}^{st}]  \right) \right)
   \end{equation}
   \[
=
\sum_{(\mathbf Q_i=\mathbf L_i\mathbf V_i)\in\Sigma_{{}^dQ}^{\SO(p,q)}}  e(L_i) \;  \widetilde  \caR^{d_i}_{\frqqq_{i},L_i,G} 
 \left(\Trans_{L_{*,H}}^{L_i} \left( [X_{L_{*,H}}^{st}] \right)   \right),
\]
et sa variante, analogue de (\ref{Formquisert}) :
 \begin{align}\label{Formquisertfi} & e(\SO(p,q))\;  \Trans_H^{\SO(p,q)} \left(  \sum_k  e(L_{k,H_1}) 
   \; \widetilde \caR^{d_{1,k}}_{\frqqq_{k,H_1}, L_{k,H_1},H_1} \left(   [X_{\U(k,c-k)}^{st}]\boxtimes  [ X_{H'_{1,k}}^{st}]   \right)
 \boxtimes   [ X_{H_2}^{st}]  \right)\\
\nonumber  &=\sum_{(\mathbf Q_i=\mathbf L_i\mathbf V_i)\in\Sigma_{{}^dQ}^{\SO(p,q)}}  e(L_i )   S_i(x_d)  \; \widetilde  \caR^{d_i}_{\frqqq_i,L_i,G} 
 \left( [X_{\U(i,c-i)}^{st}]\boxtimes  \Trans_{H'_1\times H_2}^{G'_i} 
\left(  [ X_{H'_{1,k}}^{st}]  \boxtimes [ X_{H_2}^{st}]\right)  \right).
\end{align}

  \section{Ajout d'un bloc discret de bonne parité de paramètre suffisamment grand}\label{principal}
  
\subsection{Enoncé}  
  Soient $\mathbf{G}$ l'un des groupes classiques quasi-déployés de la section \ref{grcla} et
   $\St_G :\; {}^L G \rightarrow \GL_N(\bbC)$ la représentation standard de son $L$-groupe. 
   Soit $\psi_G : \; W_\bbR \times \SL_2(\bbC) \rightarrow {}^L G $ 
un paramètre  d'Arthur pour $G$.
 Notons $\psi=\St_G \circ \psi_G$ : c'est un  paramètre  d'Arthur
de $\GL_N(\bbR)$, auquel est  associé une représentation $\Pi_\psi^ \GL=\Pi_{\phi_\psi}^ \GL$  de ce groupe.

 On suppose que $\psi$ se décompose en 
\begin{equation}\label{psidapsi}
\psi=\psi_d \oplus \psi'=(V(0,t)\boxtimes R[c]) \oplus \psi'
\end{equation}
 où $\psi_d=V(0,t)\boxtimes R[c]$ est le paramètre
d'Arthur pour $G_{N_d}$, associé à la représentation 
$\Pi_{\psi_d}^\GL=\Speh(\delta\left(\frac{t}{2}, -\frac{t}{2}\right), c)$, $t>c-1$, 
et $\psi'$ est un paramètre d'Arthur pour $G_{N'}$ de caractère infinitésimal $\lambda'=(\lambda'_1, \ldots,
 \lambda'_{N'})$. On a ici $N_d=2c$ et $N=N_d+N'$.

  On suppose que l'on a 
\begin{equation}\label{hyplemmegen2}
\frac{t-(c-1)}{2} >\vert  \lambda'_i  \vert, \qquad  (1\leq  i \leq N')  
\end{equation}
ce qui permettra d'appliquer les résultats de \cite{AMR}.

Si $G$ est  symplectique, et si $t$ est pair, (et alors  $c$ est nécessairement impair)  alors  $V(0,t)\boxtimes R[c]$ 
est à valeurs dans un groupe orthogonal et non spécial orthogonal, et   il doit en être de même de $\psi'$.
Posons alors $\epsilon_\psi=\sgn_{W_\bbR}$ et 
$\psi'_\epsilon=\psi'\otimes \epsilon_\psi=\psi'\otimes \sgn_{W_\bbR}$. Alors $\psi'_\epsilon$ se factorise en 
$\psi'_\epsilon= \Std_{G'}\circ \psi_{G'}$, où $\mathbf G'$ est un groupe symplectique de rang  $n-c$.

Dans tous les autres cas, on prend $\epsilon_\psi$ égal au caractère trivial de $W_\bbR$ et    $\psi'_\epsilon=\psi'$, et l'on a une factorisation 
$\psi'_\epsilon= \Std_{G'}\circ \psi_{G'}$, où $\mathbf G'$ est un groupe de même type que 
$\mathbf G$  (sauf les types {\bf C} et {\bf D} qui sont échangés si $c$ est impair) et de rang $n':=n-c$.

Notre but dans cette section est de décrire explicitement le paquet d'Arthur pour $G$ et ses formes intérieures pures associé au paramètre $\psi_G$, 
à partir du paquet d'Arthur pour $G'$ (et ses formes intérieures pures) associé au paramètre $\psi_{G'}$.

Si  $G'_i$ est une forme intérieure pure de $G'$, soit 
${\pi}^A(\psi_G',G'_i)$ la représentation associée au paramètre $\psi_{G'}$. On suppose que c'est une représentation unitaire
de $G'_i \times  A( \psi_{G'})$ (ceci est établi  par Arthur si $G'_i$ est quasi-déployé, et un de nos buts ici est de voir
que cette propriété passe des $G'_i$ à  $\mathbf G$ et ses formes intérieures pures).
On écrit donc 
$\displaystyle {\pi}^A(\psi_G',G'_i)= \sum_{\eta' \in \widehat{A( \psi_{G'})}} \pi( \psi_{G'},\eta',G'_i)\boxtimes \eta'$  et soit 
$\Pi(\psi_{G'},G_i')$ le paquet d'Arthur associé à $\psi_{G'}$, c'est-à-dire la réunion des 
représentations irréductibles apparaissant dans les $\pi(\psi_{G'},\eta', G_i')$ pour $\eta' \in  \widehat{A( \psi_{G'})}$.
On note $\Pi^\GL_{\psi'_\epsilon}$ la représentation unitaire irréductible de $\GL_{N'}(\bbR)$
dont $\psi'_\epsilon$ est le paramètre d'Arthur. Comme $\psi'_\epsilon=\epsilon_{\psi}\otimes \psi'$, on a  
 $\Pi^\GL_{\psi'_\epsilon}=\varepsilon_{\psi'}\otimes \Pi^\GL_{\psi'}$, où  $\varepsilon_{\psi'}$
est le caractère $\Triv_{N'}$ ou $\sgn_{N'}$ de $\GL_{N'}(\bbR)$ selon que $\epsilon_{\psi}=\Triv_{W_\bbR}$ ou $\sgn_{W_\bbR}$.

Considérons le sous-groupe parabolique standard maximal ${}^dQ={}^dL{}^dV$ de $\widehat G$ tel que le $c$-Levi associé $\mathbf L_*$
 comme dans la section \ref{cLevicla}  vérifie  $\mathbf L_*\simeq \U_c\times \mathbf G'$; et rappelons que $\U_c$ est un groupe unitaire
  quasi-déployé  de rang $c$.

Rappelons que nous avons fixé en (\ref{Li})
un système de représentants $(\mathbf Q_i=\mathbf L_i \mathbf V_i)_{i=0,\ldots,c}$ de $\Sigma^G_{{}^dQ}$
avec \[  L_i\simeq \U(i,c-i)\times G'_i \] où $G'_i$  est une forme 
intérieure pure de $G'$, et de même en (\ref{LiFi}) 
pour les formes intérieures pures $\SO(p,q)$ des groupes orthogonaux, avec un système de représentants
 $(\mathbf Q_i=\mathbf L_i \mathbf V_i)_i$ de $\Sigma^{\SO(p,q)}_{{}^dQ}$.

On définit le caractère  $\xi_{i,t}$  de $\U(i,c-i)$ par le fait que sa différentielle est $\left( \lfloor \frac{t}{2}
\rfloor, \ldots , \lfloor \frac{t}{2}\rfloor \right)$.
On peut aussi le décrire à partir de $\psi_d=V(0,t)\boxtimes R[c]$ : 
on restreint $\psi_d$ à $W_\bbC\times \SL_2(\bbC)=\bbC^\times\times \SL_2(\bbC)$, ce qui nous donne un paramètre d'Arthur 
pour $\GL_{c}(\bbC)$.
Ce paramètre est obtenu par changement de base à partir d'un paramètre pour $\U_c$, le changement de base 
${}^L\U_c\rightarrow {}^L\GL_c(\bbC)$ étant principal si $t$ est pair, 
et anti-principal si $t$ est impair ({\sl cf.} remarque \ref{bcshift}).

Rappelons le caractère $S_i$ de $\bbZ/2\bbZ$ défini en (\ref{Si}). 
\begin{lemme}\label{LemmeApsi}
On a
\begin{equation}\label{ApsiG}A(\psi_G)= \bbZ/2\bbZ  \times A(\psi_{G'}),
\end{equation}
d'où
\begin{equation}\label{ApsiGdual} \widehat{A(\psi_G})=\widehat{ \bbZ/2\bbZ } \times \widehat{ A(\psi_{G'})},
\end{equation}

Si $x \in A(\psi_G)$, écrivons $x=(x_d,x')$ via (\ref{ApsiG})
  et si $\eta \in  \widehat{A(\psi_G})$, écrivons $\eta=(\eta_d,\eta')$ via (\ref{ApsiGdual})

L'élément $s_{\psi_G} \in A(\psi_G)$ s'écrit 
\begin{equation}\label{spsiG}
s_{\psi_G}=(s_d, s_{\psi_{G'}}), \text{ avec }\qquad s_d=(-1)^{c+1} . 
\end{equation}
On a alors 
\begin{align}\label{etaspsi}  \eta(s_{\psi_G})&=\eta_d(s_d)   \eta'(s_{\psi_G'})\\  
\label{Sbsd} S_i(s_d)&=  \begin{cases} (-1)^{i(c-i)}  (-1)^{\frac{c(c-1)}{2}(c+1)}  \text{ si } n-c \text{ est pair} \\
(-1)^{i(c-i)}  (-1)^{\frac{c(c+1)}{2}(c+1)}  \text{ si } n-c \text{ est impair} \end{cases}.
\end{align}

\end{lemme}

\dem  Donnons juste quelques indications et pour fixer les idées, 
plaçons nous dans le cas {\bf A}, où $\widehat G=\SO_{2n+1}(\bbC)$ plongé dans $\GL_{2n+1}(\bbC)$ par la représentation standard.
On a de manière évidente des inclusions de groupes :
\[  \mathbf S(\Or_{2c}(\bbC)\times \Or_{2(n-c)+1}(\bbC))\subset \widehat G=\SO_{2n+1}(\bbC)\subset \GL_{2n+1}(\bbC). \]
De la décomposition $\psi=\psi_d\oplus \psi'$ de $\psi$, on déduit une décomposition du centralisateur $S_{\psi_G}$
de l'image de $\psi_G$ dans $\widehat G$ en $S_{\psi_d}\times S_{\psi'}$, où 
$S_{\psi_d}$ est le centralisateur de $\psi_d$ dans $ \Or_{2c}(\bbC)$, et $S_{\psi'}$  le centralisateur de
 $\psi'$ dans $\Or_{2(n-c)+1}(\bbC)$.
On détermine facilement  $S_{\psi_d}$ : c'est le centre de $\Or_{2c}(\bbC)$, c'est-à-dire un groupe à deux élements, que l'on a identifié
à $\bbZ/2\bbZ$ dans l'énoncé du lemme, et il est inclus dans $\SO_{2c}(\bbC)$. On en déduit que $S_{\psi'}$ est dans 
$\SO_{2(n-c)+1}(\bbC)=\widehat G'$. La décomposition (\ref{ApsiG}) en découle, avec rappelons-le, le facteur 
$\bbZ/2\bbZ$ qui désigne en fait le centre de  $\Or_{2c}(\bbC)$  (dans le cas {\bf B}, on obtient le centre de $\Sp_{2c}(\bbC))$. 
Le reste de la démonstration est élémentaire, on en laisse les détails au lecteur. 
\qed

\begin{defi}\label{defietaeta}
 On pose,   pour tout $\eta=(\eta_d,\eta')\in \widehat{A(\psi_G)}=\widehat{ \bbZ/2\bbZ } \times\widehat{ A(\psi_{G'})}$, 
 \begin{equation}\label{XAJeta}
  \pi^{\flat}(\psi_G,\eta,G):= \bigoplus_{i \in \{0,\ldots ,c\}  \vert S_i=\eta_d } 
\widetilde  \caR^{d_i  }_{\frqqq _i , L_i, G }(\xi_{i,t}  \boxtimes \pi( \psi_{G'}, \eta',  G'_i ) ) ,\end{equation}
où  $\widetilde \caR^{d_i}_\frqqq   $ est le  foncteur d'induction cohomologique normalisé de la section \ref{NormIC}
({\sl cf. } Eq. (\ref{NormIC}))  en degré 
\[ d_i=\frac{1}{2}(\dim \mathbf G-\dim \mathbf L)-(q(G)-q(L_i) ). \]

Pour les formes intérieures pures des groupes orthogonaux, on pose
\begin{equation}\label{XAJetaSOpq}
  \pi^{\flat}(\psi_G,\eta,\SO(p,q)):= \bigoplus_{  (\mathbf Q_i=\mathbf L_i\mathbf V_i)\in\Sigma_{{}^dQ}^{\SO(p,q)}    \vert S_i=\eta_d } 
\widetilde  \caR^{d_i  }_{\frqqq _i , L_i, \SO(p,q) }(\xi_{i,t}  \boxtimes \pi( \psi_{G'}, \eta',  G'_i ) )  \end{equation}
avec 
\[ d_i=\frac{1}{2}(\dim \mathbf G-\dim \mathbf L)-(q(\SO(p,q))-q(L_i) ). \] 
\end{defi}

Enonçons maintenant le résultat principal de cette section
 \begin{thm} \label{thm:main} Les représentations  $ \pi^{\flat}(\psi_G,\eta,G)$ construites ci-dessus  sont bien celles
 construites par Arthur, autrement dit
\[ \pi(\psi_G,\eta,G)= \pi^{\flat}(\psi_G,\eta,G)= \bigoplus_{i \in \{0,\ldots ,c\}  \vert S_i=\eta_d } 
\widetilde  \caR^{d_i  }_{\frqqq _i , L_i, G }(\xi_{i,t}  \boxtimes \pi( \psi_{G'}, \eta',  G'_i ) ).\]
De même, pour les formes intérieures pures $\SO(p,q)$ des groupes orthogonaux, on a 
\[ \pi(\psi_G,\eta,\SO(p,q))= \pi^{\flat}(\psi_G,\eta,\SO(p,q))\]
\[= \bigoplus_{ (\mathbf Q_i=\mathbf L_i\mathbf V_i)\in\Sigma_{{}^dQ}^{\SO(p,q)}   \vert S_i=\eta_d } 
\widetilde  \caR^{d_i  }_{\frqqq _i , L_i, \SO(p,q)}(\xi_{i,t}  \boxtimes \pi( \psi_{G'}, \eta',  G'_i ) ).\]
 \end{thm}
 
La démonstration de ce théorème sera l'objet des prochains paragraphes. 
Nous la donnons dans le cas des groupes quasi-déployés, 
elle est identique pour les formes intérieures pures $\SO(p,q)$. Elle s'appuie en grande partie sur les résultats de \cite{AMR}
auquel nous renvoyons le lecteur. L'autre ingrédient de la démonstration est le théorème \ref{ThmICetTrans}.

Commençons par vérifier que l'on obtient bien une représentation virtuelle stable en remplaçant les 
$\pi(\psi_G,\eta,G)$  par les $\pi^{\flat}(\psi_G,\eta,G)$ dans (\ref{Xpsistable}).
\begin{lemme} \label{LemXAJst}
La représentation virtuelle 
\begin{equation}   \label{XAJst}
 [ \pi^\flat(\psi_{G},G )]:=\sum_{\eta \in \widehat{A(\psi_G)} } \eta(s_{\psi_G}) [\pi^{\flat}(\psi_G,\eta,G)]
  \end{equation} est stable. 
 \end{lemme}
 
\dem  
On a  avec les notations du Lemme \ref{LemmeApsi},  $ \eta(s_{\psi_G})=\eta_d(s_d) \eta'(s_{\psi_G'}) $, d'où 
\begin{align}   \label{SAJ}
  [\pi^{\flat}(\psi_G,G )]=& \sum_{ i\in \{0,\ldots, c\} }
  \sum_{\eta' \in \widehat{A(\psi_{G'} ) } }
 S_i(s_d)   \eta' (s_{\psi_{G'}} )   \;   \widetilde  \caR^{d_i }_{\frqqq _i , L_i, G }  \left(\xi_{i,t}  \boxtimes \pi(\psi_{G'},\eta', G'_i ) \right)\\
  \nonumber 
 =&\sum_{i \in \{0,\ldots, c\}}   S_i(s_d)    \widetilde  \caR^{d_i  }_{\frqqq _i , L_i, G } \left(\xi_{i,t} \boxtimes
\left( \sum_{\eta' \in \widehat{A(\psi_{G'})}  }\eta'(s_{\psi_{G'}}) \pi(\psi_{G'},\eta',G'_i ) \right) \right)\\
\nonumber =&\sum_{i\in \{0,\ldots, c\}}  S_i (s_d)    \;  \widetilde  \caR^{d_i  }_{\frqqq _i , L_i, G } \left(\xi_{i,t}  \boxtimes
\, [\pi(\psi_{G'},G'_i)] \right).
\end{align}

Omettons le $t$ des notations en posant $\xi_i=\xi_{i,t}$.  C'est un  caractère 
de $\U(i,c-i)$ admettant une résolution de Johnson par des caractères standard. Nous renvoyons
le lecteur à     \cite{AMR}, \S 10.3 pour les notations et assertions non expliquées qui suivent.
Ces caractères standard  sont paramétrés par $\bar s\in \frI_c^{i,c-i,\pm}$ et sont notés  $X_i(\bar s)$. 
Ces paramètres ont une   $\theta$-longueur,   notée $\ell_\theta(s)$  
 qui ne dépend que de l'image $s$ de $\bar s$ dans $\frI_c$
et la formule des caractères s'écrit : 
\begin{equation}\label{formcarU}
[\xi_{i}]= \sum_{\bar s \in \frI_c^{i,c-i,\pm} } (-1)^{ q(\U(i,c-i)) } (-1)^ {\ell_ \theta(w_0)+\ell_ \theta(s)}  X_i(\bar s).
\end{equation}

Or d'après (\ref{Sbsd}) et la remarque \ref{qUiai}
  $(-1)^{ q(\U(i,c-i)) }  S_i (s_d)=(-1)^{\frac{c(c\pm1)}{2}(c+1)}$. D'autre part 
 $\ell_ \theta(w_0)=q(\U_c)$.  Comme  $q(\U_c)=(c/2)^2$ si $c$ est pair et $q(U_c)$ est pair si 
$c$ est impair, on obtient $q(\U_c)={\frac{c(c\pm1)}{2}(c+1)}\mod 2$.
On continue alors le calcul de (\ref{SAJ}) : 
\begin{align}   \label{SAJ2}
&  [\pi^{\flat}(\psi_{G},G )]=\sum_{i\in \{0,\ldots, c\}}  \sum_{\bar s \in \frI_c^{i,c-i,\pm} }  (-1)^ {l_ \theta(s)}  
\;    \widetilde  \caR^{d_i }_{\frqqq _i , L_i, G } \left(  X_i (\bar s)   \boxtimes
[\pi(\psi_{G'},G'_i )] \right)
\end{align}

  Les paramètres de Langlands  de $\U(i,c-i)$ ayant le   caractère infinitésimal voulu sont en bijection avec $\frI_c$, 
  et le pseudo-paquet paramétré par un $s\in \frI_c$ est l'ensemble  des représentations standard
  $X_i(\bar s)$ où $\bar s$ est dans  $\frI_c^{i,c-i,\pm}$ de permutation sous-jacente $s$. 
  On peut donc continuer le calcul, en posant :
  \begin{equation}\label{Xisst}[X_i(s)^{st}]=   \sum_{\bar s \in \frI_c^{i,c-i,\pm}\vert p(\bar s)=s } 
   X_i(\bar s),  \end{equation}
   et l'on obtient: 
\begin{align}   \label{SAJ4}
&  [\pi^{\flat}(\psi_{G},G )]=
  \sum_{s \in \frI_c }  (-1)^{l_ \theta(s)}  \sum_{i\in \{0,\ldots, c\}}
  \, \widetilde  \caR^{d_i }_{\frqqq _i , L_i, G }  \left([X_i (\bar s)^{st}] 
   \boxtimes [\pi(\psi_{G'},G'_i )] \right).
\end{align}

Le résultat est alors conséquence du corollaire \ref{cor:stable3}. \qed
 
Pour démontrer le théorème \ref{thm:main}, on va vérifier les identités endoscopiques qui caractérisent les 
$\pi(\psi_G,\eta,G)$.

 \subsection{Transfert endoscopique tordu}

Notre but est de démontrer l'égalité de distributions sur $\widetilde G_N$ : 
\begin{equation}\label{butfinal}
\Tr_{\theta_N} (\Pi_\psi^\GL) =\Trans_G^{\widetilde{G_N}} ([ \pi^{\flat}(\psi_{G},G ) ]) .
\end{equation}

Pour démontrer ceci, nous partons d'une part de la formule suivante  démontrée dans \cite{AMR} 
      \begin{equation}\label{iDe2}
 \Tr_{\theta_{N_d}}( \Pi_{\psi_d}^ \GL)= \sum_{s \in \frI_c } (-1)^{\ell_\theta(s)}  
 \Tr_{\theta_{N_d}}( X_\GL(s)).
  \end{equation}
  Nous renvoyons à \cite{AMR} pour la définition de  $X_\GL(s)$. C'est une représentation standard entrant dans la résolution
  de Johnson de la représentation de Speh $\Pi_{\psi_d}^\GL$.

D'autre part nous écrivons  la représentation virtuelle  $ [ \pi^{\flat}(\psi_{G},G )]$
 sous la forme : 
\begin{equation}\label{FormCarXAJ}
 [ \pi^{\flat}(\psi_{G},G )]=\sum_{\phi_G\in \Phi(G)_{\psi_G}} a^{\flat}(\phi_G) [\pi(\phi_G,G)].  
\end{equation}
Ici $\Phi(G)_{\psi_G}$ désigne le sous-ensemble de $\Phi(G)$ des paramètres de Langlands ayant le même caractère infinitésimal que celui
de $\psi_G$ (c'est donc un ensemble fini). Les $[\pi(\phi_G,G)]$, lorsque $\phi_G$ décrit $\Phi(G)_{\psi_G}$, forment une base 
de l'espace des représentations virtuelles stables de caractères infinitésimal égal à celui de $\psi_G$, ce qui justifie l'existence et l'unicité
d'une telle écriture.

De même, écrivons la représentation virtuelle  $ [ \pi(\psi_{G'},G' )]$
 sous la forme : 
\begin{equation}\label{FormCarXAJprime}
 [ \pi(\psi_{G'},G' )]=\sum_{\phi_{G'}\in \Phi(G')_{\psi_{G'}}} a(\phi_{G'}) [\pi(\phi_{G'},G')].
\end{equation}

Par transfert entre formes intérieures, pour toutes les formes intérieures $G'_i$, on a 
\begin{equation}\label{FormCarXAJi}
 [ \pi(\psi_{G'},G_i' )]=  \sum_{\phi_{G'}\in \Phi(G')_{\psi_{G'}}} a(\phi_{G'})    [\pi(\phi_{G'},G'_i)].
\end{equation}

  La définition du paquet d'Arthur $\Pi(\psi_{G'},G')$ entraîne en particulier 
  l'identité de transfert
  \begin{equation}\label{TransG'}
  \tr_{\theta_{N'}} (\Pi^ \GL_{\psi_\epsilon'})=\mathrm{Trans}_{G'}^{\widetilde{G_{N'}}}
  ([ \pi(\psi_{G'},G' )]).
  \end{equation}
  
  On en déduit par (\ref{FormCarXAJprime}) et linéarité du transfert que 
  \begin{align}\label{TransG'2}
 \tr_{\theta_{N'}} (\Pi^ \GL_{\psi'_\epsilon})&=\sum_{\phi_{G'} \in \Phi(G')_{\psi_{G'}} }   a(\phi_{G'}) \;
 \mathrm{Trans}_{G'}^{\widetilde{G_{N'}}}([\pi(\phi_{G'},G')] )\\
 \nonumber & =\sum_{\phi_{G'} \in \Phi(G')_{\psi_{G'}} }   a(\phi_{G'}) \;  \tr_{\theta_{N'}} (\Pi^\GL_{\Std_{G'}\circ \phi_{G'}}).
  \end{align}
  
On tensorise par le caractère quadratique $\epsilon_{\psi}$ introduit plus haut pour récupérer
  
   \begin{equation}\label{TransG'3}
 \tr_{\theta_{N'}} (\Pi^ \GL_{\psi'})=\sum_{\phi_{G'} \in \Phi(G')_{\psi_{G'}} }   a(\phi_{G'}) \; 
  \tr_{\theta_{N'}} (\Pi^\GL_{(\Std_{G'}\circ \phi_{G'})\otimes \epsilon_{\psi}}).
  \end{equation}
  
  Les résultats de \cite{AMR}  impliquent que 
  
    \begin{equation}
 \Tr_{\theta_N} (\Pi_\psi^\GL)= \sum_{s\in \frI_c, \, \phi_{G'} \in \Phi(G')_{\psi_{G'}}} (-1)^{\ell_\theta(s)}  a(\phi_{G'})  \; 
  \tr_{\theta_N} \left(  \Pi^\GL_{ \phi_d(s)+ (\Std_{G'}\circ \phi_{G'}   ) \otimes \epsilon_{\psi'}  } \right)  .
 \end{equation}

Dans cette équation, $\phi_d(s)$ est le paramètre de Langlands 
de la représentation standard $X_\GL(w_0s)$ de $\GL_{N_d}$, et 
$\phi_d(s)+ (\Std_{G'}\circ \phi_{G'} ) \otimes \epsilon_{\psi'}  $
 est un paramètre de Langlands  pour $G_N$.

  En appliquant le transfert à (\ref{FormCarXAJ}), on obtient :
  \begin{equation}\label{FormCarXAJ11}
 \mathrm{Trans}_G^{\widetilde G_N}( [\pi^{\flat}(\psi_{G},G )]) =\sum_{\phi_G\in \Phi(G)_{\psi_{G}}} a^{\flat}(\phi_G)
 \mathrm{Trans}_G^{\widetilde G_N}
 ( [\pi(\phi_G,G)])\end{equation}
 \[= \sum_{\phi_G\in \Phi(G)_{\psi_{G}}} a^{\flat}(\phi_G)\tr_{\theta_N}
 \left( \Pi^\GL_{\Std_G\circ \phi_G}\right). \]

En comparant ces deux dernières formules, on voit donc que l'on s'est ramené à montrer 
   \begin{equation}\label{Amontrer}
  \sum_{s\in \frI_c, \, \phi_{G'} \in \Phi(G')_{\psi_{G'}}} (-1)^{\ell_\theta(s)}  a(\phi_{G'})  \; 
  \tr_{\theta_N} \left(  \Pi^\GL_{ \phi_d(s)+ (\Std_{G'}\circ \phi_{G'}   ) \otimes \epsilon_{\psi'}  } \right)  
    = \sum_{\phi_G\in \Phi(G)_{\psi_{G}}} a^{\flat}(\phi_G)\tr_{\theta_N} \left( \Pi^\GL_{\Std_G\circ \phi_G}\right). 
  \end{equation}
  
  Nous allons définir une  application injective :
  \begin{equation}\label{iota} \iota:\;  \frI_c \times \Phi(G') _{\psi_{G'}}\rightarrow \Phi(G)_{\psi_{G}}\end{equation}
 et  il s'agit donc de montrer que si $a^{\flat}(\phi_G) \neq 0$, alors $\phi_G$ est dans l'image de cette application, 
disons d'un élément $ ( s,\phi_{G'})$ et que l'on a alors 
  \begin{align}\label{Amontrer2}
 & a^{\flat}(\phi_G) = (-1)^{\ell_\theta(s)}  a(\phi_{G'}) \\
\label{Amontrer3} &\Std_G \circ \phi_G=  \phi_d(s)+ (\Std_{G'}\circ \phi_{G'}  ) \otimes \epsilon_{\psi}  .
  \end{align}

  \medskip

  On reprend (\ref{SAJ4}) qui donne, en utilisant aussi (\ref{FormCarXAJi}) 
    \begin{align}\label{yadlespoir}
 &[\pi^\flat(\psi_G,G)] =\sum_{\phi_G\in \Phi(G)_{\psi_{G}}} a^{\flat}(\phi_G) [\pi(\phi_G,G)]    \\
 \nonumber =&  \sum_{s \in \frI_c } (-1)^{\ell_ \theta(s)}  \sum_{i\in \{0,\ldots, c\}}\; 
    \widetilde  \caR^{d_i }_{\frqqq _i , L_i, G }    \left([X_i(s)^{st}] 
   \boxtimes [\pi(\psi_{G'},G'_i )] \right)\\
  \nonumber  =&  \sum_{s \in \frI_c } 
  (-1)^{l_ \theta(s)} \sum_{i\in\{0\ldots,c\} } \sum_{\phi_{G'} \in \Phi(G')_{\psi_{G'}} } a(\phi_{G'}) \;  
  \widetilde  \caR^{d_i }_{\frqqq _i , L_i, G }  \left([X_i(s)^{st}] 
   \boxtimes [\pi(\phi_{G'},G'_i )] \right).
 \end{align}
    
  \begin{lemme}
  Pour tout $s\in \frI_c$ et pour tout $\phi_{G'}\in \Phi(G')_{\psi_{G'}}$, il existe un unique paramètre $\phi_G\in \Phi(G)_{\psi_{G}}$ tel que 
\[  \sum_{i\in \{0,\ldots, c\}}    \;   \widetilde  \caR^{d_i }_{\frqqq _i , L_i, G }    \left( [ X_i(s) ^{st}] \boxtimes
 [\pi(\phi_{G'},G'_i)] \right)=[\pi(\phi_G,G)]. \]
Ceci définit l'application $\iota$ en (\ref{iota}). On a de plus 
\[\Std_G \circ \phi_G=  \phi_d(s)+ (\Std_{G'}\circ \phi_{G'}  ) \otimes \epsilon_{\psi}. \]
  \end{lemme}
     
  \dem On a vu  dans la démonstration de la proposition \ref{cor:stable2} que le terme de gauche
  est de la forme $[\pi(\phi_G,G)]$ pour un certain paramètre
     $\phi_G \in \Phi(G)_{\psi_{G}}$,.  Il faut bien sûr tenir compte de la renormalisation des foncteurs d'induction cohomologique, 
     mais celle-ci n'affecte pas la conclusion.   Il reste à voir que $\phi_G$ vérifie (\ref{Amontrer3}).
     Soit $X_u\boxtimes X'$ une représentation standard intervenant dans $[ X_i(s) ^{st}] \boxtimes [\pi(\phi_{G'},G'_i)]$.
     Les calculs fait dans l'appendice, sections \ref{ICRS} et \ref{calculTorsionI} expriment les paramètres de Langlands 
     de la représentation standard $\widetilde  \caR^{d_i }_{\frqqq _i , L_i, G }    \left( X_u\boxtimes X' \right)$
     en fonction de ceux de $X_u$ et $X'$. Sur le facteur groupe unitaire, celui-ci ayant ses sous-groupes de Cartan connexes,
      tout est déterminé par le caractère infinitésimal. Sur le facteur groupe classique, il faut voir que la torsion 
      $(\bbC_{\gamma_{\frn_L}})^*\otimes \bbC_{\gamma_\frn}$ calculé en \ref{calculTorsionI} se compense
      dans le cas des groupes orthogonaux avec la façon dont on a normalisé les foncteurs d'induction cohomologiques
      (le caractère de torsion $\xi'$ de (\ref{NormIC})).
      Pour les groupes symplectiques, nous n'avons pas pu inclure cette torsion dans la renormalisation des foncteurs d'induction 
      cohomologique, et ceci se traduit par la torsion par $\epsilon_{\psi}$ dans (\ref{Amontrer3}).\qed

   La formule (\ref{yadlespoir}) montre alors que l'on a $ a^{\flat}(\phi_G) = (-1)^{l_\theta(s)}  a(\phi_{G'})$.
  Ceci termine la démonstration de (\ref{butfinal}).

  \subsection{Transfert endoscopique ordinaire}
  
Soit  $(H,x,\caH,\iota_{H,G}:{}^LH\rightarrow {}^LG)$  une donnée endoscopique elliptique de $G$ comme dans la section \ref{DEE}.
Ici, $x$ est un élément de $\widehat G$ et l'on suppose que $x^2=1$, $H$ est un produit de groupes classiques, 
disons $H=H_1\times H_2$, et une inclusion de $\widehat H_1 \times \widehat H_2$ dans $\widehat G$
correspondant aux valeurs propres $\pm 1$ de $x$ (mais on ne suppose pas que $H_1$ correspond à la valeur propre
$1$). 
On écrit 
\begin{equation} \label{Psietal1}
 \psi_G=\xi_{H,G} \circ \psi_H, \qquad  \psi_H=\psi_{H_1}\times \psi_{H_2}, \quad \psi=\Std_G\circ \psi_G, \end{equation}
\begin{equation} \label{Psietal2} \psi= \psi_d+\psi'=\psi_1+\psi_2 = \psi_d+\psi'_1+\psi_2 , \qquad  \psi'=\psi_1'+\psi_2.\end{equation}
Les indices $1$ et $2$ sont
 déterminés par le fait que 
 \begin{equation}\label{Psietal3} \Std_{H_1} \circ \psi_{H_1}= \psi_{1}^{\epsilon}  , \qquad   
\Std_{H_2} \circ \psi_{H_2}=\psi_{2}^{\epsilon}. \end{equation}
Ici, $\psi_{1}^{\epsilon}=\psi_1$ sauf dans le cas où $\mathbf H_1\times \mathbf H_2=\Sp_{2a}\times \SO^{qd}_{2b}$ où 
$\psi_{1}^{\epsilon}=\sgn_{W_\bbR}\otimes \psi_1$, et 
 $\psi_{2}^{\epsilon}=\psi_2$ sauf   dans le cas où $\mathbf H_1\times \mathbf H_2=\SO^{qd}_{2b}\times \Sp_{2a}$ où 
$\psi_{2}^{\epsilon}=\sgn_{W_\bbR}\otimes \psi_2$   et $\psi_{1}^{\epsilon}=\psi_1$.
Ces derniers cas sont ceux où les paramètres $\psi_1$ et $\psi_2$ ne sont pas à valeurs dans $\SO(2a+1,\bbC)$ mais seulement dans 
$\Or(2a+1,\bbC)$. On corrige donc en tensorisant par $\sgn_{W_\bbR}$ pour obtenir un paramètre $\psi_{1}^{\epsilon}$ ou $\psi_{2}^{\epsilon}$
à valeurs dans  $\SO(2a+1,\bbC)$. Cette correction apparaît aussi dans le dernier cas de la remarque \ref{plongH}.

\medskip

Autrement dit, on se retrouve dans la situation de la section \ref{cLeviH}, avec 
$\mathbf  L_*=\U_c\times \mathbf G'$, $\mathbf H=\mathbf H_1\times \mathbf H_2$, 
$\mathbf  L_{*,H}=\U_c\times \mathbf H'_1\times \mathbf H_2$, et 
$\mathbf  L_{*,H}$ est un groupe endoscopique pour $\mathbf L_*$ (via le fait que 
 $\mathbf H'_1\times \mathbf H_2$ est un groupe endoscopique pour $\mathbf G'$, et c'est aussi un 
 $c$-sous-groupe de Levi de $\mathbf H$   (via le fait que 
 $\U_c\times \mathbf H'_1$ est un   $c$-sous-groupe de Levi de $\mathbf H_1$).  

\bigskip 


Si $G$ est symplectique et $c$ est impair, on a introduit le  paramètre $\psi'_\epsilon=\sgn_{W_\bbR}\otimes \psi'$,
 pour pouvoir factoriser
   \begin{equation}\label{Psietal5} \psi'_\epsilon=\Std_{G'}\circ \psi_{G'}.\end{equation}
 On a alors  
   \begin{equation}\label{Psietal6} \psi'_\epsilon=\sgn_{W_\bbR}\otimes \psi'_1 +  \sgn_{W_\bbR}\otimes \psi_2=\psi'_{1,\epsilon}+\psi_{2,\epsilon}.
 \end{equation}
 \medskip 
 
 On applique la discussion ci-dessus avec $(\mathbf G', \mathbf H'_1\times \mathbf H_2, \psi_{G'})$ plutôt que 
 $(\mathbf G,\mathbf H_1\times \mathbf H_2,\psi_G)$.
 Dans (\ref{Psietal3}), à  la place de $\psi_{H_1}$ et $\psi_{H_2}$, on récupère des paramètres pour $\mathbf H'_1$ et $\mathbf H_2$ 
 que l'on va noter 
 \begin{equation}\label{psiH1psiH2}  \psi'_{H'_1}\quad \text{ et  } \quad  \psi'_{H_2}. \end{equation}
 
 On peut aussi appliquer  la discussion ci-dessus avec $(\mathbf H_1,\mathbf H_1', \psi_{H_1})$ plutôt que $(\mathbf G,
 \mathbf G',\psi_G)$, et l'on récupère un paramètre  pour $\mathbf H'_1$ que l'on va noter 
  \begin{equation}\label{psiH1'} \psi_{H'_1}.\end{equation}

 Dans le cas des groupes orthogonaux, on a 
 \begin{equation}\label{psipsi'0}
   \Std_{H_1'}\circ \psi'_{H'_1}=\psi'_1 =  \Std_{H_1'}\circ \psi_{H'_1}, \quad \psi'_{H'_1}=\psi_{H'_1}, \quad \psi'_{H_2}=\psi_{H_2}. 
   \end{equation}

Dans le cas des groupes symplectiques, le mieux est de voir ce qui se passe dans tous les cas. Remarquons que l'on a
 $\sgn_{W_\bbR}\otimes \psi_d\simeq \psi_d$. 

$c$ pair, $\mathbf H_1\times \mathbf H_2=\Sp_{2a}\times \SO_{2b}^d$, ou $\mathbf H_1\times \mathbf H_2= \SO_{2b}^d\times \Sp_{2a}$.
Dans ce cas, aucune torsion n'apparaît sur les paramètres, et comme pour  les groupes orthogonaux on a 
\begin{equation}\label{psipsi'} 
\Std_{H_1'}\circ \psi'_{H'_1}=\psi'_1=  \Std_{H_1'}\circ \psi_{H'_1},\quad  \psi'_{H'_1}=\psi_{H'_1}, \quad \psi'_{H_2}=\psi_{H_2}.\end{equation}

 $c$ pair, $\mathbf H_1\times \mathbf H_2=\Sp_{2a}\times \SO_{2b}^{qd}$, il apparaît la torsion $\psi_1^\epsilon$ pour définir
 $\psi_{H_1}$, et la même torsion ${\psi'_1}^\epsilon$  pour définir  $\psi'_{H'_1}$ avec $\psi_1^\epsilon=\psi_d+{\psi'_1}^\epsilon$.
  On a donc  
  \begin{equation}\label{psipsi'2}   \Std_{H_1'}\circ \psi'_{H'_1}={\psi'_1}^\epsilon=  \Std_{H_1'}\circ \psi_{H'_1},
  \quad  \psi'_{H'_1}=\psi_{H'_1},  \quad \psi'_{H_2}=\psi_{H_2}.  
  \end{equation}

 $c$ pair, $\mathbf H_1\times \mathbf H_2= \SO_{2b}^d\times \Sp_{2a}$ . Il apparaît la torsion $\psi_2^\epsilon$ pour définir
 $\psi_{H_2}$, et la même torsion $\psi_2^\epsilon$  pour définir  $\psi'_{H_2}$.   On a donc 
  \begin{equation}\label{psipsi'3}  \Std_{H_1'}\circ \psi'_{H'_1}={\psi'_1}=  \Std_{H_1'}\circ \psi_{H'_1}, 
  \quad \psi'_{H'_1}=\psi_{H'_1}, \quad \psi'_{H_2}=\psi_{H_2}.
  \end{equation}

$c$ impair, $\mathbf H_1\times \mathbf H_2=\Sp_{2a}\times \SO_{2b}^d$.
On a la torsion $\psi'_\epsilon$ qui définit $\psi_{G'}$ et la même torsion sur $\psi'_{1,\epsilon}$  pour définir  $\psi'_{H'_1}$ et $\psi_{H'_1}$.
 On a donc 
\begin{equation}\label{psipsi'4} 
\Std_{H_1'}\circ \psi'_{H'_1}=\psi'_{1,\epsilon}=  \Std_{H_1'}\circ \psi_{H'_1}, \quad \psi'_{H'_1}=\psi_{H'_1}, \quad \psi'_{H_2}=\sgn_{W_\bbR}\otimes \psi_{H_2}. \end{equation}

$c$ impair,  $\mathbf H_1\times \mathbf H_2= \SO_{2b}^d\times \Sp_{2a}$.
On a la torsion $\psi'_\epsilon$ qui définit $\psi_{G'}$, et la même torsion sur $\psi'_{1,\epsilon}$  pour définir  $\psi'_{H'_1}$, 
 mais que l'on ne retrouve pas pour définir $\psi_{H'_1}$ car $\mathbf H_1$ est orthogonal.
On a donc  
 \begin{equation}\label{psipsi'5} 
\Std_{H_1'}\circ \psi'_{H'_1}=\psi'_{1,\epsilon}, \quad   \Std_{H_1'}\circ \psi_{H'_1}=\psi'_1,
 \quad \psi'_{H'_1}=\sgn_{W_\bbR}\otimes \psi_{H'_1}, \quad \psi'_{H_2}=\sgn_{W_\bbR}\otimes \psi_{H_2}. \end{equation}

 $c$ impair, $\mathbf H_1\times \mathbf H_2=\Sp_{2a}\times \SO_{2b}^{qd}$. Il apparaît la torsion $\psi'_\epsilon$ qui définit $\psi_{G'}$ et 
 la torsion $\psi_1^\epsilon$ pour définir $\psi_{H_1}$. La torsion qui définit $\psi'_{H'_1}$ et $\psi_{H'_1}$ est la même. 
  On a donc  
  \begin{equation}\label{psipsi'6}   \Std_{H_1'}\circ \psi'_{H'_1}={\psi'_{1,}}_\epsilon^\epsilon=\psi'_1= \Std_{H_1'}\circ \psi_{H'_1}
  \quad  \psi'_{H'_1}=\psi_{H'_1},  \quad \psi'_{H_2}=\sgn_{W_\bbR}\otimes \psi_{H_2}.  
  \end{equation}

 $c$ impair, $\mathbf H_1\times \mathbf H_2= \SO_{2b}^{qd}\times \Sp_{2a}$. Il apparaît la torsion  
 $\psi'_\epsilon$ qui définit $\psi_{G'}$  et la même torsion sur $\psi'_{1,\epsilon}$  pour définir  $\psi'_{H'_1}$, 
 mais que l'on ne retrouve pas pour définir $\psi_{H'_1}$ car $\mathbf H_1$ est orthogonal. Il y  a aussi  la torsion $\psi_2^\epsilon$ qui  définit
 $\psi_{H_2}$. 
    On a donc 
  \begin{equation}\label{psipsi'7} 
  \Std_{H_1'}\circ \psi'_{H'_1}=\psi'_{1,\epsilon}, \quad   \Std_{H_1'}\circ \psi_{H'_1}=\psi'_1,
 \quad \psi'_{H'_1}=\sgn_{W_\bbR}\otimes \psi_{H'_1}, \quad \psi'_{H_2}=\sgn_{W_\bbR}\otimes \psi_{H_2}. \end{equation}

\bigskip 
 Reprenons la démonstration de l'identité de transfert endoscopique. Il s'agit ici de montrer
\begin{equation}\label{butordi}  
\mathrm{Trans}_H^G([\pi^{\flat}(\psi_H,H)])= \sum_{\eta \in \widehat{A(\psi_G)}} \eta(s_{\psi_G}x) \; 
[\pi^\flat(\psi_G,\eta,G)].
  \end{equation}
  
  On remarque qu'il n'y a pas de signe de Kottwitz ici car $H$ et $G$ sont quasi-déployés.
Commençons par   le membre de droite :
\begin{align}
  &\sum_{\eta \in \widehat{A(\psi_G)}} \eta(s_{\psi_G}x) \; [\pi^\flat(\psi_G,\eta,G)]  \\
 \nonumber &=\sum_{(\eta_d,\eta') \in \widehat{\bbZ/2\bbZ} \times \widehat{A(\psi_{G'})} } \sum_{i\in\{0,\ldots, c\}\vert S_i=\eta_d}
    \eta_d(s_dx_d)  \eta'(s_{\psi_{G'}}x') \;  \ \widetilde  \caR^{d_i }_{\frqqq _i , L_i, G }\left(   \xi_{i,t} \boxtimes \pi(\psi_{G'},\eta',G'_i)  \right)\\
\nonumber    &= \sum_{i\in\{0,\ldots, c\}} S_i(s_d x_d) \;   \widetilde  \caR^{d_i }_{\frqqq _i , L_i, G } \left(   \xi_{i,t} \boxtimes \left(
    \sum_{\eta' \in  A(\psi_{G'}) } \eta'(s_{\psi_{G'}}x')  [ \pi (\psi_{G'},\eta',G'_i)] \right) \right)
  \end{align}
  
  Comme $(H'_1\times H_2,\ldots)$ est une donnée endoscopique pour $G'_i$, on a une identité de transfert
  \begin{equation}
   e(G'_i) \sum_{\eta' \in \widehat{ A(\psi_{G'}) }} \eta'(s_{\psi_{G'}}x')  [ \pi(\psi_{G'},\eta',G'_i)]  =
    \mathrm{Trans}_{H'_1\times H_2}^{G'_i}([\pi(\psi'_{H'_1}\times \psi'_{H_2},H'_1 \times H_2)])
  \end{equation}
  et ainsi, en utilisant aussi (\ref{formcarU}) et (\ref{Xisst}):
 \begin{align}
  &\sum_{\eta \in \widehat{A(\psi_G)}} \eta(s_{\psi_G}x) \;  [\pi^\flat(\psi_G,\eta,G)] \\
\nonumber&= \sum_{i\in\{0,\ldots, c\}} S_i(s_d x_d) e(G'_i)\;   \widetilde  \caR^{d_i }_{\frqqq _i , L_i, G } \left(   \xi_{i,t}
\boxtimes \left(
   \mathrm{Trans}_{H'_1\times H_2}^{G'_i}([\pi(\psi'_{H'_1}\times \psi'_{H_2},H'_1\times H_2)])
  \right) \right)\\
     \nonumber&= \sum_{i\in\{0,\ldots, c\}}    \sum_{\bar s \in \frI_c } (-1)^{l_ \theta(\bar s)}
    S_i(x_d) \, 
    \widetilde  \caR^{d_i }_{\frqqq _i , L_i, G } \left(  [X_i(\bar s)^{st}]      \boxtimes \left(
   \mathrm{Trans}_{H'_1\times H_2}^{G'_i}\left([\pi(\psi'_{H'_1}, H'_1)] \boxtimes  
   [\pi(\psi'_{H_2},H_2)]\right)
  \right) \right)
  \end{align}

D'autre part, pour réécrire le membre de gauche, on écrit en utilisant (\ref{SAJ4}) avec $H_1$ à la place de $G$ :
\begin{align}
&[\pi(\psi_H,H)]=[\pi(\psi_{H_1},H_1)] \boxtimes [\pi(\psi_{H_2},H_2)]\\
\nonumber
  &=\left( \sum_{\bar s \in \frI_c }   (-1)^{l_ \theta(\bar s)}  \sum_{k \in \{0,\ldots, c\}}
\widetilde \caR^{d_{1,k}}_{\frqqq_{k,H_1}, L_{k,H_1},H_1} \left([X_k(\bar s)^{st}] 
   \boxtimes [\pi(\psi_{H'_1},H'_k)] \right)  \right) \boxtimes [\pi(\psi_{H_2},H_2)]
\end{align}

On obtient la conclusion voulue grâce à la relation (\ref{Formquisert}).
 
Pour les groupes symplectiques, on vérifie que les torsions sur les paramètres pour passer de $\psi'_{H'_1}$ et $\psi'_{H_2}$
 à   $\psi_{H'_1}$ et $\psi_{H_2}$ données ci-dessus coïncident bien avec les torsions 
$\tilde \varepsilon_1$ et $ \tilde\varepsilon_2$ de la formule 
 (\ref{Formquisert}).
\qed

\section{Paquets d'Arthur pour les paramètres de bonne parité}\label{PABPR}

\subsection{Enoncé du théorème principal}
Soit $G$ un groupe classique, et soit $\psi_G$ un paramètre d'Arthur pour $G$.
Supposons que $\psi_G$ soit de bonne parité, on décompose donc $\psi=\Std_G\circ \psi_G$ en 
\[  \psi=  \left(\bigoplus_{r=1,\ldots, R} V(0,t_r)\boxtimes R[a_r]\right) \oplus   \left(  \bigoplus_{m=1,\ldots ,M}
 W(0,\epsilon_m) \boxtimes R[a'_m] \right) =\psi_d\oplus \psi_{bp,u}\]
où $\psi_{bp,u}: =  \bigoplus_{m=1,\ldots, M} W(0,\epsilon_m) \boxtimes R[a'_m]   $ est unipotent de bonne parité et 
$\psi_d=\bigoplus_{r=1,\ldots, R} V(0,t_r)\boxtimes R[a_r]$ est la partie discrète de bonne parité.
On note $n$ le rang de $\mathbf G$, et l'on pose $c=\sum_{r=1}^R a_r$, $n'=n-c$.

Quitte à remplacer $\pi_G$ par un paramètre dans la même classe de conjugaison, on peut supposer, et c'est ce que l'on fait, que 
\begin{equation}\label{CondFair}
\forall r\in \{1, \ldots R-1\}, \quad t_r\geq t_{r+1}, \qquad \text{ et si } t_r=t_{r+1} \text { alors } a_r\geq a_{r+1}.
\end{equation}

Dans le cas {\bf A} (groupes symplectiques), soit $f$ le nombre de  $r\in \{1,\ldots R\}$ tel que $t_r$ est pair.
Si $f$ est impair, on pose $\epsilon_\psi=\sgn_{W_\bbR}$, et dans tous les autres cas (donc 
en particulier pour les groupes orthogonaux), on pose $\epsilon_\psi=\Triv_{W_\bbR}$, le caractère trivial de 
$W_\bbR$. On pose aussi $\psi_{bp,u,\epsilon}=\psi_{bp,u}\otimes \epsilon_{\psi}$.
Alors $\psi_{bp,u,\epsilon}$ se factorise en $\Std_{G'}\circ \psi_{G'}$ pour un groupe classique 
quasi-déployé $\mathbf G'$ de  même type que $\mathbf G$. Le paramètre  $\psi_{G'}$ est unipotent de bonne parité
 pour $\mathbf G'$ et ses formes intérieures.
Pour une forme intérieure pure $\mathbf G'_{int}$ de $\mathbf G'$, on a d'après les résultats rappelés dans la section 
\ref{parunip}, la représentation ${\pi}^A(\psi_{G'}, G'_{int})$ de $G'_{int}\times A(\psi_{G'})$.
Notre but est ici de donner une formule pour la représentation ${\pi}^A(\psi_G,G)$ en fonction des 
${\pi}^A(\psi_{G'}, G'_i)$ lorsque $\mathbf G'_{int}$ décrit l'ensemble des formes intérieures pures de
$\mathbf G'$ et des couples $(t_r,a_r)$, $r=1,\ldots R$.

\medskip

Considérons un sous-groupe de parabolique  standard ${}^dQ={}^dL{}^dV$ de $\widehat G$ comme dans la section \ref{cLevicla}, à ceci près
que ce n'est plus un sous-groupe parabolique maximal, le $c$-Levi associé $\mathbf L_*$ est ici de la forme
\[   \mathbf L_*=\U_{a_1}\times \cdots \times \U_{a_R}\times \mathbf G' .\]

On fixe comme en (\ref{LiFi})
un système de représentants $(\mathbf Q_{\underline i}=\mathbf L_{\underline{i}} \mathbf V_{\underline i } )_{\underline i}$ 
de $\Sigma^G_{{}^dQ}$
avec \[  L_{ \underline{i} } \simeq \U(i_1,a_1-i_1)\times \cdots \times    \U(i_R,a_R-i_R)   \times  G'_{\underline i} \] où $G'_{\underline i}$ 
 est une forme  intérieure pure de $G'$,  et $\underline{i}$ parcourt l'ensemble  $I$ défini comme suit  :
 
 \begin{defi} \label{def:I}Notons $I$  l'ensemble des $ (i_1,\ldots i_R)\in \bbN^R$  vérifiant : 
 
-  dans le cas {\bf A }  (groupes symplectiques) :  pour tout $r \in \{1,\ldots,R \}$,  $i_r\in \{0,a_r\}$.
 
- dans le cas des  groupes $\SO(p,q)$ :  pour tout $r\in \{1,\ldots,R \}$, 
 $i_r\in \{0,a_r\}$  et $\displaystyle \sum_{r=1,\ldots,R}i_r\leq p$, $\displaystyle \sum_{r=1,\ldots,R} a_r- i_r\leq q$.
  \end{defi} 
 Dans le cas des groupes $\SO(p,q)$,  on a  alors $\mathbf G'_{\underline i}= \SO(p',q')$ avec 
 \[ \displaystyle p'=p-\sum_{r=1,\ldots,R ,}i_r,  \quad \displaystyle q'=q- \sum_{r=1,\ldots,r,} a_r- i.\]

On définit le caractère  $\xi_{i_r,t_r}$  de $\U(i_r,a_r-i_r)$ par le fait que sa différentielle est $\left( \lfloor \frac{t_r}{2}
\rfloor, \ldots , \lfloor \frac{t_r}{2}\rfloor \right)$.
Définissons maintenant les foncteurs d'induction normalisés dans le cadre présent, généralisant ce que l'on a fait dans la section \ref{SecNormIC}.
On part du foncteur d'induction cohomologique non normalisé ${}^n\caR^{d_{\underline i}}_{\frqqq _{\underline i} , L_{\underline i}, G }$, en degré
\[d_{\underline i}=\frac{1}{2} (\dim \mathbf G -\dim \mathbf L_*)-(q(G)-q(L_{\underline i})).\]
Si $Z=Z_{u,1}\boxtimes \cdots \boxtimes Z_{u,R}\boxtimes Z'$ est une représentation de $L_i$, on pose 
\[  \widetilde {\caR}^{d_{\underline i}}_{\frqqq _{\underline i} , L_{\underline i}, G }(Z)=
{}^n\caR^{d_{\underline i}}_{\frqqq _{\underline i} , L_{\underline i}, G } \left(
\boxtimes_{r=1,\ldots R} (Z_{u,r}\otimes \xi_{u,r})\boxtimes (Z'\otimes \xi')\right),   \]
où le caractère de torsion $\xi_{u,r}$ du facteur $\U(i_r,a_r-i_r)$ est défini comme dans la définition \ref{def:xiu}, 
en remplaçant $c$ par $a_r$ et $n$ par $\displaystyle n-\sum_{j=1}^{r-1}a_j$.
On remarque que la condition de bonne parité du paramètre fait que l'on peut enlever les parties entières dans la formule
exactement lorsque $t_r$ est pair. De même, le caractère $\xi'$ de $G_i'$ est 
trivial dans le cas {\bf A} des groupes symplectique, et vaut $\sgn_{G'_i}^c $ dans le cas des groupes orthogonaux.


La condition (\ref{CondFair}) entraine que les inductions cohomologiques 
que nous considérons sont dans le weakly fair range (\cite{KV}, Definition 0.49). Par conséquent
 ces inductions cohomologiques sont nulles en tout degré sauf au plus un et envoient une représentation 
 unitaire sur une représentation unitaire. Mais elles n'ont pas de raison de conserver l'irréductibilité contrairement à ce qui se passe dans le good range.


Si ${\pi}=\tau\boxtimes \rho $ est une représentation de $L_{\underline i}\times A(\psi_{G'})$, 
on note 
\begin{equation}\label{def:Rgras}
\widetilde {{\mathcal R}}_{\underline i}     (\tau\boxtimes \rho)
= \widetilde  \caR^{d_{\underline i}  }_{\frqqq _{\underline i} , L_{\underline i}, G }(\tau )\boxtimes \rho.
\end{equation}
On obtient une représentation de $G \times  A(\psi_{G'})$

Comme dans le  lemme \ref{LemmeApsi}, on a une décomposition
$S_{\psi_G}=S_{\psi_d}\times S_{\psi_{bp,u}}$, où $S_{\psi_d}$ est le centralisateur de 
$\psi_d$ dans $\Or(2c,\bbC)$ ou  $\Sp(2c,\bbC)$ et  $S_{\psi_{bp,u}}$
est le centralisateur de  $\psi_{bp,u}$ dans $\Or(2n',\bbC)$ ou $\Sp(2n',\bbC)$ (les groupes symplectiques apparaissent dans le cas {\bf B} de 
la section \ref{grcla}, les groupes orthogonaux dans les autres cas), et cette décomposition induit un isomorphisme
(avec les notations évidentes)
$A(\psi_G)\simeq A(\psi_d)\times A(\psi_{bp,u})$ avec $A(\psi_{bp,u})\simeq A(\psi_{G'})$.
Notons  $Z_r$  le centre de $\Or(2a_r,\bbC)$,  un groupe à deux éléments dont on note $z_r$ l'élément non trivial.
  On a alors une   une surjection
\begin{equation}\label{Apsi6}  \prod_{r=1, \ldots ,R} Z_r   \rightarrow  A(\psi_{d}) .
\end{equation}
qui induit une surjection
\begin{equation}\label{Apsi5}  \prod_{r=1, \ldots ,R} Z_r \times A(\psi_{G'}) \rightarrow  A(\psi_{G}) .
\end{equation}
Le noyau de la surjection (\ref{Apsi6}) est engendré par les éléments $z_r z_{r+1}$ où $r$ est tel que 
$t_r=t_{r+1}$ et $a_r=a_{r+1}$. 

Soient $p,q,p',q'$ des entiers. Le groupe $\U(p+p',q+q')$ admet un $c$-sous groupe de Levi de la forme 
$\U(p,q)\times \U(p',q')$, et l'on a donc un  foncteur d'induction cohomologique
$\caR_{U(p,q)\times \U(p',q'),\U(p+p',q+q')}^{d_{p,q,p',q'}}$. Nous laissons au lecteur le soin de déterminer
le $c$-sous parabolique absent de la notation et le \og bon \fg \,  degré $d_{p,q,p',q'}$. Nous renormalisons 
ce foncteur comme dans la partie \ref{SecNormIC}, pour obtenir un foncteur 
$\widetilde \caR_{U(p,q)\times \U(p',q'),\U(p+p',q+q')}^{d_{p,q,p',q'}}$

\begin{lemme} [\cite{Trap}, Lemma 8.7, 9.3 ] Soit $r\in\{1,\ldots, R\}$  tel que 
$t_r=t_{r+1}$ et $a_r=a_{r+1}$. Alors 
\[ \widetilde \caR^{d_{i_r,a_r-i_r,i_{r+1}, a_{r+1}-i_{r+1}}}_{\U(i_r,a_r-i_r)\times \U(i_{r+1},a_{r+1}-i_{r+1}),\U(i_r+ i_{r+1},a_r+a_{r+1}-i_r-i_{r+1} )}
\left( \xi_{i_r,t_r}\boxtimes \xi_{i_{r+1}, t_{r+1}}\right) =0\] 
sauf si $i_r=a_{r+1}-i_{r+1}$ et $i_{r+1}=a_r-i_r$.
\end{lemme}
Il faut remarquer que ce résultat est valable dans tout le \og mediocre range \fg, donc en particulier dans le weakly fair range
où nous l'appliquons ici.

\medskip

Notons $S_{\underline i}$ le caractère de $\displaystyle \prod_{r=1, \ldots ,R} Z_r$, donné sur le $r$-ième  facteur
 $Z_r=\bbZ/2\bbZ$ par le caractère 
$S_{i_r}$ analogue au caractère  défini en  (\ref{Si}), avec $c$ remplacé par $a_r$ et $n$ remplacé par $\displaystyle n-\sum_{j=1}^{r-1} a_j$.

Nous avons introduit toutes les notations nécessaires pour énoncer notre résultat principal.
\begin{thm}\label{thm:Princ} Avec les notations qui précèdent
\begin{equation}\label{formulefinale}
{\pi}^A(\psi_G,G)= \bigoplus_{\underline i\in I}   S_{\underline i}\boxtimes  \widetilde {{\mathcal R}}_{\underline i}  \left(
\left( \boxtimes_{r=1,\ldots R} \; \xi_{i_r,t_r}\right)\boxtimes {\pi}^A(\psi_{G'}, G'_i)
\right).
\end{equation}
\end{thm}

La somme est sur l'ensemble $I$ de la définition \ref{def:I}. Expliquons le terme de droite : 
${\pi}^A(\psi_{G'}, G'_i)$ est une représentation de $G'_i\times A(\psi_{G'})$, et  
$\left(\boxtimes_{r=1,\ldots R} \; \xi_{i_r,t_r}\right)\boxtimes {\pi}^A(\psi_{G'}, G'_i)$ est une représentation de 
$L_i\times  A(\psi_{G'})$, donc 
$ \widetilde {{\mathcal R}}_{\underline i}  \left(
\left( \boxtimes_{r=1,\ldots R} \; \xi_{i_r,t_r}\right)\boxtimes {\pi}^A(\psi_{G'}, G'_i)
\right)$  est une représentation de 
$G \times  A(\psi_{G'})$. 
D'après le lemme ci-dessus, et par transitivité de l'induction cohomologique,  si $ \widetilde {{\mathcal R}}_{\underline i}  \left(
\left( \boxtimes_{r=1,\ldots R} \; \xi_{i_r,t_r}\right)\boxtimes {\pi}^A(\psi_{G'}, G'_i)\right)\neq 0$, alors on calcule facilement que 
$S_{\underline i}$ est trivial  sur le noyau de (\ref{Apsi6}) et l'on utilise  (\ref{Apsi5}) pour voir les termes de la somme
du membre de droite de (\ref{formulefinale}) comme des représentations de $G\times A(\psi_G)$.

\medskip

\dem On suppose dans un premier temps que $\psi_G$ vérifie la condition de régularité suivante :
\begin{equation}\label{conditionregul}  \forall r=1,\ldots ,R-1, \;
t_r-a_r+1>t_{r+1}+a_{r+1}-1 \; \text{ et }  t_{R} -a_{R}+1> \max_{m=1,\ldots M } \{a_m-1\}
\end{equation}
Le résultat  s'obtient alors  immédiatement  à partir  de la section précédente 
par une récurrence  sur  $R$. La condition de régularité (\ref{conditionregul}) assure qu'à chaque étape, la condition
 (\ref{hyplemmegen2}) est satisfaite et que les inductions cohomologiques se font à chaque fois dans le good range.
 On utilise aussi la transitivité de l'induction cohomologique. De plus, dans ce cas, notons que (\ref{Apsi5}) est un isomorphisme.

Pour passer du cas régulier au cas général, on utilise les résultats de \cite{MR4}, Théorèmes  4.1.1 et 4.5.1. \qed

\begin{rmq}
Les coefficients qui interviennent dans la décomposition donnée sont égaux à un. Mais malheureusement sans 
la condition de régularité de la preuve on ne sait pas que les induites cohomologiques sont irréductibles et inéquivalentes. 
Par contre, on le sait sous la condition de régularité car on est alors dans le good range. Sous la condition de régularité,
 on a donc démontré que $\pi^A(\psi)$ est une somme sans multiplicité de représentations irréductibles unitaires 
 et pour $\pi$ une représentation irréductible de $G$  intervenant dans la restriction de $\pi^A(\psi)$ à
  $G$, le coefficient est un caractère de $A(\psi)$. C'est la propriété de multiplicité un utile pour le calcul des multiplicités globales.
\end{rmq}

\subsection{Au sujet des choix\label{choix}}
Comme il est bien connu, le choix de la normalisation locale des facteurs de transfert est un point délicat même si ici, cela 
a pu se faire grâce aux travaux de Kottwitz et Shestad de façon assez simple. La vox populi dit que la normalisation se fait à l'aide 
d'un modèle de Whittaker, c'est un peu abusif car il faut aussi le choix d'une forme intérieure quasi-déployée de référence.
 Par exemple, dans le cas du groupe très simple $\SO_{2n+1}$, il n'y a qu'un seul modèle de Whittaker mais il y a deux choix de formes
  bilinéaires symétriques non dégénérées de dimension $2n+1$ ayant un espace isotrope de dimension $n$. Pour ce groupe c'est là qu'est
   le choix. Donc pour tous les groupes considérés ici, il y a deux choix possibles, soit il y a deux choix de modèles de Whittaker,
    c'est le cas des groupes symplectiques et des groupes orthogonaux pairs dont la forme intérieure quasi-déployée est en
     fait déployée et dans les cas restant il y a deux choix de formes bilinéaires symétriques non dégénérées maximalement isotropes.

Les choix faits ont été explicités dans le paragraphe 6.1 voir en particulier la remarque \ref{ChoiX}, et on peut les justifier par 
 leur compatibilité à l'induction cohomologique.  Par exemple pour le groupe $\SO_{2n+1}$ les formes bilinéaires considérées pour 
déterminer la forme intérieure sont celles qui ont un nombre impair de plus. Ces choix ont en revanche  le  défaut  de ne pas être compatible 
à l'induction parabolique. Pour  éviter ce problème, on peut considérer comme forme quasi-déployée de référence, 
celle qui correspond à la forme bilinéaire $\displaystyle \sum_{i=1}^{n+1}x_i^2-\sum_{i=n+1}^{2n+1}x_i^2$. Ce qui change alors est le calcul 
du caractère $S_{\underline{i}}$ du théorème principal. Exprimons sans démonstration quel caractère on trouve avec cet autre choix. 
Pour cela on a besoin de la notation pour tout $r\in [1,R]$, $a_{<r}=\sum_{j<r}a_j$  avec  $a_{<1}=0$. Alors $S_{\underline{i}}$ vaut 
sur $Z_r$ le caractère ${\sgn}^{\epsilon_r}$ où 
$$\epsilon_r=i_r(a_{<r}+1)+(a_r-i_r)a_{< r}+a_r(a_r+1)/2.$$

\section{Appendice}

Dans cet appendice, nous regroupons des résultats sur l'induction cohomologique et les représentations standard
dont nous avons besoin dans le corps de l'article. Nous suivons principalement \cite{AV}.

\subsection{Représentations standard}
Soient $\mathbf G$ un groupe algébrique connexe réductif défini sur $\bbR$ et $\tau$ l'involution de Cartan de $\mathbf G$.
Soit $\mathbf T$ un sous-groupe de Cartan de $\mathbf G$ défini sur $\bbR$ et $\tau$-stable, $T=\mathbf T(\bbR)$.
Soit $\mathbf B$ un sous-groupe de Borel de $\mathbf G$(seulement défini sur $\bbC$), contenant $\mathbf T$, 
 soit $R^+=R(\frt,\frb)$ le système de racines
positives défini par $\mathbf B$,  et soit $\rho$ la demi-somme des racines de $R^+$. Notons $E_\rho$ la représentation
spécifique de dimension $1$ de $T^\rho$ définie par $\rho$.

Comme dans \cite{AV}, Definition 8.11,  on introduit le revêtement à deux feuillets $T^\rho$
 (qui contrairement à ce que suggère la notation, ne dépend 
pas du  choix d'un système de racines positives,  {\sl cf.} \cite{AV}, Lemma 8.15).

\begin{defi}[\cite{AV}, Defiition 8.18]\label{bonneposition}
Considérons une paire de Borel  $(\mathbf B,\mathbf T)$ comme ci-dessus, et soit $\Lambda$ un caractère 
spécifique   de $T^\rho$ dont on note $\lambda$ la différentielle.
On dit que $\mathbf B$ et $\Lambda$ sont en bonne position relative si
 pour toute racine $\alpha\in R(\frt,\frb)=R^+$,  entière pour $\lambda$ (c'est-à-dire $\bil{\Lambda}{\check \alpha}\in \bbZ$), 
 l'une des  conditions suivantes est réalisée :
 
 $(i)$ $\alpha$ est imaginaire (c'est-à-dire $\tau \alpha=\alpha$) et $\bil{\lambda }{\check \alpha}\geq 0$,
 
$(ii)$ $\alpha$ est complexe (c'est-à-dire $\tau \alpha\neq \pm \alpha$), 
 $\tau \alpha \in R^+$, $\bil{\lambda }{\check \alpha}\geq 0$ ou  $\bil{\lambda }{\tau \check \alpha}\geq 0$,

$(iii)$ $\alpha$ est complexe,  $-\tau \alpha \in R^+$, $\bil{\lambda }{\check \alpha}\leq 0$ ou $\bil{\lambda }{-\tau \check \alpha}\leq 0$,

$(iv)$ $\alpha$ est réelle (c'est-à-dire $\tau \alpha=-\alpha$) et $\bil{\lambda }{\check \alpha}\leq 0$.

Si tel est le cas, on pose alors 
\begin{equation}\label{IBL}
I(\mathbf B,\Lambda)={}^\sharp \caR_\frb^d(\bbC_\Lambda).
\end{equation}

Ici, $d=d(G,T)=\dim \frn\cap \frk$ où $\frn$ est le radical nilpotent de $\frb$ et  le foncteur ${}^\sharp \caR_\frb^d$ 
est le foncteur d'induction cohomologique de Zuckerman en degré $d$, normalisé comme dans 
\cite{AV} ou \cite{KV}, c'est-à-dire que l'on tensorise $\bbC_\Lambda$ par le caractère $E_\rho$ de $T^\rho$ avant d'induire.
Le caractère infinitésimal est préservé.

On appelle $I(\mathbf B,\Lambda)$ une limite de représentation standard de $G$. Si de plus $\bil{\lambda }{\check \alpha}\neq 0$,
pour toute racine imaginaire $\alpha$, on dit que  $\Lambda$ est $M$-régulier et que $I(\mathbf B,\Lambda)$ une représentation standard  de $G$. 
\end{defi}

\begin{lemme}[\cite{AV}, Lemma 8.20] Si $(\mathbf B,\Lambda)$ et $(\mathbf B',\Lambda)$ sont en bonne position relative, 
et que les systèmes de racines positives imaginaires  et réelles définis par $\mathbf B$ et $\mathbf B'$ sont les mêmes,  alors 
$I(\mathbf B,\Lambda)=I(\mathbf B',\Lambda)$.
\end{lemme}

\begin{rmqs}\label{compat}1. Etant  donnés un système de racines imaginaires positives $R_{Im}^+$ et un système de racines réelles positives
$R^+_{Re}$ vérifiant  $\bil{\lambda }{\check \alpha}\geq 0$ pour toute racine $\alpha\in R_{Im}^+$ entière et 
$\bil{\lambda }{\check \alpha}\leq 0$ pour toute racine $\alpha\in R_{Re}^+$ entière, il existe
un sous-groupe de Borel $\mathbf B $ en bonne position relativement à $\Lambda$ tel que   
$R^+=R(\frt,\frb)$ contienne $R_{Im}^+$ et $R_{Re}^+$. 

--- 2. Etant donnés $R_{Im}^+$ et $R^+_{Re}$, le lemme montre que $I(\mathbf B,\Lambda)$ ne dépend pas du choix 
de $\mathbf B$ pourvu qu'il soit en bonne position relativement à $\Lambda$ et que $R_{Im}^+$ et $R^+_{Re}$
soient inclus dans $R(\frt,\frb)$. On peut donc noter 
\[I(\mathbf B,\Lambda)=I(R_{Im}^+,R^+_{Re},\Lambda) . \]
Si de plus $\Lambda$ est $M$-régulier, c'est-à-dire si $\bil{\check \alpha}{\lambda}\neq 0$ pour toute racine imaginaire $\alpha$, alors
$R_{Im}^+$ est déterminé par $\Lambda$, et l'on peut simplement noter 
\[I(\mathbf B,\Lambda)=I(R^+_{Re},\Lambda) .\]
\end{rmqs}

Nous allons maintenant rappeler le  lien avec la classification  de Langlands. 
\medskip 

{\bf Séries discrètes relatives.}
Commençons par le cas des séries discrètes relatives. On suppose donc que $T$ est compact modulo le centre.
On fixe un système de racines positives $\Sigma^+=\Sigma^+_{Im}$ de $\frt$ dans $\frg$.
 Soit $\rho$ la demi-somme des racines positives. Harish-Chandra
a montré que les séries discrètes relatives de $G$ sont paramétrées par les couples
$(\lambda,\chi)$ où $\lambda \in \frt^*$ et $\chi$ est un caractère de $Z(G)$ vérifiant deux conditions : 
\begin{equation}\label{regu1} \bil{\lambda}{\check \alpha}\neq 0, \quad (\alpha\in \Sigma^+),\end{equation}
\begin{equation}\label{regu2}  \text{Il existe un caractère $\Gamma$ de $T$ tel que
$\Gamma_{\vert Z(G)}=\chi$ et  $d\Gamma=\lambda-\rho$}.\end{equation}
Le caractère 
\[\Lambda=\Gamma\otimes \bbC_\rho\]
est alors un caractère spécifique de $T^\rho$ dont la différentielle est $\lambda$. Soit $\mathbf B$ le sous-groupe de 
Borel contenant $\mathbf T$ tel que $\lambda$ est dominant pour le système de racines positives défini par $\mathbf B$
(le système $\Sigma^+$ choisi au départ n'a finalement pas d'importance).

Remarquons que $\Lambda$ détermine $(\lambda,\chi)$ et que $\bbC^\rho$ étant une représentation
spécifique de $T^\rho$, $\Lambda$ est spécifique si et seulement si $\Gamma$ est un caractère de $T$.
La série discrète relative de paramètre d'Harish-Chandra $\Lambda$ est alors $I(\mathbf B,\Lambda)={}^\sharp \caR_\frb^S(\Lambda)$.

De même, les limites de séries discrètes relatives ${}^\sharp \caR_\frb^S(\Lambda)$ sont obtenues en abandonnant la condition (\ref{regu1}).
Le sous-groupe de Borel $\mathbf B$ n'est alors plus déterminé par $\lambda$. Il se peut que 
$I(\mathbf B,\Lambda)={}^\sharp \caR_\frb^S(\Lambda)$ soit nul. Ceci arrive si et seulement si 
l'une des racines simples $\alpha$ de $\mathbf T$ dans $\mathbf B$ telle  que $\bil{\lambda}{\check \alpha}= 0$  est imaginaire compacte.

Pour obtenir une paramétrisation des limites de séries discrètes relatives, on peut mettre la condition 
 $\bil{\lambda}{\check \alpha}\neq  0$ pour toute racine simple
imaginaire compacte  $\alpha$ de $\mathbf T$ dans $\mathbf B$, ou bien considérer par abus de 
langage que la représentation nulle est une limite de séries discrètes relative.
\bigskip 

{\bf Cas général.}
On revient à  $T,\Lambda$ comme au-début de la section.
On fixe aussi un système de racines imaginaires positives $R_{Im}^+$ et un système de racines réelles positives
$R^+_{Re}$ comme dans la remarque \ref{compat}.
 On note $\rho_{Im}$ la demi-somme des racines dans $R_{Im}^+$, et $T^{\rho_{Im}}$ le revêtement à deux feuillets de 
$T$ construit avec $\rho_{Im}$. Soit $M$ le sous-groupe de Levi de $G$ défini par $T$ et les racines imaginaires
($M=G^A$ si $T=T_cA$ est la décomposition de $T$ selon $\tau$).
Choisissons un sous-groupe parabolique $P=MN$ de $G$ vérifiant les trois conditions suivantes. 
Premièrement, $\frn$ est défini sur $\bbR$ (il est donc stable par $\sigma_G$, ou de manière équivalente,  stable par $-\tau$).
 Deuxièmement, $R_{Re}^+\subset  R(\frt,\frn)$. Troisièmement,
 pour toute racine $\alpha \in R(\frt,\frn)$ telle que 
$\bil{\lambda}{\check \alpha}\in \bbZ$, on a $\Re e(\bil{\check \alpha}{\tau\lambda-\lambda})\geq 0$.
Remarquons que cette condition  pour les racines réelles est  $\Re e(\bil{\check \alpha}{\lambda})\leq 0$, 
qui est déjà supposée satisfaite,  et pour les racines
complexes,  on est dans le cas $(iii)$ de la définition \ref{bonneposition}. Or si 
$\Re e(\bil{\check \alpha}{\tau\lambda-\lambda})\geq 0$, au moins l'un des entiers $\bil{\check \alpha}{\tau\lambda}$ ou
$-\bil{\check \alpha}{\lambda}$ est positif ou nul, ce qui montre que la condition $(iii)$ est vérifiée.
Un tel sous-groupe parabolique existe.
Soit alors $\mathbf B$ le sous-groupe de Borel contenu dans $\mathbf P$ et tel que $R^+_{Im}\subset R(\frt,\frb)$.
Il est alors clair que $\mathbf B$ et $\Lambda$ sont en bonne position relative. On peut obtenir 
$I(\mathbf B,\Lambda)$ par induction par étapes en passant  par $P$.
Soit $\rho_\frn$ la demi-somme des racines de $\frt$ dans $\frn$, et $T^{\rho_\frn}$ le revêtement à deux feuillets de 
$T$ construit avec $\rho_{\frn}$.
Le caractère $2\rho_\frn$ de $T$ est à valeurs réelles, et donc sa valeur absolue admet une unique 
racine carrée positive, que nous allons noter $\vert \rho_\frn \vert$. Son relèvement à $T^{\rho_\frn}$
coïncide sur la composante neutre avec le caractère $\rho_\frn$. Notons $\gamma_\frn$ le caractère
quotient $\frac{\rho_\frn}{\vert \rho_\frn\vert}$ de $T^{\rho_\frn}$ (c'est un caractère à valeurs dans les racines quatrième de l'unité).
On pose 
\[\Lambda_M =\Lambda\otimes  \bbC_{\gamma_\frn}.\]
On vérifie que c'est un caractère du revêtement $T^{\rho_{Im}}$ dont la différentielle est $\lambda$.
Soit $X_M$ la limite de série discrète relative de $M$ de paramètre $(\Lambda_M,\mathbf B_M=\mathbf B\cap \mathbf M)$
(ou encore de paramètre $(\lambda_M, R^+_{Im})$.
 On a alors  ({\sl cf.} \cite{AV}, Eq. (8.22))
\begin{equation}\label{IBLl}
I(\mathbf B,\Lambda)=\Ind_P^G(X_M).
\end{equation}

{\bf Autre réalisation des représentations standard.}
On repart de  $T,\Lambda, R_{Im}^+, R^+_{Re}$ comme ci-dessus.
Soit $\mathbf L$ le sous-groupe de Levi défini par $\mathbf T$ et les racines réelles. Il est défini sur $\bbR$, 
c'est un $c$-Levi au sens de la section \ref{cLevi}, 
et l'on note $L$ le groupe de ses points réels. Choisissons un sous-groupe parabolique 
$\mathbf Q=\mathbf L\mathbf V$ vérifiant les conditions suivantes. Premièrement, $\mathbf Q$ est $\tau$-stable.
Deuxièmement, $R_{Im}^+\subset  R(\frt,\frv)$. Troisièmement, pour chaque racine complexe $\alpha\in R(\frt,\frv)$ telle que 
$\bil{\check \alpha}{\lambda}\in \bbZ$, 
$\Re e(\bil{\check \alpha}{\tau\lambda+\lambda})\geq 0$. 
 Soit $\mathbf B$ le sous-groupe de Borel contenant $\mathbf T$, contenu dans 
$\mathbf Q$ et tel que $R^+_{Re}\subset R(\frt,\frb)$. Alors $\Lambda$ et $\mathbf B$ sont en bonne position relative.
Posons $\frb_L=\frb\cap \frl$; c'est une sous-algèbre de Borel de $\frl$. Notons $\frn_L$ son radical unipotent
et $\rho_{\frn_L}$ la demi-somme des racines de $\frt$ dans $\frn_L$.
On un un revêtement à deux feuillets $T^{\rho_{\frn_L}} $ et des caractères $\rho_{\frn_L}$, $\vert \rho_{\frn_L}\vert$ et 
$\gamma_{\frn_L}=\frac{\rho_{\frn_L}}{\vert \rho_{\frn_L}\vert}$ de ce revêtement.
Le caractère 
\[\Lambda_L=  \Lambda\otimes (\bbC_{\rho_\frv})^*\otimes \gamma_{\frn_L}\]
se factorise par $T$. Formons la série principale
\[X_L=\Ind_{B_L}^L(\Lambda_L) . \]
On a alors  ({\sl cf.} \cite{AV}, Eq. (8.23))
\[I(\mathbf B,\Lambda)=\caR_\frqqq^d (X_L) \]
où le foncteur d'induction cohomologique est ici celui normalisé comme dans \cite{Vgreen}, et $d=d(G,L)=\dim (\frv\cap \frk)$.

\subsection{Induites cohomologiques de représentations standard} \label{ICRS} 

On se place maintenant dans la situation suivante. 
Soit $\mathbf Q=\mathbf L\mathbf V$ un sous-groupe parabolique $\tau$-stable de $\mathbf G$, avec $\mathbf L$ défini sur 
$\bbR$. Soit $I(\mathbf B_L, \Lambda_L)$ une représentation standard de $L$. Ici $\Lambda_L$ est un caractère spécifique
d'un revêtement à deux feuillets $T^{\rho_L}$ d'un sous-groupe de Cartan $T$ de $L$ et $\mathbf B_L$ est un sous-groupe de 
Borel de $\mathbf L$ en bonne position relativement à $\Lambda_L$.
Posons $d=d(G,L)=\dim (\frv\cap \frk)$. On suppose que la différentielle $\lambda_L$ de $\Lambda_L$ est dans le  good range,
 et l'on forme
\begin{equation}  \label{RIBL}  \caR^d_\frqqq(I(\mathbf B_L,\Lambda_L)). \end{equation}

D'après la section précédente, on peut ecrire $I(\mathbf B_L,\Lambda_L)$ sous la forme
\begin{equation}  \label{RIBL2}  
I(\mathbf B_L,\Lambda_L)=\caR_{\frqqq_1}^{d_1} ( \Ind_{B_{L_1}}^{L_1}(\Lambda_{L_1})  )
\end{equation}
où $L_1$ est le sous-groupe de Levi de $L$ défini par $T$ et les racines réelles de $\frt$ dans $\frl$, 
$\mathbf Q_1=\mathbf L_1\mathbf V_1$ est un sous-groupe parabolique $\tau$-stable de $\mathbf L$
vérifiant $R_{L,Im}^+\subset  R(\frt,\frv_1)$ et   pour chaque racine complexe $\alpha\in R(\frt,\frv_1)$ telle que 
$\bil{\check \alpha}{\lambda_L}\in \bbZ$,  $\Re e(\bil{\check \alpha}{\tau\lambda_L+\lambda_L})\geq 0$.
D'autre part $d=d(L,L_1)=\dim(\frv_1\cap \frk)$.
  Le sous-groupe de Borel  $\mathbf B_{L}$  de $\mathbf L$ est celui contenant $\mathbf T$, contenu dans 
$\mathbf Q_1$ et tel que $R^+_{L,Re}\subset R(\frt,\frb_L)$. 
Posons $\frb_{L_1}=\frb_L\cap \frl_1$; c'est une sous-algèbre de Borel de $\frl_1$. Notons $\frn_1$ son radical unipotent
et $\rho_{\frn_1}$ la demi-somme des racines de $\frt$ dans $\frn_1$.
On un un revêtement à deux feuillets $T^{\rho_{\frn_1}} $ et des caractères $\rho_{\frn_1}$, $\vert \rho_{\frn_1}\vert$ et 
$\gamma_{\frn_1}=\frac{\rho_{\frn_1}}{\vert \rho_{\frn_1}\vert}$ de ce revêtement.
On a alors
\begin{equation}  \label{RIBL3}\Lambda_{L_1}=  \Lambda_L \otimes (\bbC_{\rho_{\frv_1}})^*\otimes \gamma_{\frn_1}.
\end{equation}
Compte tenu de (\ref{RIBL}) et (\ref{RIBL2}), on a
\begin{equation}  \label{RIBL4} 
 \caR^d_\frqqq(I(\mathbf B_L,\Lambda_L))=\caR^d_\frqqq \left( \caR_{\frqqq_1}^{d_1} \left( \Ind_{B_{L_1}}^{L_1} \left( 
 \Lambda_L\otimes (\bbC_{\rho_{\frv_1}})^*\otimes \gamma_{\frn_1}
 \right) \right) \right).
 \end{equation}
Posons $\mathbf Q'=\mathbf L_*\mathbf V_*\mathbf V=\mathbf L_*\mathbf V'$. C'est un sous-groupe parabolique $\tau$-stable
de $\mathbf G$ et par transitivité de l'induction cohomologique  ({\sl cf.} \cite{Vgreen}, Cor. 6.3.10), on obtient
\begin{equation}  \label{RIBL5} 
  \caR^d_\frqqq(I(\mathbf B_L,\Lambda_L))= \caR^{d+d_1}_{\frqqq'} \left( \Ind_{B_{L_1}}^{L_1} \left( 
 \Lambda_L\otimes (\bbC_{\rho_{\frv_1}})^*\otimes \gamma_{\frn_1}
 \right) \right) .
 \end{equation}
Posons 
\begin{equation}  \label{RIBL6} 
 \Lambda=\Lambda_L\otimes \bbC_{\rho_\frv}
\end{equation}
C'est un caractère du revêtement $T^\rho$ de $T$.
On a $\Lambda_L=\Lambda\otimes (\bbC_{\rho_\frv})^*$, et donc 
\begin{equation}  \label{RIBL7} 
\caR^d_\frqqq(I(\mathbf B_L,\Lambda_L))= \caR^{d+d_1}_{\frqqq'} \left( \Ind_{B_{L_1}}^{L_1} \left( 
 \Lambda \otimes (\bbC_{\rho_{\frv'}})^*\otimes \gamma_{\frn_1}
 \right) \right) .
\end{equation}
Soit $\mathbf B$ le sous-groupe de Borel de $\mathbf G$ contenant $\mathbf T$, contenu dans $\mathbf Q'$
et  tel que $R^+_{L,Re}=R^+_{Re}\subset R(\frt,\frb)$ (c'est-à-dire $\mathbf B=\mathbf B_L\mathbf V$). 
On vérifie facilement que $\mathbf B$ et $\Lambda$ sont en bonne position relative, et l'on a finalement
\begin{equation}  \label{RIBL8}  \caR^d_\frqqq(I(\mathbf B_L,\Lambda_L))=
I(\mathbf B,\Lambda). \end{equation}

Réécrivons maintenant ceci en termes  de  représentations standard  de Langlands. On commence par 
$I(\mathbf B_L,\Lambda_L)$ que l'on écrit de manière analogue à (\ref{IBLl})
\[I(\mathbf B_L,\Lambda_L)=\Ind_{P_L=M_LN_L}^L  (X_{M_L}). \]
Ici, $M_L$ est le sous-groupe de Levi de $L$ défini par $T$ et les racines imaginaires de $\frt$ dans $\frl$, 
$P_L=M_LN_L$ est un sous-groupe parabolique vérifiant les conditions voulues de positivité, $X_{M_L}$ est la série 
discrète relative de $M_L$ de paramètre d'Harish-Chandra
\[ \Lambda_{M_L}= \Lambda_L \otimes \bbC_{\gamma_{\frn_L}}. \]
On obtient alors, avec (\ref{IBLl}):
\begin{prop} \label{TorsionI}Avec les notations ci-dessus :
\begin{equation}  \label{RIBL9}
 \caR^{S}_{\frqqq}\left(\Ind_{P_L=M_LN_L}^L  (X_{M_L})\right) =  \Ind_{P=MN}^G(X_M).
\end{equation}
où $X_M$ est la série discrète relative de $M$ de paramètre d'Harish-Chandra
\begin{equation}  \label{RIBL10}\Lambda_M=\Lambda \otimes \bbC_{\gamma_{\frn}}
 =  \Lambda_L\otimes \bbC_{\rho_\frv} \otimes \bbC_{\gamma_{\frn}}
= \Lambda_{M_L}\otimes (\bbC_{\gamma_{\frn_L}})^* \otimes \bbC_{\rho_\frv} \otimes \bbC_{\gamma_{\frn}} .
\end{equation}
Ainsi, dans le \og good range\fg, l'induction cohomologique envoie une limite de représentation standard sur une limite de représentation
standard, et  la formule \ref{RIBL10} donne le paramètre d'Harish-Chandra de la représentation induite en fonction de celui de la
représentation induisante.
\end{prop}
\bigskip

\subsection{Calcul du terme $(\bbC_{\gamma_{\frn_L}})^* \otimes \bbC_{\rho_\fru} \otimes \bbC_{\gamma_{\frn}} $ pour les groupes classiques}
\label{calculTorsionI}
Nous nous replaçons dans le contexte de la section précédente, en particulier de la proposition \ref{TorsionI}, que nous particularisons au cas
 des groupes classiques de la section \ref{grcla}, avec pour sous-groupe parabolique $\mathbf Q$ l'un de ceux de la section \ref{cLevicla}.
 En particulier, $\mathbf L$ est un $c$-Levi maximal  de $\mathbf G$, et l'on a 
 \[L =\U(i,c-i)\times G_i' \]
 où $G'_i$ est une forme intérieure pure d'un groupe $\mathbf G'$ de même type que $\mathbf G$.

 Le caractère $\bbC_{\rho_V}$  est sans mystère. Tout d'abord $\rho_V=\rho_G-\rho_L$, et l'on calcule en coordonnées usuelles,
  dans chacun des cas : 
  
  {\bf cas A} : $\rho_V=\rho_G-\rho_L=(n,n-1,\ldots ,1)-\left( \frac{c-1}{2}, \frac{c-3}{2}, \ldots ,-  \frac{c-1}{2}, n-c, \ldots ,1\right)=
  \left( n- \frac{c-1}{2}, \ldots ,  n- \frac{c-1}{2}, 0, \ldots 0\right)$.
 
  {\bf cas B} : $\rho_V=\rho_G-\rho_L=(\frac{2n-1}{2},\frac{2n-3}{2},\ldots ,\frac{1}{2})-\left( \frac{c-1}{2}, \frac{c-3}{2}, \ldots ,-  \frac{c-1}{2},
  \frac{2(n-c)-1}{2}, n-c, \ldots ,\frac{1}{2}\right)=
  \left(  \frac{2n-c}{2}, \ldots ,   \frac{2n-c}{2}, 0, \ldots 0\right)$.
 
  {\bf cas C et D} : $\rho_V=\rho_G-\rho_L=(n-1,n-2,\ldots ,0)-\left( \frac{c-1}{2}, \frac{c-3}{2}, \ldots ,-  \frac{c-1}{2}, n-c-1, \ldots ,0\right)=
  \left( n-1- \frac{c-1}{2}, \ldots ,  n-1- \frac{c-1}{2}, 0, \ldots 0\right)$.

 On constate que $\rho_\frv$ ne s'intègre pas nécessairement en un caractère de $L$ (il y a une condition sur la parité de $c$ dans chaque cas
 pour que les coordonnées de $\rho_V$ soit des entiers), mais qu'il  définit bien un 
 caractère central d'un revêtement double de $L$, trivial sur le facteur $G'_i$.
Comme les sous-groupes de Cartan des groupes unitaires sont connexes, $\rho_V$ détermine $\bbC_{\rho_V}$.

Le terme $(\bbC_{\gamma_{\frn_L}})^*  \otimes \bbC_{\gamma_{\frn}} $ est plus difficile à calculer.
Traitons d'abord  le cas des groupes symplectiques. Rappelons que nous sommes parti  d'un sous-groupe de Cartan 
 de $\mathbf G$ que nous appelons ici $\mathbf D$ plutôt que $\mathbf T$ comme dans la section précédente. 
 A conjugaison près, on peut supposer $D$ de la forme 
\[D\simeq  \U(1)^{r_1} \times (\bbC^\times)^{m_1}\times  \U(1)^{r_2} \times (\bbC^\times)^{m_2}\times (\bbR^\times)^{s_2}, 
\]
avec $c=r_1+2m_1$ et $r_1+2m_1+r_2+2m_2+s_2=n$.

La sous-algèbre parabolique $\frn$ est la somme  d'espaces radiciels $\frg^\alpha_\bbR$, avec $\alpha$ racine réelle, et 
d'espaces $\left(\frg^\alpha \oplus\frg^{\sigma(\alpha)}\right)^\sigma$, où $(\alpha,\sigma(\alpha))$ sont des couples de racines complexes.
Comme toutes les racines réelles de $\mathbf D$ dans $\mathbf G$ sont aussi dans $\mathbf L$, un supplémentaire
de $\frn_L$ dans $\frn$ est constitué de la somme d'espaces  
$\left(\frg^\alpha \oplus\frg^{\sigma(\alpha)}\right)^\sigma$, où $\alpha,\sigma(\alpha)$ est un  couple de racines complexes qui n'est pas dans $L$.
Les racines de  $\mathbf D$ dans $\mathbf G$ qui ne sont pas dans  $\mathbf L$ sont les 
$\pm (e_i+e_j)$, $1\leq i<j\leq c$, $\pm(e_i\pm e_j)$, $1\leq i\leq c<j\leq n$, $\pm 2e_i$, $c<i\leq n$.

On a 
\[\sigma(e_i)=-\epsilon_i \text{ si }  1\leq i\leq r_1, \,  c+1\leq i\leq c+r_2, \]
\[\sigma(e_i)=e_i \text{ si }    c+r_2+2m_2+1\leq i\leq  n, \]
et on peut supposer, quitte à conjuguer, que 
\begin{align*}&\sigma(e_{r_1+2i-1})=-e_{r_1+2i} \text{ si }  1\leq i\leq m_1, \\ 
&\sigma(e_{c+r_2+2i-1})=-e_{c+r_2+2i} \text{ si }  1\leq i\leq m_2.
\end{align*}
Fixons $i\leq r_1$.  On regarde les couples de racines complexes $(\alpha,\sigma(\alpha))$ dans $\frn$ qui ne sont pas dans 
$\frn_L$ pour calculer leur contribution à $\gamma_\frn \gamma_{\frn_L}^*$, au besoin
en groupant plusieurs couples.  
\[(e_i+e_{r_1+2j-1} ,- e_i-e_{r_1+2j}), \; 1\leq j\leq m_1, \; \frac{\alpha+\sigma(\alpha)}{2}=
\frac{e_{r_1+2j-1}-e_{r_1+2j}}{2},\]
\[(e_i+e_{r_1+2j} ,- e_i-e_{r_1+2j-1}), \; 1\leq j\leq m_1, \; \frac{\alpha+\sigma(\alpha)}{2}=
\frac{-e_{r_1+2j-1}+e_{r_1+2j}}{2},\]
La contribution de ces deux couples va donc être triviale.

\[(e_i+e_{c+r_2+2j-1} ,- e_i-e_{c+r_2+2j}), \; 1\leq j\leq m_2, \; \frac{\alpha+\sigma(\alpha)}{2}=
\frac{e_{c+r_2+2j-1}-e_{c+r_2+2j}}{2},\]
\[(e_i-e_{c+r_2+2j-1} , -e_i+e_{c+r_2+2j}), \; 1\leq j\leq m_2, \; \frac{\alpha+\sigma(\alpha)}{2}=
\frac{-e_{c+r_2+2j-1}+e_{c+r_2+2j}}{2},\]
 ou bien 
\[(-e_i+e_{c+r_2+2j-1} , e_i-e_{c+r_2+2j}), \; 1\leq j\leq m_2, \; \frac{\alpha+\sigma(\alpha)}{2}=
\frac{e_{c+r_2+2j-1}-e_{c+r_2+2j}}{2},\]    
La contribution de ces deux couples va donc être triviale, ou bien $\frac{e_{c+r_2+2j-1}-e_{c+r_2+2j}  }
{\vert e_{c+r_2+2j-1}-e_{c+r_2+2j}\vert}$.   Mais comme la racine réelle $e_{c+r_2+2j-1}-e_{c+r_2+2j}$
vit sur un facteur $\bbC^\times$ qui est connexe, c'est de toutes façon trivial.

 \[(e_i+e_{j} , -e_i+e_{j}), \; c+r_2+2m_2+1\leq j\leq n, \; \frac{\alpha+\sigma(\alpha)}{2}=
e_j,\]
La contribution d'un tel  couple à $\gamma_\frn\gamma_{\frn_L}^*$ est donc $\frac{e_j}{\vert e_j\vert}$.

\bigskip

Fixons $1\leq i\leq m_1$. Les indices  $r_1+2i-1$  et  $r_1+2i$ vont  contribuer de la manière suivante à $\rho_\frn-\rho_{\frn_L}$
(sans recompter les contributions étant déjà apparues ci-dessus) :  
\[(2e_{r_1+2i-1} , -2e_{r_1+2i}), \; 1\leq i \leq m_1, \; \frac{\alpha+\sigma(\alpha)}{2}=
e_{r_1+2i-1}-e_{r_1+2i},\]
La contribution est donc  $\frac{e_{c+r_2+2j-1}-e_{c+r_2+2j}  }{\vert e_{c+r_2+2j-1}-e_{c+r_2+2j}\vert}$.   
Mais comme la racine réelle $e_{c+r_2+2j-1}-e_{c+r_2+2j}$
vit sur un facteur $\bbC^\times$ qui est connexe, c'est trivial.

\[(e_{r_1+2i-1} + e_{r_1+2i'-1},-e_{r_1+2i} -e_{r_1+2i'}), \; 1\leq i<i' \leq m_1, \]
\[ \frac{\alpha+\sigma(\alpha)}{2}=
\frac{e_{r_1+2i-1}-e_{r_1+2i}+e_{r_1+2i'-1}-e_{r_1+2i'}}{2},\]

\[(e_{r_1+2i-1} + e_{r_1+2i'},-e_{r_1+2i} - e_{r_1+2i'-1}), \; 1\leq i<i' \leq m_1, \]
\[ \frac{\alpha+\sigma(\alpha)}{2}=
\frac{e_{r_1+2i-1}-e_{r_1+2i}-e_{r_1+2i'-1}+e_{r_1+2i'}}{2},\]
 La contribution de ces deux couples à $\gamma_\frn\gamma_{\frn_L}^*$ est donc 
 $\frac{e_{r_1+2i-1}-e_{r_1+2i}}{\vert e_{r_1+2i-1}-e_{r_1+2i}}$.
 Mais comme la racine réelle $ e_{r_1+2i-1}-e_{r_1+2i} $
vit sur un facteur $\bbC^\times$ qui est connexe, c'est trivial.
D'autre choix de racines positives peuvent conduire à une contribution 
$\frac{e_{r_1+2i'-1}-e_{r_1+2i'}}{\vert e_{r_1+2i'-1}-e_{r_1+2i'}}$, mais on a la même conclusion.

\[(e_{r_1+2i-1} + e_j,-e_{r_1+2i} -e_j, \; c+1\leq j \leq c+r_2, \]
\[ \frac{\alpha+\sigma(\alpha)}{2}=
\frac{e_{r_1+2i-1}-e_{r_1+2i}}{2},\]

\[(e_{r_1+2i-1} - e_j,-e_{r_1+2i} +e_j, \; c+1\leq j \leq c+r_2, \]
\[ \frac{\alpha+\sigma(\alpha)}{2}=
\frac{e_{r_1+2i-1}-e_{r_1+2i}}{2},\]

La contribution de ces deux couples à $\gamma_\frn\gamma_{\frn_L}^*$ est donc 
 $\frac{e_{r_1+2i-1}-e_{r_1+2i}}{\vert e_{r_1+2i-1}-e_{r_1+2i}}$.
 Mais comme la racine réelle $ e_{r_1+2i-1}-e_{r_1+2i} $
vit sur un facteur $\bbC^\times$ qui est connexe, c'est trivial.

\[(e_{r_1+2i-1} + e_{c+r_2+2'-1},-e_{r_1+2i} -e_{c+r_2+2i'}), \; 1\leq i<i' \leq m_2, \]
\[ \frac{\alpha+\sigma(\alpha)}{2}=
\frac{e_{r_1+2i-1}-e_{r_1+2i}+e_{r_1+2i'-1}-e_{r_1+2i'}}{2},\]

\[(e_{r_1+2i-1} + e_{c+r_2+2i'},-e_{r_1+2i} - e_{c+r_2+2i'-1}), \; 1\leq i<i' \leq m_2, \]
\[ \frac{\alpha+\sigma(\alpha)}{2}=
\frac{e_{r_1+2i-1}-e_{r_1+2i}-e_{c+r_2+2i'-1}+e_{c+r_2+2i'}}{2},\]
 La contribution de ces deux couples à $\gamma_\frn\gamma_{\frn_L}^*$ est donc 
 $\frac{e_{r_1+2i-1}-e_{r_1+2i}}{\vert e_{r_1+2i-1}-e_{r_1+2i}}$.
 Mais comme la racine réelle $ e_{r_1+2i-1}-e_{r_1+2i} $
vit sur un facteur $\bbC^\times$ qui est connexe, c'est trivial.
D'autre choix de racines positives peuvent conduire à une contribution 
$\frac{e_{c+r_2+2i'-1}-e_{c+r_2+2i'}}{\vert e_{c+r_2+2i'-1}-e_{c+r_2+2i'}}$, mais on a la même conclusion.

\[(e_{r_1+2i-1} + e_j,-e_{r_1+2i} +e_j, \; c+r_2+2m_2+1\leq j \leq n, \]
\[ \frac{\alpha+\sigma(\alpha)}{2}=
\frac{e_{r_1+2i-1}-e_{r_1+2i}}{2}+e_j,\]

\[(e_{r_1+2i-1} - e_j,-e_{r_1+2i} -e_j, \; c++r_2+2m_2+1\leq j \leq n, \]
\[ \frac{\alpha+\sigma(\alpha)}{2}=
\frac{e_{r_1+2i-1}-e_{r_1+2i}}{2}-e_j,\]

ou bien
\[(-e_{r_1+2i-1} + e_j,e_{r_1+2i} +e_j, \; c++r_2+2m_2+1\leq j \leq n, \]
\[ \frac{\alpha+\sigma(\alpha)}{2}=
\frac{-e_{r_1+2i-1}+e_{r_1+2i}}{2}+e_j,\] 

La contribution de ces deux couples à $\gamma_\frn\gamma_{\frn_L}^*$ est donc 
 $\frac{e_{r_1+2i-1}-e_{r_1+2i}}{\vert e_{r_1+2i-1}-e_{r_1+2i}}$.
 Mais comme la racine réelle $ e_{r_1+2i-1}-e_{r_1+2i} $
vit sur un facteur $\bbC^\times$ qui est connexe, c'est trivial.
 ou bien $\frac{2e_j}{\vert 2 e_j\vert }$. Mais là encore, la contribution est triviale à cause du $2$.

\bigskip 

Finalement, en regardant toutes ces contributions, on voit que $\gamma_\frn\gamma_{\frn_L}^*$
est le caractère quadratique de $D$ qui met sur le $j$-ième facteur $\bbR^\times$ le
signe $\left(\frac{e_j}{\vert e_j\vert}\right)^{r_1}= \left(\frac{e_j}{\vert e_j\vert}\right)^{r_1+2m_1}=\left(\frac{e_j}{\vert e_j\vert}\right)^{c} $. 

Si $c$ est pair, la torsion $\gamma_\frn\gamma_{\frn_L}^*$ est donc triviale. Si $c$ est impair, 
notons $\sgn_D$ ce caractère quadratique non trivial de $D$ qui met le caractère $\sgn$ sur chacun des facteurs $\bbR^\times$.

\bigskip 

Le calcul pour les groupes orthogonaux est similaire. Pour le cas {\bf B} des groupes orthogonaux impairs, 
on remplace les racines longues $2e_i$ dans le calculs ci-dessus par les racines courtes $e_i$, et dans les cas 
 {\bf C} et  {\bf D}  des groupes orthogonaux impairs, on enlève ces racines longues du calcul.  Comme elles ne contribuaient pas 
et que dans le cas {\bf B}, les racines courtes $e_i$ ne vont pas contribuer non plus, on obtient exactement le même résultat.

\bigskip 

\bigskip
\bibliographystyle{smfalpha}
\bibliography{MR3}

  \end{document}